\documentclass[leqno, notitlepage]{amsart}

\usepackage{amssymb}
\usepackage[all, poly]{xy}

\DeclareMathOperator{\Hom}{Hom}

\DeclareMathOperator{\End}{End}
\DeclareMathOperator{\Aut}{Aut}

\DeclareMathOperator{\In}{In}
\DeclareMathOperator{\Out}{Out}
\DeclareMathOperator{\tr}{tr}
\DeclareMathOperator{\Coker}{Coker}
\DeclareMathOperator{\rank}{rank}
\DeclareMathOperator{\im}{im}
\DeclareMathOperator{\Id}{Id}

\newcommand{\sss}{\scriptscriptstyle}

\swapnumbers
\theoremstyle{plain} 
\newtheorem*{theorem}{Theorem}
\newtheorem*{corollary}{Corollary}

\newtheorem*{proposition}{Proposition}
\newtheorem{numtheorem}[subsection]{Theorem}

\theoremstyle{remark}
\newtheorem*{remark}{Remark}

\numberwithin{equation}{subsection}

\numberwithin{enumi}{subsection}
\newenvironment{alphenum}{
\begin{enumerate}
\setcounter{enumi}{\value{equation}}
}{
\setcounter{equation}{\value{enumi}}
\end{enumerate}}

\begin{document}

\title[Tensor product varieties]{
Tensor product varieties and crystals \\
ADE case}

\author{ Anton Malkin }
\address{ Department of Mathematics, Yale University,
P.O. Box 208283, New Haven, CT 06520-8283 }
\email{anton.malkin@yale.edu}


\begin{abstract}

Let $\mathfrak{g}$ be a simple simply laced Lie
algebra. 
In this paper two families of varieties associated
to the Dynkin graph of $\mathfrak{g}$
are described:
``tensor product'' and ``multiplicity''
varieties. These varieties are closely
related to Nakajima's quiver varieties 
and should play an
important role in the geometric constructions of 
tensor products and intertwining operators.
In particular it is shown that
the set of irreducible components of a tensor product
variety can be equipped with a structure of 
$\mathfrak{g}$-crystal isomorphic to the crystal of
the canonical basis of the tensor
product of several simple finite dimensional
representations of $\mathfrak{g}$, and that
the number of irreducible components of a multiplicity
variety is equal to the multiplicity of a certain 
representation in the tensor product of several
others. Moreover the decomposition of
a tensor product into a direct sum is
described geometrically (on the level
of crystals). 

\end{abstract}

\maketitle

\setcounter{section}{-1}

\section{Introduction}

\subsubsection{}

The purpose of this paper is to define
and study two (closely related) families
of quasi-projective varieties associated to
a simply laced Dynkin graph $D$:
the \emph{tensor product varieties} and 
the \emph{multiplicity varieties}.

\subsubsection{}

In this paper for the sake of clarity
the Dynkin graph $D$
is assumed to be of ADE type. This
condition can be relaxed by adjusting
slightly the definitions of the varieties
involved. 

Let $\mathfrak{g}$ be the simple Lie
algebra associated to $D$ and 
$\mathfrak{g}' = \mathfrak{g}
\oplus \mathfrak{t}$, where
$\mathfrak{t}$ is a Cartan
subalgebra of $\mathfrak{g}$.
Thus $\mathfrak{g}'$ is a reductive
Lie algebra.

The precise definition of the
tensor product and the multiplicity
varieties is rather involved 
(see \ref{DefinitionOfTS}). For the
purpose of this Introduction it is
enough to know that a multiplicity
variety 
${^n \mathfrak{S}} 
(\mu^0; \mu^1, \ldots , \mu^n)$
is associated to a set
$\mu^0, \mu^1, \ldots , \mu^n$, where
$\mu^k$ belongs to a certain
subset of the weight lattice
of the reductive Lie algebra
$\mathfrak{g}'$ (the set of ``integrable
positive'' weights).
A tensor product variety
${^n \mathfrak{T}} 
(\mu^1, \ldots , \mu^n; \nu)$
is associated to a set
$\mu^1, \ldots , \mu^n, \nu$, where 
$\mu^1, \ldots , \mu^n$ are
as above, and
$\nu$ is a weight of $\mathfrak{g}'$.
This way of writing parameters
of ${^n \mathfrak{S}}$ and 
${^n \mathfrak{T}}$ is for the purposes
of the Introduction only.
Notation in the main body
of the paper is different.

In the case $n =1$ the variety
${^1 \mathfrak{S}} (\mu^0; \mu^1)$ is
empty unless $\mu^0 = \mu^1$ in
which case it coincides with
Nakajima's quiver variety
$\mathfrak{M}^{reg}_0 (v,w)$ 
(cf. \cite{Nakajima1994, Nakajima1998}).
Here $v$ and $w$ can be expressed 
in terms of $\mu^0$.
Similarly
${^1 \mathfrak{T}} (\mu^1; \nu)$
coincides with
a fiber $\mathfrak{M} (\mu^1, \nu)=
\mathfrak{M} (v, v_0, w)$
of Nakajima's resolution of
singularities of the singular 
quiver variety (where again
$v$, $v_0$, and $w$ can be expressed 
in terms of $\mu^1$ and $\nu$).

\subsubsection{}
\label{IntroNakajimaCrystals}

Multiplicity varieties and 
tensor product varieties 
have pure dimensions. Let 
${^n \mathcal{T}} 
(\mu^1, \ldots , \mu^n ; \nu)$
(resp. ${^n \mathcal{S}} 
(\mu^0; \mu^1, \ldots , \mu^n )$,
$\mathcal{M} (\mu, \nu)$) be
the set of irreducible components
of the tensor product variety
${^n \mathfrak{T}} 
(\mu^1, \ldots , \mu^n ; \nu)$
(resp. the multiplicity variety
${^n \mathfrak{S}} 
(\mu^0; \mu^1, \ldots , \mu^n )$,
Nakajima's variety
$\mathfrak{M} (\mu, \nu)$), and
let
${^n \mathcal{T}} 
(\mu^1, \ldots , \mu^n ) =
\bigsqcup_{\nu} {^n \mathcal{T}} 
(\mu^1, \ldots , \mu^n ; \nu) \; $,
$\mathcal{M} (\mu ) =
\bigsqcup_{\nu} \mathcal{M} (\mu , \nu)$.

Nakajima \cite{Nakajima1998}
(based on an idea of Lusztig 
\cite[12]{Lusztig1991a}) introduced
a structure of a $\mathfrak{g}'$-crystal
on the set $\mathcal{M} (\mu )$
and it was shown by Saito 
(based on his joint work
with Kashiwara \cite{KashiwaraSaito})
that this crystal is isomorphic to the
crystal of the canonical basis of
the irreducible representation
of $\mathfrak{g}'$ with highest 
weight $\mu$. Strictly speaking the 
above mentioned authors consider
$\mathfrak{g}$-crystals, but the extension
to $\mathfrak{g}'$ is straightforward and
the weight lattice of 
$\mathfrak{g}'$ appears more naturally in
geometry of quiver varieties.

The main result of this paper is the
construction of two bijections between sets
of irreducible components 
(cf. \ref{DefinitionOfAlpha},
\ref{SecondBijection}):
\begin{equation}\nonumber
\begin{split}
\alpha_n : \quad {^n \mathcal{T}} 
(\mu^1, \ldots , \mu^n )
&\xrightarrow{\sim}
\mathcal{M} ( \mu_1 ) 
\times \ldots \times
\mathcal{M} ( \mu_n ) \; ,
\\
\beta_n : \quad {^n \mathcal{T}} 
(\mu^1, \ldots , \mu^n )
&\xrightarrow{\sim}
\bigsqcup_{\mu_0}
{^n \mathcal{S}} 
(\mu^0; \mu^1, \ldots , \mu^n )
\times
\mathcal{M} ( \mu^0 ) \; .
\end{split}
\end{equation}
Moreover one has the following theorem
(cf. Theorems \ref{MainTheoremForQ},
\ref{GPrimeSection}).
\begin{theorem}\label{IntroTheorem}
The composite bijection
\begin{equation}\nonumber
\tau_n = 
\beta_n \circ \alpha_n^{-1} : \quad 
\mathcal{M} ( \mu_1 ) 
\otimes \ldots \otimes
\mathcal{M} ( \mu_n )
\xrightarrow{\sim}
\bigoplus_{\mu_0}
{^n \mathcal{S}} 
(\mu^0; \mu^1, \ldots , \mu^n )
\otimes
\mathcal{M} ( \mu^0 ) 
\end{equation}
is an isomorphism of 
$\mathfrak{g}'$-crystals if
one endows the sets
$\mathcal{M} ( \mu_1 ), \ldots , 
\mathcal{M} ( \mu_n )$ with Nakajima's
crystal structure, and considers
the set ${^n \mathcal{S}} 
(\mu^0; \mu^1, \ldots , \mu^n )$ as
a trivial crystal.
\end{theorem}
In the above theorem $\otimes$ 
(resp. $\oplus$) denotes
the tensor product (resp. direct sum)
of crystals, which
coincides with the direct product
(resp. disjoint union)
on the level of sets.

It follows in particular, that the set 
${^n \mathcal{T}} 
(\mu^1, \ldots , \mu^n )$ can be equipped
with the structure of a 
$\mathfrak{g}'$-crystal using ether
of the bijections $\alpha_n$ or $\beta_n$,
and this crystal is isomorphic to
the crystal of the canonical basis
in the tensor product of $n$ irreducible
representations of $\mathfrak{g}'$ with
highest weights $\mu^1, \ldots , \mu^n$.

\subsubsection{}

The proof of the Theorem 
\ref{IntroTheorem} uses double reduction.
First it is shown that the statement 
for any $n$ follows from the 
corresponding statement for $n=2$.
Then one uses restriction to
Levi factors of parabolic
subalgebras of $\mathfrak{g}$ to
reduce the problem to $\mathfrak{sl}_2$
case. When $n=2$ and 
$\mathfrak{g}=\mathfrak{sl}_2$ the 
Theorem becomes an elementary linear
algebraic statement.

\subsubsection{}
The existence of the crystal
isomorphism $\tau_n$ allows one
to apply a theorem of Joseph
\cite[Proposition 6.4.21]{Joseph1995}
about the uniqueness of the family of
crystals closed with respect to tensor
products to give another proof of
isomorphism between $\mathcal{M} (\mu)$
and the crystal of the canonical basis
of the irreducible
representation of $\mathfrak{g}'$
with highest weight $\mu$, and,
more importantly, to prove that
the set ${^n \mathcal{S}} 
(\mu^0; \mu^1, \ldots , \mu^n )$ 
of irreducible components of
a multiplicity variety has the cardinal equal
to the multiplicity of the irreducible
representation of $\mathfrak{g}'$ with
highest weight $\mu^0$ in a direct sum
decomposition of the tensor product
of $n$ irreducible representations with 
highest weights $\mu^1 , \ldots , \mu^n$.
It is also equal to the corresponding
multiplicity for the tensor product of
representations of the simple Lie algebra 
$\mathfrak{g}$ if one restricts the 
weights to the Cartan subalgebra of
$\mathfrak{g}$. 

In other words one has
the following corollary of the Theorem
\ref{IntroTheorem}.

\begin{corollary}\label{IntroCorollary}
Let $L ( \mu )$ denote the irreducible 
representation of $\mathfrak{g}'$ with
highest weight $\mu$.
Then
\begin{equation}\nonumber
| {^n \mathcal{S}} (
\mu^0; \mu^1, \ldots , \mu^n ) | =
\dim_{\mathbb{C}} \Hom_{\mathfrak{g}'}
\bigl( L (\mu^0), L (\mu^1) 
\otimes \ldots
\otimes L (\mu^n) \bigr) \; .
\end{equation}
\end{corollary}

\subsubsection{}

The paper is organized as follows. Section 
\ref{CrystalsSection} contains a review
of Kashiwara's theory of crystals. Section
\ref{QuiverTensorSection} begins
with a description of Nakajima's
quiver varieties. Then the
tensor product and the multiplicity varieties
are defined in \ref{DefinitionOfTS}. The
reminder of Section 
\ref{QuiverTensorSection} is devoted to
various properties of these varieties.
In particular the tensor decomposition
bijection is described in
\ref{TensorDecomposition}.
In Section \ref{QuiverCrystals} it is
shown (following Lusztig and Nakajima)
how one can use restrictions to Levi
subalgebras of parabolic subalgebras 
of $\mathfrak{g}$ to 
define crystal structures
on the sets of irreducible components 
of the varieties involved. Then
a variant of 
Theorem \ref{IntroTheorem}
(Theorem \ref{MainTheoremForQ}) is proven
using the Levi restriction.

\subsubsection{}
\label{IntroAppendix}

Definitions of various varieties
given in the main body of the paper are
rather messy due to abundance of 
indices. Unfortunately, one needs 
this notation in order to follow the
arguments used in the proofs. However
in the Appendix a more conceptual
definition of the multiplicity varieties
is given. Namely they appear to be
closely related to the Hall-Ringel 
algebra (cf. \cite{Ringel1988})
associated to a certain algebra
$\Tilde{\mathcal{F}}$ introduced by Lusztig
\cite{Lusztig1998, Lusztig2000a}.
The Appendix also contains a diagram
of varieties called the tensor product 
diagram which can be (conjecturally) used
to equip a certain category of
perverse sheaves on a variety $Z_D$
described by Lusztig 
\cite{Lusztig1998} with a structure
of a Tannakian category. 

\subsubsection{}

Tensor product varieties ${^n \mathfrak{T}} ( \mu^1 , \ldots , \mu^n ;
v)$ are also described in a recent paper of H. Nakajima 
\cite{Nakajima2001p}.
Actually, Nakajima only considers some special values of the highest
weights $\{ \mu^i \}$ that correspond to Lagrangian fibers in quiver
varieties. However it is not restrictive for representation theory
applications, and generalization to arbitrary $\{ \mu^i \}$ is
straightforward.

Nakajima studies deeper geometric structures on a
tensor product variety than just the set of its irreducible components,
namely the Borel-Moore homology and equivariant $K$-theory, which allows
him to obtain some very interesting results in the representation theory
of quantum affine algebras at zero level.

On the contrary the author of the
present paper is mainly interested in the geometric description of the
multiplicities in tensor products, and more generally in monoidal
structures on various geometric objects on the singular quiver
variety $Z$ (cf. the Appendix). The importance of tensor product varieties
for him is in their relation to resolutions of singularities of
multiplicity varieties.

\subsubsection{}
Various special cases of tensor product
and multiplicity varieties have appeared in
the literature before (and certainly served
as important sources of motivation for
the author). The most important case is
related to a Dynkin graph of type $A$.
It is known (cf. \cite{Nakajima1994})
that some special cases of quiver varieties
associated to this Dynkin graph are related
to partial resolution of singularities of
the nilpotent cone in $\mathfrak{gl}_{k}$.
Similarly, the multiplicity varieties in 
this case are related to Spaltenstein 
varieties for $\mathfrak{gl}_{k}$
(a Spaltenstein variety is a variety
consisting of all parabolic subgroups $P$
of a given type that contain a given
nilpotent operator $t \in \mathfrak{gl}_{k}$, 
and such that the projection of
$t$ to the Levi factor $L$ of $P$ belongs to
a given nilpotent orbit in $L$). A theorem due
to Hall (cf. \cite{Hall1959}, 
\cite[Chapter II]{Macdonald1995})
implies that the number of irreducible
components of a Spaltenstein variety is equal
to a certain multiplicity in the tensor product
decomposition for $\mathfrak{gl}_{\sss N}$ 
($k$ has no relation to $N$), which
is a special case of Corollary 
\ref{IntroCorollary}. To the best of the
author's knowledge the tensor product varieties
are new even in the $\mathfrak{gl}_{\sss N}$-case. All
known proofs of the Hall theorem are either
combinatorial (through the Littlewood-Richardson
Rule), or use the relation between the tensor 
product for $GL$ and the restriction for
symmetric groups combined with results
of Borho and MacPherson 
\cite{BorhoMacPherson}. 
The existence of the tensor product varieties
together with bijections $\alpha_n$ and
$\beta_n$ provides a direct proof of
the Hall theorem by showing the role
of Spaltenstein varieties in the
geometric theory of the tensor product.

In the special case of the tensor product
of $n$ fundamental representations 
of $\mathfrak{gl}_{\sss N}$ 
the tensor product variety
is a certain Lagrangian subvariety in the
cotangent bundle of the variety studied
by Grojnowski and Lusztig in
\cite{GrojnowskiLusztig}. However in general
tensor product varieties cannot be
represented in this form (and thus cannot be
used to prove positivity properties of
canonical bases).

The $\mathfrak{gl}_{\sss N}$ 
case is studied in 
\cite{Malkin2000a} without mentioning
quiver terminology (purely in the
language of flags and nilpotent orbits).

\subsubsection{}

Recently Lusztig \cite{Lusztig2000a}
described a locally closed subset 
of the variety $Z$ (no relation with
$Z_D$ of \ref{IntroAppendix})
constructed by
Nakajima (cf. \cite{Nakajima1998}),
such that irreducible components
of this subset (conjecturally) form
a crystal isomorphic to the tensor
product of two irreducible representations
of $\mathfrak{g}$. The relation of this
construction to the tensor product variety
(for $n=2$) will be discussed elsewhere.
In particular, using constructions 
of this paper one can show that the 
set of irreducible
components of the Lusztig's variety 
is in a natural bijection with the
set of irreducible components of
the corresponding tensor product variety.
Note however
that Lusztig's construction does not
produce a direct sum decomposition
for a tensor product, nor it can
be generalized to a case of more
than two multiples. On the other hand
it provides a geometric analogue of
a factor of the universal enveloping
algebra of $\mathfrak{g}$ isomorphic
to the tensor product of a highest weight 
and a lowest weight representations.

\subsubsection{}

An important source of inspiration for the 
author was a work of
A. Braverman and D. Gaitsgory 
\cite{BravermanGaitsgory1999}
where they constructed
geometric crystals (together
with a crystal tensor product)
via the Affine Grassmannian.
Note that Nakajima's construction
of crystals (cf. 
\ref{IntroNakajimaCrystals})
uses a fiber of a 
resolution of singularities 
of a singular variety 
(quiver variety $\mathfrak{M}_0 (v,w)$), 
while  
Braverman and Gaitsgory use cells of a
perverse stratification
of a singular variety 
(the Affine Grassmannian).
Similar difference exists between
the constructions of the crystal
tensor product.

\subsubsection{}
Though the ground field is $\mathbb{C}$
throughout the paper the multiplicity 
varieties are defined over arbitrary 
fields, and it is shown in 
\cite{Malkin2001a} that the 
number of $\mathbb{F}_q$-rational
points of ${^n \mathcal{S}} (
\mu^0; \mu^1, \ldots , \mu^n )$ is
given by a polynomial in $q$ with
the leading coefficient equal to 
$\dim_{\mathbb{C}} \Hom_{\mathfrak{g}'}
(L (\mu^0), L (\mu^1) \otimes \ldots
\otimes L (\mu^n))$ 
(cf. Corollary \ref{IntroCorollary}).
This statement is a direct generalization
of the Hall theorem 
(cf. \cite{Hall1959}, 
\cite[Chapter II]{Macdonald1995})
giving the number
of $\mathbb{F}_q$-rational points in
a Spaltenstein variety for 
$\mathfrak{gl}_{\sss N}$.

\subsubsection{}

Throughout the paper the following 
conventions are used:
the ground field is $\mathbb{C}$;
``closed'', ``locally closed'', etc.,
refer to the Zariski topology;
``fibration'' means 
``locally trivial fibration'', where 
``locally'' refers to the Zariski
topology, however trivialization
is analytic (not regular).

\subsubsection{Acknowledgements}
This paper owns its very existence
to Igor Frenkel and Lev Rozansky. 
All ideas described below crystallized 
during our numerous discussions   
in the summer 2000.

\tableofcontents

\section{Crystals}
\label{CrystalsSection}

Crystals were unearthed by M. Kashiwara 
\cite{Kashiwara1990,
Kashiwara1991, Kashiwara1994}.
An excellent survey of crystals as well
as some new results are given by 
A. Joseph in
\cite[Chapters 5, 6]{Joseph1995}.

\subsection{Weights and roots}
\label{WeightsRoots}
Let $\mathfrak{g}$ be a reductive or 
a Kac-Moody Lie algebra, 
$I$ be the set of vertices of the
Dynkin graph of $\mathfrak{g}$.
It is assumed throughout the paper
that $\mathfrak{g}$ is simply laced.
The weight and coweight lattices are 
identified with $\mathbb{Z} [I]$, so that
the natural pairing between them 
(denoted $<,>$) becomes
\begin{equation}\nonumber
<v, u> = 
\sum_{i \in I} 
v_i u_i \; ,
\end{equation}
where $v_i$ denotes 
the $i$-th component of 
$v \in \mathbb{Z} [I]$.

Let $A$ be the Cartan matrix of
$\mathfrak{g}$. Then the simple root
$\Hat{i}$,
corresponding to a simple weight 
$i \in I$ is given by 
$A i = \sum_{j \in I} A_{ji} j
\in \mathbb{Z}_{\geq 0} [I]$. 

\subsection{Definition of 
$\mathfrak{g}$-crystals}
\indent\par\noindent
A $\mathfrak{g}$-\emph{crystal} is a tuple 
$(\mathcal{A}, wt, 
\{ \varepsilon_i \}_{i \in I} ,
\{ \varphi_i \}_{i \in I} ,
\{ \Tilde{e}_i \}_{i \in I} ,
\{ \Tilde{f}_i \}_{i \in I} )$, where
\begin{itemize}
\item
$\mathcal{A}$ is a set,
\item
$wt$ is a map from $\mathcal{A}$ 
to $\mathcal{Q}_{\mathfrak{g}}$,
\item
$\varepsilon_i$ and $\varphi_i$ are maps from
$\mathcal{A}$ to $\mathbb{Z}$,
\item
$\Tilde{e}_i$ and $\Tilde{f}_i$ are maps 
from $\mathcal{A}$ to 
$\mathcal{A} \cup \{ 0 \}$. 
\end{itemize}
These data should satisfy the following axioms:
\begin{itemize}
\item
$wt (\Tilde{e}_i a) = wt (a) + \Hat{i}$, 
\quad
$\varphi_i (\Tilde{e}_i a) = \varphi_i (a) + 1$, 
\quad
$\varepsilon_i (\Tilde{e}_i a) = 
\varepsilon_i (a) - 1$, 
\quad
for any $i \in I$ and $a \in \mathcal{A}$ 
such that $\Tilde{e}_i a \neq 0$,
\item
$wt (\Tilde{f}_i a) = wt (a) - \Hat{i}$, 
\quad
$\varphi_i (\Tilde{f}_i a) = \varphi_i (a) - 1$, 
\quad
$\varepsilon_i (\Tilde{f}_i a) = 
\varepsilon_i (a) + 1$,
\quad
for any $i \in I$ and $a \in \mathcal{A}$ 
such that $\Tilde{f}_i a \neq 0$,
\item
$(wt (a))_i
= \varphi_i (a) - \varepsilon_i (a)$
for any $i \in I$ and $a \in \mathcal{A}$,
\item
$\Tilde{f}_i a = b$ if and only if 
$\Tilde{e}_i b = a$, where
$i \in I$, and $a,b \in \mathcal{A}$.
\end{itemize}
The maps $\Tilde{e}_i$ and $\Tilde{f}_i$
are called Kashiwara's operators, and
the map $wt$ is called the weight function.

\begin{remark}
In a more general definition 
of crystals the maps $\varepsilon_i$ and 
$\varphi_i$ are allowed to have infinite values.
\end{remark}

A $\mathfrak{g}$-crystal 
$(\mathcal{A}, wt, 
\{ \varepsilon_i \}_{i \in I} ,
\{ \varphi_i \}_{i \in I} ,
\{ \Tilde{e}_i \}_{i \in I} ,
\{ \Tilde{f}_i \}_{i \in I} )$
is called \emph{trivial} if
\begin{equation}\nonumber
\begin{split}
wt (a) = 0 
&\text{ for any $a \in \mathcal{A}$,} \\  
\varepsilon_i (a) = \varphi_i (a) = 0 
&\text{ for any $i \in I$ and
$a \in \mathcal{A}$,} \\ 
\Tilde{e}_i a = \Tilde{f}_i a = 0 
&\text{ for any $i \in I$ and 
$a \in \mathcal{A}$.} \\ 
\end{split}
\end{equation}
Any set $\mathcal{A}$ can be equipped with the
trivial crystal structure as above.

A $\mathfrak{g}$-crystal 
$(\mathcal{A}, wt, 
\{ \varepsilon_i \}_{i \in I} ,
\{ \varphi_i \}_{i \in I} ,
\{ \Tilde{e}_i \}_{i \in I} ,
\{ \Tilde{f}_i \}_{i \in I} )$
is called \emph{normal} if
\begin{equation}\nonumber
\begin{split}
\varepsilon_i (a) &= \max 
\{ n \; | \; \Tilde{e}_i^n a \neq 0 \} \; , 
\\ \varphi_i (a) &= \max 
\{ n \; | \; \Tilde{f}_i^n a \neq 0 \} \; ,
\end{split}
\end{equation}
for any $i \in I$ and 
$a \in \mathcal{A}$. Thus in a normal 
$\mathfrak{g}$-crystal the maps 
$\varepsilon_i$ and $\varphi_i$
are uniquely determined by the action
of $\Tilde{e}_i$ and $\Tilde{f}_i$.
A trivial crystal is normal.

In the rest of the paper all 
$\mathfrak{g}$-crystals are assumed normal,
and thus the maps
$\varepsilon_i$ and $\varphi_i$ are usually
omitted.

By abuse of notation 
a $\mathfrak{g}$-crystal
$(\mathcal{A}, wt, 
\{ \varepsilon_i \}_{i \in I} ,
\{ \varphi_i \}_{i \in I} ,
\{ \Tilde{e}_i \}_{i \in I} ,
\{ \Tilde{f}_i \}_{i \in I} )$ is 
sometimes denoted simply by 
$\mathcal{A}$.

An \emph{isomorphism} of two 
$\mathfrak{g}$-crystals 
$\mathcal{A}$ and
$\mathcal{B}$
is a bijection between the sets 
$\mathcal{A}$ and $\mathcal{B}$ 
commuting with the action of
the operators $\Tilde{e}_i$ and
$\Tilde{f}_i$, and the functions 
$wt$, $\varepsilon_i$, and $\varphi_i$. 

The \emph{direct sum} 
$\mathcal{A} \oplus \mathcal{B}$
of two $\mathfrak{g}$-crystals 
$\mathcal{A}$ and $\mathcal{B}$
is their disjoint union as sets
with the maps $\Tilde{e}_i$,
$\Tilde{f}_i$, $wt$,
$\varepsilon_i$, and $\varphi_i$ 
acting on each
component of the union separately. 

\subsection{Tensor product of 
$\mathfrak{g}$-crystals}
\label{CrystalTensorProduct}

The \emph{tensor product} 
$\mathcal{A} \otimes \mathcal{B}$
of two $\mathfrak{g}$-crystals 
$\mathcal{A}$ and $\mathcal{B}$
is their direct product as sets
equipped with the 
following crystal structure:
\begin{equation}
\label{DefinitionOfCrystalProduct}
\begin{split}
wt ((a,b)) &= wt (a) + wt (b)  \; ,
\\
\varepsilon_i ((a,b)) &=
\max \{ \varepsilon_i (a), \varepsilon_i (a)
+ \varepsilon_i (b) - \varphi_i (a) \} \; ,
\\
\varphi_i ((a,b)) &=
\max \{ \varphi_i (b), \varphi_i (a)
+ \varphi_i (b)- \varepsilon_i (b) \} \; ,
\\
\Tilde{e}_i ((a,b)) &=
\begin{cases}
(\Tilde{e}_i a \; , \; b)
&\text{ if } 
\varphi_i (a) \geq \varepsilon_i (b) \; ,\\
(a \; , \; \Tilde{e}_i b)
&\text{ if } 
\varphi_i (a) < \varepsilon_i (b) \; ,
\end{cases} 
\\
\Tilde{f}_i ((a,b)) &=
\begin{cases}
(\Tilde{f}_i a \; , \;  b)
&\text{ if } 
\varphi_i (a) > \varepsilon_i (b) \; ,\\
(a \; , \; \Tilde{f}_i b) 
&\text{ if } 
\varphi_i (a) \leq \varepsilon_i (b) \; .
\end{cases}
\end{split}
\end{equation}
Here $(a,0)$, $(0,b)$, and $(0,0)$ are
identified with $0$.
One can check that the set 
$\mathcal{A} \times \mathcal{B}$ with the
above structure satisfies all the axioms
of a (normal) crystal, and
that the tensor product of crystals
is associative. 

\begin{remark}
Tensor product of crystals 
is not commutative.
\end{remark}

\begin{remark}
The tensor product 
(in any order) of a crystal 
$\mathcal{A}$ with a trivial crystal 
$\mathcal{B}$ is isomorphic to the direct sum 
$\bigoplus_{b \in \mathcal{B}} \mathcal{A}_b$
where each $\mathcal{A}_b$ is isomorphic to 
$\mathcal{A}$.
\end{remark}

\subsection{Highest weight crystals and 
closed families}
A crystal $\mathcal{A}$ is a 
\emph{highest weight}
crystal with highest weight 
$\lambda \in \mathcal{Q}_{\mathfrak{g}}$ if
there exists an element 
$a_{\lambda} \in \mathcal{A}$ such that:
\begin{itemize}
\item 
$wt ( a_{\lambda} ) = \lambda$ ,
\item
$\Tilde{e}_i a_{\lambda} = 0$ 
for any $i \in I$,
\item
any element of $\mathcal{A}$ can be obtained
from $a_{\lambda}$ by successive applications
of the operators $\Tilde{f}_i$.
\end{itemize}

Consider a family of highest weight 
normal crystals
$\{ \mathcal{A} (\lambda ) \}_{\lambda 
\in \mathcal{J}}$ labeled by a set
$\mathcal{J} \subset 
\mathcal{Q}_{\mathfrak{g}}$ (the highest
weight of $\mathcal{A} ( \lambda )$ is
$\lambda$). The family
$\{ \mathcal{A} (\lambda ) \}_{\lambda 
\in \mathcal{J}}$ is called 
\emph{strictly closed} 
(with respect to tensor products)
if the tensor product of any two members of
the family is isomorphic to a direct sum 
of members of the family:
\begin{equation}\nonumber
\mathcal{A} ( \mu^1 ) \otimes
\mathcal{A} ( \mu^2 ) = 
\bigoplus_{\lambda \in \mathcal{J}}
\mathcal{U} ((\mu^1 , \mu^2 ), \lambda )
\otimes \mathcal{A} ( \lambda ) \; ,
\end{equation}
where 
$\mathcal{U} ((\mu^1 , \mu^2 ), \lambda )$
is a set equipped 
with the trivial crystal structure.
The family
$\{ \mathcal{A} (\lambda ) \}_{\lambda 
\in \mathcal{J}}$ is called 
\emph{closed} 
if the tensor product 
$\mathcal{A} (\mu^1) \otimes
\mathcal{A} (\mu^2)$ 
of any two members of the family
contains $\mathcal{A}_{\mu^1+\mu^2}$
as a direct summand. Any strictly closed
family is closed. The converse is also true,
as a corollary of Theorem 
\ref{JosephTheorem}.

Let $\mathcal{Q}^{+}_{\mathfrak{g}}
\subset \mathcal{Q}_{\mathfrak{g}}$ be
the set of highest weights of integrable 
highest weights modules of $\mathfrak{g}$
(in the reductive case a module is called
integrable if it is derived from a 
polynomial representation of the corresponding
connected simply connected reductive group, 
in the Kac-Moody
case the highest weight should be a 
positive linear combination of the fundamental
weights). The original motivation for the 
introduction  of crystals was the discovery 
by M. Kashiwara \cite{Kashiwara1991} 
and G. Lusztig  \cite{Lusztig1991a}
of canonical (or crystal) bases in 
integrable highest 
weight modules of a (quantum) Kac-Moody algebra. 
These bases have many favorable properties, one
of which is that 
as sets they are equipped with a crystal 
structure. In other words to each 
irreducible integrable highest weight module
$L (\lambda )$ corresponds a normal crystal
$\mathcal{L} (\lambda )$
(crystal of the canonical basis). In this 
way one obtains a strictly closed
family of crystals
$\{ \mathcal{L} (\lambda ) \}_{\lambda 
\in \mathcal{Q}^+_{\mathfrak{g}}}$, 
satisfying the following two properties:
\begin{itemize}
\item
the cardinality of $\mathcal{L} (\lambda )$
is equal to the dimension of $L (\lambda)$;
\item
one has the following tensor product 
decompositions for 
$\mathfrak{g}$-modules and
$\mathfrak{g}$-crystals:
\begin{equation}\nonumber
\begin{split}
L ( \mu^1 ) \otimes
L ( \mu^2 ) &= 
\bigoplus_{\lambda 
\in \mathcal{Q}^+_{\mathfrak{g}}}
C ((\mu^1 , \mu^2 ), \lambda )
\otimes L ( \lambda ) \; ,
\\
\mathcal{L} ( \mu^1 ) \otimes
\mathcal{L} ( \mu^2 ) &= 
\bigoplus_{\lambda 
\in \mathcal{Q}^+_{\mathfrak{g}}}
\mathcal{C} ((\mu^1 , \mu^2 ), \lambda )
\otimes \mathcal{L} ( \lambda ) \; ,
\end{split}
\end{equation}
where the cardinality of the set
(trivial crystal)
$\mathcal{C} 
((\mu^1 , \mu^2 ), \lambda )$
is equal to the dimension of the linear space
(trivial $\mathfrak{g}$-module)
$C ((\mu^1 , \mu^2 ), \lambda )$.
\end{itemize}

\noindent
The aim of this paper is to 
construct another strictly closed 
family of $\mathfrak{g}$-crystals 
$\{ \mathcal{M} (\lambda ) \}_{\lambda 
\in \mathcal{Q}^+_{\mathfrak{g}}}$
(for $\mathfrak{g}$ being $gl_{\sss N}$
or a symmetric Kac-Moody algebra), using
geometry associated to $\mathfrak{g}$. The
following crucial theorem ensures that this
family is isomorphic to the family
$\{ \mathcal{L} (\lambda ) \}_{\lambda 
\in \mathcal{Q}^+_{\mathfrak{g}}}$
of crystals of canonical bases 
(two families of crystals labeled
by the same index set 
are called isomorphic
if the corresponding members of the
families are isomorphic as crystals). 
\begin{numtheorem}\label{JosephTheorem}
There exists a unique 
(up to an isomorphism) closed family of
$\mathfrak{g}$-crystals labeled
by $\mathcal{Q}^+_{\mathfrak{g}}$. 
\end{numtheorem}
\begin{proof}
For the proof in Kac-Moody case see
\cite[Proposition 6.4.21]{Joseph1995}. 
The statement for a reductive $\mathfrak{g}$ 
easily follows from the statement for
the factor of $\mathfrak{g}$ by its center.
\end{proof} 

\section{Nakajima varieties and 
tensor product varieties}
\label{QuiverTensorSection}

\subsection{Oriented graphs and path 
algebras}
\label{PathAlgebraSubsection}

Let $I$ be the set of vertices 
of the Dynkin graph of a simple 
simply laced Lie 
algebra $\mathfrak{g}$. Let $H$ be the set 
of pairs consisting of an
edge of the Dynkin graph of $\mathfrak{g}$ 
and an orientation of this edge. 
The target (resp. source) vertex of 
$h \in H$ is denoted by $\In (h)$ 
(resp. $\Out (h)$). Thus $(I,H)$ is
an oriented graph (note that it has twice 
as many edges as the Dynkin graph of 
$\mathfrak{g}$ has). For $h \in H$ let 
$\Bar{h}$ be the same edge
but with opposite orientation (i.e. 
$\Bar{h}$ is the unique element of $H$ 
such that
$\In (\Bar{h}) = \Out (h)$ and 
$\Out (\Bar{h}) = \In (h)$).

Let $X$ be a symmetric $I \times I$ matrix 
uniquely defined by the following equation
\begin{equation}\label{DefinitionOfX}
< X v , u > =
\sum_{h \in H} v_{\In (h)} 
u_{\Out (h)} \; ,
\end{equation}
for any $v$, $u \in \mathbb{Z} [I]$.
Note that $X = 2 \Id - A$, where
$A$ is the Cartan matrix of $\mathfrak{g}$.

Let $\mathcal{F}$ be the path algebra of 
the oriented graph $(I,H)$ over 
$\mathbb{C}$. Fix a function 
$\varepsilon : H \rightarrow \mathbb{C}^*$ 
such that
$\varepsilon (h) + 
\varepsilon (\Bar{h}) = 0$ for any $h \in H$. 
Let 
\begin{equation}\label{DefinitionOfThetaI}
\theta_i = 
\sum_{\substack{h \in H \\ \In (h) = i}}
\varepsilon (h) h \Bar{h} 
\in \mathcal{F} \; ,
\end{equation} 
and
let $\mathcal{F}_0$ be the factor algebra
of $\mathcal{F}$ by the two-sided ideal 
generated by elements $\theta_i$ for all 
$i \in I$. The algebra $\mathcal{F}_0$ 
is called the 
\emph{preprojective algebra}. 
It was introduced by 
Gelfand and Ponomarev.

Note that an $\mathcal{F}$-module
is just a $\mathbb{Z} [I]$-graded 
$\mathbb{C}$-linear space $V$ 
together with a collection of linear maps
$x= \{ x_h \in \Hom_{\mathbb{C}}
( V_{\Out (h)}, V_{\In (h)} ) 
\}_{h \in H}$.
This $\mathcal{F}$-module is denoted
by $(V, x)$. It is always assumed 
below that $V$ is finite dimensional.

Given a finite dimensional
$\mathbb{Z}[I]$-graded
$\mathbb{C}$-linear space $V$,
a set of linear maps 
$x = \{ x_h \in \Hom_{\mathbb{C}}
(V_{\Out (h)},
V_{\In (h)}) \}_{h \in H}$
such that $(V, x)$ is
an $\mathcal{F}$-module is called
a \emph{representation}
of $\mathcal{F}$ in $V$.

Let $(V, x)$ be an $\mathcal{F}$-module
and $W$ be a $\mathbb{Z} [I]$-graded
subspace of $V$ such that 
$x_h W_{\Out (h)} \subset
W_{\In (h)}$ for any $h \in H$. 
Then $(W, \; x|_{\sss W})$ (resp.
$(V/W, \; x|_{\sss V/W})$)
is an $\mathcal{F}$-submodule 
(resp. factor module) of
$(V , x)$. Sometimes just
$W$ (resp. $V/W$)
is called a submodule (resp. 
factor module) because the 
representation of $\mathcal{F}$
is uniquely defined by the restriction.

An $\mathcal{F}$-module $(V, x)$
is also an $\mathcal{F}_0$-module if
and only if
$\sum_{\substack{h \in H \\ \In (h) = i}}
\varepsilon (h) x_h x_{\Bar{h}} = 0$
for any $i \in I$.

A $k$-tuple $h_1 , h_2 , \ldots , h_k$ of
elements of $H$ is called a \emph{path}
of length $k$ if 
$\In ( h_i ) = \Out ( h_{i-1} )$ for 
$i = 2,  \ldots , k$.
An $\mathcal{F}$-module $(V , x)$
is called \emph{nilpotent} if 
there exists $N \in \mathbb{Z}_{> 0}$
such that 
$x_{h_1} x_{h_2} \ldots x_{h_N}=0$
for any path 
$h_1 , h_2 , \ldots , h_N$
of length $N$. 

The following proposition is due to
Lusztig \cite{Lusztig1991a}.

\begin{proposition}\indent\par
\begin{alphenum}
\item\label{PathAlgebraNilpotentClosed}
An $\mathcal{F}$-module $(V, x)$
is nilpotent if and only if 
$x_{h_1} x_{h_2} 
\ldots x_{h_{|\dim \mathbf{V}|}}=0$
(i.e. any path of length 
$| \dim V |$ is in the kernel
of the representation $x$).
\item\label{PathAlgebraNilpotentADE}
Any $\mathcal{F}_0$-module is nilpotent
as an $\mathcal{F}$-module.. 
\end{alphenum}
\end{proposition}
\begin{proof}
The ``if'' part of 
\ref{PathAlgebraNilpotentClosed} follows from 
the definition of a nilpotent $\mathcal{F}$-module.

The ``only if'' part of 
\ref{PathAlgebraNilpotentClosed} is proved in
\cite[Lemma 1.8]{Lusztig1991a}. 

\ref{PathAlgebraNilpotentADE} is proved in
\cite[Section 12]{Lusztig1991a}. Recall that the
Dynkin graph is assumed to be 
of finite (ADE) type.

\end{proof}

\subsection{Quiver varieties}

Quiver varieties were introduced by Nakajima
\cite{Nakajima1994, Nakajima1998}
as a generalization of the moduli space of 
Yang-Mills instantons on an ALE space
(cf. \cite{KronheimerNakajima}). On the
other hand they also generalize Spaltenstein
variety for $GL (N)$ 
(the variety of all parabolic subgroups $P$ 
such that a given unipotent element 
$u \in GL (N)$ belongs to the unipotent
radical of $P$). Nakajima has shown 
in \cite{Nakajima1994, Nakajima1998, 
Nakajima2001} that the quiver varieties 
play the same role in
the representation theory of (quantum)
Kac-Moody algebras as Spaltenstein varieties
do in the representation theory of 
$GL (N)$ (cf. 
\cite{Ginzburg1991, KazhdanLusztig1987}).

Let $D$, $V$ be finite dimensional 
$\mathbb{Z} [I]$-graded
$\mathbb{C}$-linear spaces.
An \emph{ADHM datum}\footnote{ADHM 
stands for Atiyah, Drinfeld, 
Hitchin, Manin
\cite{AtiyahDrinfeldHitchinManin}} 
(on $D$, $V$) is a triple 
$(x, p, q)$, where

$x = \{ x_h \}_{h \in H}$ is a representation
of $\mathcal{F}$ in $V$,

$p \in \Hom_{\mathbb{C}^I} 
(D, V)$ is a $\mathbb{Z} [I]$-graded
$\mathbb{C}$-linear map from 
$D$ to $V$,

$q \in \Hom_{\mathbb{C}^I} (V, D)$
is a $\mathbb{Z} [I]$-graded
$\mathbb{C}$-linear map from 
$V$ to $D$,

and $x$, $p$, $q$ satisfy the following
equation (in $\End_{\mathbb{C}} V_i$)
\begin{equation}\label{ADHMequation}
\sum_{\substack{
h \in H \\ \In (h) = i }}
\varepsilon (h) x_h x_{\Bar{h}} - 
p_i q_i = 0
\end{equation}
for all $i \in I$.
The set of all ADHM data on $D$, $V$
form an affine variety denoted by
$\Lambda_{D, V}$. 

Let $D = \{ 0 \}$. Then $p=q=0$ and
$x$ gives a representation of the 
preprojective algebra 
$\mathcal{F}_0$ in $V$. 
The variety $\Lambda_{\{ 0 \}, V}$
is denoted simply by $\Lambda_V$,
and an element $(x,0,0)$ of $\Lambda_V$ 
is usually written simply as $x$. Thus 
$\Lambda_V$ is the variety
of all representations of the preprojective
algebra $\mathcal{F}_0$ in $V$. 
It follows
from \ref{PathAlgebraNilpotentADE} that
if $x \in \Lambda_V$ then
$(V, x)$ is a nilpotent $\mathcal{F}$-module.
The variety $\Lambda_V$
was introduced by Lusztig 
\cite[Section 12]{Lusztig1991a}. He also
defined an open subset ${^{\delta} \Lambda}_V$
of $\Lambda_V$ given by
\begin{equation}\label{DeltaLambdaV}
{^{\delta} \Lambda}_V = \{
x \in \Lambda_V \; | \;
\dim \Coker (\bigoplus_{\substack{h \in H \\
\In (h) = i}} x_h ) \leq \delta_i 
\text{ for any $i \in I$ }\} \: , 
\end{equation}
where $\delta \in \mathbb{Z}_{\geq 0} [I]$.
The following proposition is
proven in \cite[Section 12]{Lusztig1991a}.

\begin{proposition} Let $v = \dim V$. Then
\begin{alphenum}
\item 
${^{\delta} \Lambda}_V$ is either 
empty or has pure dimension 
$\frac{1}{2}<X v, v>$,
\item
\label{DimOfLambdaNilpotent}
$\Lambda_V$ has
pure dimension 
$\frac{1}{2}<X v, v>$.
\end{alphenum}
\end{proposition}

\subsection{Stability and $G_V$-action}
\label{Stability}
Let $(x, p, q) \in \Lambda_{D, V}$ 
and $E$ be a graded subspace of $V$.
Then $\overline{E}$ (resp.
$\underline{E}$) denotes 
the smallest $\mathcal{F}$-submodule
of $(V,x)$ containing $E$
(resp. the largest $\mathcal{F}$-submodule
of $(V,x)$ contained in $E$). 

Note that the triple
$(x|_{\underline{q^{-1} (0)}}, 0, 0)$ 
(resp. 
$(x|_{V/\overline{p(D)}}, 0, 0)$)
belongs to 
$\Lambda_{\underline{q^{-1} (0)}}$ 
(resp. 
$\Lambda_{V/\overline{p(D)}}$).
Nakajima \cite{Nakajima1994, Nakajima1998}
defined open (possibly empty)
subsets of the variety
$\Lambda_{D, V}$
given as sets of points satisfying some
stability conditions. Two examples of such 
subsets are  
$\Lambda^{s}_{D, V}$
(the set of \emph{stable} points) and
$\Lambda^{\ast s}_{D, V}$
(the set of \emph{$\ast$-stable} points).
The set 
$\Lambda^{s}_{D, V}$ 
(resp. 
$\Lambda^{\ast s}_{D, V}$)
is the set of all triples
$(x, p, q) \in 
\Lambda_{D, V}$ 
such that $\overline{p(D)}= V$
(resp. $\underline{q^{-1} (0)}=0$). Let
$\Lambda^{s, \ast s}_{D, V}=
\Lambda^{s}_{D, V} 
\cap
\Lambda^{\ast s}_{D, V}$.

Let $G_V = \Aut_{\mathbb{C}^I} V = 
\prod_{i \in I} GL ( V_i)$ be 
the group of graded automorphisms
of $V$. The group $G_V$ acts on 
$\Lambda_V$ and $\Lambda_{D, V}$
as follows:
\begin{equation}\nonumber
\begin{split}
(g, x) &\rightarrow x^g
\\
(g, (x,p,q)) &\rightarrow (x^g,p^g,q^g) \; ,
\end{split}
\end{equation}
where $g \in G_V$, 
$x \in \Lambda_V$
(resp. $(x,p,q) \in
\Lambda_{D, V}$),
and
\begin{equation}\nonumber
x^g_h = g_{\In (h)} x_h g_{\Out (h)}^{-1} \; , 
\qquad p^g_i = g_i p_i \; , \qquad 
q^g_i = q_i g_i^{-1} \; ,
\end{equation} 
for $h \in H$, $i \in I$.

The subsets 
$\Lambda^{s}_{D, V}$,
$\Lambda^{\ast s}_{D, V}$,
and $\Lambda^{s, \ast s}_{D, V}$, are
preserved by the $G_V$-action.

The following proposition due to
Nakajima and Crawley-Boevey
explains the importance
of the open sets 
$\Lambda^{s}_{D, V}$ and
$\Lambda^{\ast s}_{D, V}$.

\begin{proposition}
Let $d = \dim D$, $v = \dim V$. Then
\begin{alphenum}
\item\label{DimOfLambdaStable}
$\Lambda^{s}_{D, V}$
is empty or an irreducible 
smooth variety of dimension
\begin{equation}\nonumber
\dim \Lambda^{s}_{D, V} =
<X v , v> + 2 <d, v> - <v, v> \; ; 
\end{equation}
\item\label{DimOfLambdaAstStable}
$\Lambda^{\ast s}_{D, V}$
is empty or an irreducible 
smooth variety of dimension
\begin{equation}\nonumber
\dim \Lambda^{\ast s}_{D, V} =
<X v , v> + 2 <d, v> - <v, v> \; ; 
\end{equation}
\item\label{IrrStableAstStable}
$\Lambda^{s, \ast s}_{D, V}$
is empty or an irreducible 
smooth variety of dimension
\begin{equation}\nonumber
\dim
\Lambda^{s, \ast s}_{D, V} =
<X v , v> + 2 <d, v> - <v, v> \; ; 
\end{equation}
\item\label{LambdaSAstSEmpty}
$\Lambda^{s, \ast s}_{D, V}$
is non-empty if and only if
$d - 2v + Xv \in \mathbb{Z}_{\geq 0}[I]$;
\item\label{StableFreeG} 
$G_V$-action is free on 
$\Lambda^{s}_{D, V}$,
$\Lambda^{\ast s}_{D, V}$, and
$\Lambda^{s, \ast s}_{D, V}$.
\end{alphenum}
\end{proposition}
See \eqref{DefinitionOfX} for the
definition of the matrix $X$.
\begin{proof}
The smoothness part of
\ref{DimOfLambdaAstStable}
and the statement \ref{StableFreeG} 
in the case of $\Lambda^{\ast s}_{D, V}$
are proven by Nakajima 
\cite[Lemma 3.10]{Nakajima1998}.
The fact that $\Lambda^{\ast s}_{D, V}$
is connected is proven by 
Crawley-Boevey 
\cite[Remarks after the 
Introduction]{CrawleyBoevey2000}).

The proof of \ref{DimOfLambdaStable} is 
analogous to the one of 
\ref{DimOfLambdaAstStable}, or one
can deduce the former from the latter by 
transposing $(x,p,q)$. Either argument  also 
proves the part of \ref{StableFreeG} concerning
$\Lambda^{s}_{D, V}$.

\ref{IrrStableAstStable} (and the corresponding
part of \ref{StableFreeG}) follows from the 
definition of 
$\Lambda^{s, \ast s}_{D, V}$ 
($\Lambda^{s, \ast s}_{D, V} = 
\Lambda^{s}_{D, V} \cap 
\Lambda^{\ast s}_{D, V}$),
\ref{DimOfLambdaStable} and
\ref{DimOfLambdaAstStable}. 

\ref{LambdaSAstSEmpty} is proven by 
Nakajima \cite[Sections 10.5 -- 
10.9]{Nakajima1998}
(see also \cite[Proposition 
1.11]{Lusztig2000a}).
\end{proof}

It follows from \ref{StableFreeG} and 
\ref{DimOfLambdaStable} that
$\mathfrak{M}^{s}_{D, V} =
\Lambda^{s}_{D, V} /
G_V$ is naturally an algebraic 
variety of dimension 
\begin{equation}\nonumber
\dim \mathfrak{M}^s_{D, V}=
<X v, v> + 2 <d, v> - 2 <v, v> \; ,
\end{equation}
and the natural projection
$\Lambda^{s}_{D, V}
\rightarrow \mathfrak{M}^s_{D, V}$
is a principal $G_V$-bundle.
The same holds for
$\mathfrak{M}^{\ast s}_{D, V}=
\Lambda^{\ast s}_{D, V} / G_V$ and
$\mathfrak{M}^{s, \ast s}_{D, V}=
\Lambda^{s, \ast s}_{D, V} /G_V$. 
The varieties 
$\mathfrak{M}^{s}_{D, V}$,
$\mathfrak{M}^{\ast s}_{D, V}$,
and $\mathfrak{M}^{s, \ast s}_{D, V}$ 
are called \emph{quiver varieties}
(cf. \cite{Nakajima1994, Nakajima1998}).

Since
$\mathfrak{M}^{s}_{D,V}$,
$\mathfrak{M}^{\ast s}_{D,V}$, and
$\mathfrak{M}^{s, \ast s}_{D,V}$,
are $G_V$-orbit spaces, they do not
depend on $V$ as long as $\dim V =v$
is fixed. Also if
$\dim D = \dim D'$ the above varieties
are isomorphic (non canonically) to
the corresponding varieties defined using
$D'$ instead of $D$. As 
$D$ never vary,
the following notation is often used
(cf. \cite{Nakajima1994, Nakajima1998}):
$\mathfrak{M}^{s} (d,v) =
\mathfrak{M}^{s}_{D,V}$,
$\mathfrak{M}^{\ast s} (d,v)=
\mathfrak{M}^{\ast s}_{D,V}$, 
$\mathfrak{M}^{s, \ast s} (d,v)=
\mathfrak{M}^{s, \ast s}_{D,V}$.

\subsection{The variety
$\mathfrak{M}^s (d,v_0,v)$}
\label{MdvvSection}

Lusztig 
\cite[Section 1]{Lusztig2000a} 
introduced two kinds of subsets of
$\Lambda_{D,V}$. Namely, 
given a $\mathbb{Z} [I]$-graded
subspace $U \subset V$ let 
\begin{equation}\nonumber
\Lambda_{D, V, U} =
\{ (x,p,q) \in
\Lambda_{D, V} \; | \;
\underline{ q^{-1} (0)} = U \} \; ,
\end{equation}
and given 
$v_0 \in \mathbb{Z}_{\geq 0} [I]$ let
\begin{equation}\nonumber
\Lambda_{D,V,v_0}  =
\{ (x,p,q) \in \Lambda_{D,V} \; | \;
\dim \underline{q^{-1} (0)} = 
v - v_0 \} \; .
\end{equation}
The following proposition is proven
in \cite[1.8]{Lusztig2000a}.
\begin{proposition}
$\Lambda_{D, V, U}$ and
$\Lambda_{D,V,v_0}$ are locally closed
in $\Lambda_{D, V}$.
\end{proposition}
Let 
$\Lambda^s_{D, V, U} =
\Lambda_{D, V, U} \cap \Lambda^s_{D, V}$,
$\Lambda^s_{D,V,v_0} =
\Lambda_{D,V,v_0}  \cap 
\Lambda^s_{D, V}$, and
\begin{equation}\nonumber
\mathfrak{M}^s (d,v_0,v) =
\Lambda^s_{D,V,v_0} / G_V  \; .
\end{equation}

\subsection{The set 
$\mathcal{M} (d, v_0, v)$} 
\label{DefinitionOfMCal} 
Let $U$, $T$ be complimentary 
$\mathbb{Z} [I]$-graded subspaces of $V$. 
Then one has a natural map 
(cf. \cite{Lusztig2000a})
\begin{equation}\nonumber
\gamma : \Lambda_{D, V, U}
\rightarrow 
\Lambda_U \times 
\Lambda^{\ast s}_{D, T} \; ,
\end{equation}
given by
\begin{equation}\nonumber
\gamma ((x,p,q)) = (x^{\sss UU},
(x^{\sss TT}, p^{\sss TD} ,
q^{\sss DT})) \;,
\end{equation}
where 
the notation $p^{\sss TD}$ means the 
$D \rightarrow T$ component of the 
block-matrix 
$p: D \rightarrow U \oplus T$, and,
similarly, $x^{\sss UU}$, 
$x^{\sss TT}$, $q^{\sss DT}$ are
block-components of $x$ and $q$.
The map $\gamma$ is well-defined because
for $(x,p,q) \in
\Lambda_{D, V, U}$, the subspace $U$ is
equal to $\underline{q^{-1} (0)}$ 
(i.e. it is the maximal 
$\mathcal{F}$-submodule contained in 
the kernel of $q$). 
The following proposition
is due to Lusztig
\cite[Section 1]{Lusztig2000a}.

\begin{proposition}
Let $d = \dim D$, $v = \dim V$,
$u = \dim U$, $t = \dim T = v - u$.
Then
\begin{alphenum}
\item\label{LambdaDVUGamma}
the map $\gamma$ is a vector bundle
with fibers of dimension
$<d, u> + <X t, u> - <t, u>$;
\item\label{LambdaDVUOpen} 
if $\Lambda^s_{D, V, U}$ is non-empty 
then
$\gamma (\Lambda^s_{D, V, U}) = 
{^{\delta} \Lambda}_U \times
\Lambda^{s, \ast s}_{D, T}$,
where $\delta = d - 2t + Xt$ (the matrix $X$
is defined in \eqref{DefinitionOfX}),
and 
$\Lambda^s_{D, V, U}$ is an open dense
subset in $\gamma^{-1} 
({^{\delta} \Lambda}_U \times
\Lambda^{s, \ast s}_{D, T})$;
\item\label{LambdaDVUDimBig}
$\Lambda^s_{D, V, U}$ is empty or has pure
dimension
\begin{equation}\nonumber
\dim \Lambda^s_{D, V, U} =
\frac{1}{2}<X u, u> + <X t, v> +
<d, v> + <d, t> - <t, v> \; ;  
\end{equation}
\item\label{LambdaDVUDimSmall}
$\Lambda^s_{D,V,v_0}$ is empty or 
has pure dimension
\begin{equation}\nonumber
\dim \Lambda^s_{D,V,v_0}=
\frac{1}{2}<X v, v> + 
\frac{1}{2}<X v_0, v_0> +
<d, v> + <d, v_0> - <v_0, v_0> \; ;  
\end{equation}
\item\label{LambdaDVUDimM}
$\mathfrak{M}^s (d, v_0, v)$ is empty or 
has pure dimension
\begin{multline}\nonumber
\dim \mathfrak{M}^s (d, v_0, v) 
= \\ =
\frac{1}{2}<X v, v> +
\frac{1}{2}<X v_0, v_0> + <d, v> +
<d, v_0> - <v_0, v_0> - <v, v> 
\\ =
\frac{1}{2} 
\dim \mathfrak{M}^s (d, v) +
\frac{1}{2}
\dim \mathfrak{M}^{s, \ast s} (d, v_0) \; .  
\end{multline}
\end{alphenum}
\end{proposition}
\begin{proof}
\ref{LambdaDVUGamma} and 
\ref{LambdaDVUOpen} are proven in
\cite[Proposition 1.16]{Lusztig2000a} (see
the proof of Proposition
\ref{InductionOfPi} for a similar
argument).

\ref{LambdaDVUDimBig} follows
from \ref{LambdaDVUGamma}, 
\ref{LambdaDVUOpen}, 
\ref{IrrStableAstStable}, and
\ref{DimOfLambdaNilpotent}.

Let $Gr_w^V$ denote
the variety of all $\mathbb{Z} [I]$-graded
subspaces of graded dimension
$w$ in $V$. It is a 
$G_V$-homogeneous 
variety with connected stabilizer of a point,
and $\dim Gr_w^V = <w, v-w>$.
The map
$\Lambda^s_{D,V,v_0}\rightarrow
Gr_{v - v_0}^V$ given by
$(x,p,q) \rightarrow
\underline{q^{-1} (0)}$ is a locally 
trivial fibration over 
$Gr_{v - v_0}^V$ with the fiber over 
$U \in Gr_{v - v_0}^V$ equal to
$\Lambda^s_{D, V, U}$. 
Now \ref{LambdaDVUDimSmall} follows
from \ref{LambdaDVUDimBig}.

\ref{LambdaDVUDimM} follows from
\ref{LambdaDVUDimSmall} and 
\ref{StableFreeG}.
\end{proof}

\noindent
The statement \ref{LambdaDVUDimM}
also follows from 
\cite[Proof of Theorem 7.2]{Nakajima1998}.

Let $\mathcal{M} (d, v_0, v)$ 
be the set of irreducible 
components of the variety
$\mathfrak{M}^s ( d, v_0, v )$. Since
$\Lambda^s_{D,V,v_0}$ is the total space 
of a principal $G_V$-bundle over
$\mathfrak{M}^s (d, v_0, v)$, while
$\Lambda^s_{D, V, U }$ 
(where $\dim U = v - v_0$) is a fiber
of the fibration 
$\Lambda^s_{D,V,v_0}\rightarrow
Gr_{v - v_0}^V$ over
a simply connected homogeneous space, one
has natural bijections between
the sets of irreducible components of 
$\Lambda^s_{D,V,v_0}$ and
$\Lambda^s_{D, V, U}$ 
and the set
$\mathcal{M} (d, v_0, v)$.
Abusing notation the former two sets
(of irreducible components) will be also
denoted by 
$\mathcal{M} (d, v_0, v)$.

Let $\mathcal{M} (d, v_0)=
\bigsqcup_{v \in \mathbb{Z} [I]}
\mathcal{M} (d, v_0, v)$.
In \ref{CrystalMDV0} the set
$\mathcal{M} (d, v_0)$ is equipped 
with a structure of $\mathfrak{g}$-crystal,
which is shown (in \ref{Corollary}) to
coincide with the crystal of
the canonical basis of a highest weight
representation of $\mathfrak{g}$.

\subsection{Tensor product varieties and
multiplicity varieties}
\label{DefinitionOfTS}

Given an $n$-tuple $\mathbf{v}$ of elements of 
$\mathbb{Z} [I]$, $\mathbf{v}^k$ denotes the
$k$-th component of $\mathbf{v}$. Thus 
$\mathbf{v} = 
(\mathbf{v}^1, \ldots , \mathbf{v}^n)$,
$\mathbf{v}^k \in \mathbb{Z} [I]$
($1 \leq k \leq n$), 
$\mathbf{v}^k_i \in \mathbb{Z}$
($i \in I$).

Let $D$ be a $\mathbb{Z} [I]$-graded
$\mathbb{C}$-linear space, and 
$\mathbf{D} =
\{ 0 = \mathbf{D}^0 \subset
\mathbf{D}^1 \subset \mathbf{D}^2
\subset \ldots \subset 
\mathbf{D}^n = D \}$ be a partial
$\mathbb{Z} [I]$-graded flag in $D$.
Then $\dim \mathbf{D}$ denotes an
$n$-tuple of elements of $\mathbb{Z} [I]$ 
given as $(\dim \mathbf{D})^k = 
\dim \mathbf{D}^k - \dim \mathbf{D}^{k-1}$
(for $1 \leq k \leq n$).

Recall that if 
$(x, p, q) \in \Lambda_{D, V}$ 
and $E$ is a $\mathbb{Z} [I]$-graded 
subspace of $V$ then
$\overline{E}$ (resp.
$\underline{E}$) denotes 
the smallest $\mathcal{F}$-submodule
of $(V,x)$ containing $E$
(resp. the largest $\mathcal{F}$-submodule
of $(V,x)$ contained in $E$). 

Let $D$, $V$ be $\mathbb{Z} [I]$-graded
$\mathbb{C}$-linear spaces,
$\mathbf{D} =
\{ 0 = \mathbf{D}^0 \subset
\mathbf{D}^1 \subset \mathbf{D}^2
\subset \ldots \subset 
\mathbf{D}^n = D \}$ be a partial
$\mathbb{Z} [I]$-graded flag in $D$.
Then ${^n \Pi}^s_{D, \mathbf{D}, V}$
denotes the set of all
$(x,p,q) \in \Lambda^s_{D,V}$ 
such that 
\begin{equation}\label{PiSubmodules}
\overline{p (\mathbf{D}^k)} 
\subset
\underline{q^{-1} (\mathbf{D}^k)}   
\end{equation}
for all $k = 1, \ldots , n$.

Note that if $(x,p,q) \in
{^n \Pi}^s_{D, \mathbf{D}, V}$
then $\overline{p (D)} = V$ 
(because $(x,p,q) \in \Lambda^s_{D,V}$),
and the triple 
\begin{equation}\nonumber
\bigl(
x|_{\sss \overline{p (\mathbf{D}^k)}/
(\overline{p (\mathbf{D}^k)} \cap 
\underline{q^{-1} 
(\mathbf{D}^{k-1})})} \; , \; 
p|_{\sss \mathbf{D}^k} \mod 
\underline{q^{-1} 
(\mathbf{D}^{k-1})} \; , \;
q|_{\sss \overline{p (\mathbf{D}^k)}} 
\mod \mathbf{D}^{k-1} \bigr) 
\end{equation}
belongs to $\Lambda^{s, \ast s}_{
\mathbf{D}^k/\mathbf{D}^{k-1} \; , \;
\overline{p (\mathbf{D}^k)}/
(\overline{p (\mathbf{D}^k)} \cap 
\underline{q^{-1} (\mathbf{D}^{k-1})})}$
for all $k = 1, \ldots , n$.

Let $\mathbf{v}$ be an $n$-tuple
of elements of 
$\mathbb{Z}_{\geq 0} [I]$.
Then
\begin{equation}\nonumber
{^n \Pi}^s_{D, \mathbf{D}, V, 
\mathbf{v}} =  \{ (x,p,q) \in 
{^n \Pi}^s_{D, \mathbf{D}, V} \; | \;
\dim \biggl(
\overline{p (\mathbf{D}^k)}/
\bigl(\overline{p (\mathbf{D}^k)} \cap 
\underline{q^{-1} (\mathbf{D}^{k-1})}
\bigr) \biggr) = \mathbf{v}^k \} \; .
\end{equation}

Let $\Tilde{\mathbf{v}}$ 
be another $n$-tuple of
elements of $\mathbb{Z}_{\geq 0} [I]$ such that
$\sum_{k=1}^n \mathbf{v}^k +
\sum_{k=1}^n \Tilde{\mathbf{v}}^k = 
\dim V$. Then
\begin{multline}\nonumber
{^n \Pi}^s_{D, \mathbf{D}, V, \mathbf{v}, 
\Tilde{\mathbf{v}}} = \\ =
\{ (x,p,q) \in
{^n \Pi}^s_{D, \mathbf{D}, V, \mathbf{v}}
\; | \;
\dim \biggl(
\bigl(\overline{p (\mathbf{D}^k)} \cap 
\underline{q^{-1} (\mathbf{D}^{k-1})}
\bigr)/
\overline{p (\mathbf{D}^{k-1})} \biggr) =
\Tilde{\mathbf{v}}^k \} \; .
\end{multline}

Let $\mathbf{V} =
( 0 = \mathbf{V}^0
\subset
\Tilde{\mathbf{V}}^1
\subset
\mathbf{V}^1
\subset
\Tilde{\mathbf{V}}^2
\subset
\ldots
\subset
\Tilde{\mathbf{V}}^n
\subset
\mathbf{V}^n = V )$ be
a $2n$-step $\mathbb{Z} [I]$-graded
partial flag in $V$. Then
\begin{equation}\nonumber
{^n \Pi}^s_{D, \mathbf{D}, V, 
\mathbf{V}}  =
\{ (x,p,q) \in
{^n \Pi}^s_{D, \mathbf{D}, V}
\; | \;
\overline{p (\mathbf{D}^k)} = 
\mathbf{V}^k \; , \;
\underline{q^{-1} (\mathbf{D}^{k-1})}
\cap
\overline{p (\mathbf{D}^k)} =
\Tilde{\mathbf{V}}^k \} \; .
\end{equation}

Note that 
${^n \Pi}^s_{D, \mathbf{D}, V, 
\mathbf{V}} \subset
{^n \Pi}^s_{D, \mathbf{D}, V, \mathbf{v}, 
\Tilde{\mathbf{v}}}$, where
$\mathbf{v}^k = 
\dim \mathbf{V}^k -
\dim \Tilde{\mathbf{V}}^k$,
$\Tilde{\mathbf{v}}^k = 
\dim \Tilde{\mathbf{V}}^k -
\dim \mathbf{V}^{k-1}$.

The proof of the following 
proposition is analogous to
the proof of 
Proposition \ref{MdvvSection}.
\begin{proposition}The sets
${^n \Pi}^s_{D, \mathbf{D}, V, 
\mathbf{v}}$, 
${^n \Pi}^s_{D, \mathbf{D}, V, 
\mathbf{v}, \Tilde{\mathbf{v}}}$, and
${^n \Pi}^s_{D, \mathbf{D}, V, 
\mathbf{V}}$,
are locally closed in 
$\Lambda_{D, V}$.
\end{proposition}

Let ${^n \Pi}^{s, \ast s}_{D, 
\mathbf{D}, V, \mathbf{v}}=
{^n \Pi}^s_{D, \mathbf{D}, V, 
\mathbf{v}} \cap
\Lambda^{s, \ast s}_{D, V}$
(a locally closed subset of 
$\Lambda_{D, V}$).
The sets 
${^n \Pi}^s_{D, \mathbf{D}, V, 
\mathbf{v}}$ and
${^n \Pi}^{s, \ast s}_{D, 
\mathbf{D}, V, 
\mathbf{v}}$ are invariant 
under the action of 
$G_V$, and, moreover, this
action is free (because of
\ref{StableFreeG}). Hence
${^n \mathfrak{T}}_{D, \mathbf{D}, 
V, \mathbf{v}} =
{^n \Pi}^s_{D, \mathbf{D}, V, 
\mathbf{v}} /G_V$ and
${^n \mathfrak{S}}_{D, \mathbf{D}, 
V, \mathbf{v}} =
{^n \Pi}^{s, \ast s}_{D, \mathbf{D}, V, 
\mathbf{v}} /G_V$
are naturally quasi-projective
varieties.

The variety 
${^n \mathfrak{T}}_{D, \mathbf{D}, 
V, \mathbf{v}}$ is called 
\emph{tensor product variety}
because of its role 
in the geometric description of the  
tensor product of finite dimensional
representations of $\mathfrak{g}$.
In particular it is shown below
(cf. \ref{CorollaryTensorVariety})
that the set of irreducible components
of ${^n \mathfrak{T}}_{D, \mathbf{D}, 
V, \mathbf{v}}$ is in a bijection with
a weight subset of 
the crystal of the canonical basis
of a product of $n$ representations
of $\mathfrak{g}$.

The variety 
${^n \mathfrak{S}}_{D, \mathbf{D}, 
V, \mathbf{v}}$ is called 
\emph{multiplicity variety}
because as proven below 
(cf. \ref{ADEHallCorollary}) 
the number of its irreducible components
is equal to the multiplicity of
a certain representation of
$\mathfrak{g}$ in a tensor product
of $n$ representations.

The varieties 
${^n \mathfrak{T}}_{D, \mathbf{D}, 
V, \mathbf{v}}$ and 
${^n \mathfrak{S}}_{D, \mathbf{D}, 
V, \mathbf{v}}$
do not depend on the choice of
$V$ as long as $\dim V = v$
is fixed, and they do not
depend up to a non canonical
isomorphism on $D$ and $\mathbf{D}$.
As $D$ and $\mathbf{D}$ never vary in
this paper, a simpler notation is used:
${^n \mathfrak{T}} (d, \mathbf{d},
v, \mathbf{v})=
{^n \mathfrak{T}}_{D, \mathbf{D}, 
V, \mathbf{v}}$ and 
${^n \mathfrak{S}} (d, \mathbf{d}, 
v, \mathbf{v})=
{^n \mathfrak{S}}_{D, \mathbf{D}, 
V, \mathbf{v}}$, where 
$d = \dim D$, 
$\mathbf{d}=\dim \mathbf{D}$ 
(cf. the end of
subsection \ref{Stability}).

\subsection{Inductive construction of 
tensor product varieties}
\label{InductionOfPi}
Let $\mathbf{V} =
( 0 = \mathbf{V}^0
\subset
\Tilde{\mathbf{V}}^1
\subset
\mathbf{V}^1
\subset
\Tilde{\mathbf{V}}^2
\subset
\ldots
\subset
\Tilde{\mathbf{V}}^n
\subset
\mathbf{V}^n = V )$ be
a $2n$-step $\mathbb{Z} [I]$-graded
partial flag in $V$,
and 
$\mathbf{D} =
( 0 = \mathbf{D}^0
\subset
\mathbf{D}^1
\subset
\ldots
\subset
\mathbf{D}^n = D )$ be
an $n$-step $\mathbb{Z} [I]$-graded
partial flag in $D$. Assume $n \geq 2$
and choose an integer $k$ such that
$0 < k < n$. Let $U$ (resp. $C$) be a 
$\mathbb{Z} [I]$-graded subspace in 
$V$ (resp. in $D$) 
complimentary to $\mathbf{V}^k$
(resp. $\mathbf{D}^k$).
Let $\mathbf{V}' =
( 0 = \mathbf{V}^0
\subset
\Tilde{\mathbf{V}}^1
\subset
\mathbf{V}^1
\subset
\Tilde{\mathbf{V}}^2
\subset
\ldots
\subset
\Tilde{\mathbf{V}}^k
\subset
\mathbf{V}^k  )$ 
(resp. $\mathbf{D}' =
( 0 = \mathbf{D}^0
\subset
\mathbf{D}^1
\subset
\ldots
\subset
\mathbf{D}^k  )$)
be
a $2k$-step $\mathbb{Z} [I]$-graded
partial flag in $\mathbf{V}^k$
(resp. a $k$-step 
$\mathbb{Z} [I]$-graded
partial flag in $\mathbf{D}^k$),
obtained by restricting the flag
$\mathbf{V}$ (resp. $\mathbf{D}$).
Similarly, let
$\mathbf{U} =
( 0 = \mathbf{U}^0
\subset
\Tilde{\mathbf{U}}^1
\subset
\mathbf{U}^1
\subset
\Tilde{\mathbf{U}}^2
\subset
\ldots
\subset
\Tilde{\mathbf{U}}^{n-k}
\subset
\mathbf{U}^{n-k} = U)$ 
(resp. $\mathbf{C} =
( 0 = \mathbf{D}^0
\subset
\mathbf{D}^1
\subset
\ldots
\subset
\mathbf{D}^{n-k} = D)$)
be
a $2(n-k)$-step $\mathbb{Z} [I]$-graded
partial flag in $U$
(resp. an $(n-k)$-step 
$\mathbb{Z} [I]$-graded
partial flag in $C$),
obtained by taking intersections of
the subspaces of the flag $\mathbf{V}$ 
(resp. $\mathbf{D}$) with
$U$ (resp. $C$).

One has a regular map
\begin{equation}\nonumber
\rho_2: {^n \Pi}^s_{D, \mathbf{D}, V, 
\mathbf{V}} 
\rightarrow
{^ k\Pi}^s_{\mathbf{D}^k, \mathbf{D}', 
\mathbf{V}^k, \mathbf{V}'}
\times
{^{n-k} \Pi}^s_{C, \mathbf{C}, U, 
\mathbf{U}}
\end{equation}
given by
\begin{equation}\nonumber
\rho_2 (x,p,q) =(
(x^{\sss \mathbf{V}^k \mathbf{V}^k},  
p^{\sss \mathbf{V}^k \mathbf{D}^k} ,
q^{\sss \mathbf{D}^k \mathbf{V}^k}),
(x^{\sss UU}, p^{\sss UC} ,
q^{\sss CU})) \; ,
\end{equation}
where for example $p^{\sss UC}$ is the
$C \rightarrow U$ component of the
block-matrix 
$p: D=\mathbf{D}^k \oplus C 
\rightarrow \mathbf{V}^k \oplus U =V$, 
and the other maps in the RHS are 
defined similarly.

\begin{proposition}
The map $\rho_2$ is a vector bundle with
fibers of dimension
\begin{equation}\nonumber
< X u , v-u > +
< c , v-u > +
< d-c , u > -
< u , v-u > \; ,
\end{equation}
where $v = \dim V$, $d = \dim D$, 
$u = \dim U$, $c = \dim C$.
\end{proposition}
\begin{proof}
The fiber of the map $\rho_2$ over a point
$((x^{\sss \mathbf{V}^k \mathbf{V}^k},  
p^{\sss \mathbf{V}^k \mathbf{D}^k} ,
q^{\sss \mathbf{D}^k \mathbf{V}^k}),
(x^{\sss UU}, 
p^{\sss UC} , q^{\sss CU}))$
in 
${^k \Pi}^s_{\mathbf{D}^k, \mathbf{D}', 
\mathbf{V}^k, \mathbf{V}'}
\times
{^{n-k} \Pi}^s_{C, \mathbf{C}, U, 
\mathbf{U}}$
consists of all linear maps
\begin{equation}\nonumber
\begin{split}
x_h^{\sss \mathbf{V}^k U} 
&\in
\Hom_{\mathbb{C}}
(U_{\Out (h)}, 
\mathbf{V}^k_{\In (h)}) \; ,
\\
p_i^{\sss \mathbf{V}^k C} 
&\in
\Hom_{\mathbb{C}}
(C_i, \mathbf{V}^k_i ) \; ,
\\
q_i^{\sss \mathbf{D}^k U} 
&\in
\Hom_{\mathbb{C}}
( U_i , \mathbf{D}^k_i ) \; ,
\end{split}
\end{equation}
(where $i \in I$, $h \in H$) 
subject to the following 
linear equations
\begin{equation}\nonumber
\sum_{\substack{
h \in H \\ \In (h) = i }}
\varepsilon (h)
( x^{\sss \mathbf{V}^k \mathbf{V}^k}_h
x^{\sss \mathbf{V}^k U}_{\Bar{h}} +
x^{\sss \mathbf{V}^k U}_h
x^{\sss U U}_{\Bar{h}}) - 
p^{\sss \mathbf{V}^k C}_i q^{\sss CU}_i -
p^{\sss \mathbf{V}^k \mathbf{D}^k}_i 
q^{\sss \mathbf{D}^k U}_i =0
\end{equation}
for any $i \in I$. Now the proposition 
follows from the fact that the linear map
\begin{multline}\label{DimPiSurjMap}
(\oplus_{h \in H}
\Hom_{\mathbb{C}}
(U_{\Out (h)},
\mathbf{V}^k_{\In (h)}))
\oplus
(\oplus_{i \in I} 
\Hom_{\mathbb{C}}
(C_i, \mathbf{V}^k_i )) 
\oplus
(\oplus_{i \in I}
\Hom_{\mathbb{C}}
( U_i , \mathbf{D}^k_i ))
\rightarrow \\ \rightarrow
\oplus_{i \in I}
\Hom_{\mathbb{C}}
(U_i , \mathbf{V}^k_i ) 
\end{multline}
given by
\begin{equation}\nonumber
\{ x^{\sss \mathbf{V}^k U} ,
p^{\sss \mathbf{V}^k C} ,
q^{\sss \mathbf{D}^k U} \}
\rightarrow
\sum_{\substack{
h \in H \\ \In (h) = i }}
\varepsilon (h)
( x^{\sss \mathbf{V}^k \mathbf{V}^k}_h
x^{\sss \mathbf{V}^k U}_{\Bar{h}} +
x^{\sss \mathbf{V}^k U}_h
x^{\sss U U}_{\Bar{h}}) - 
p^{\sss \mathbf{V}^k C}_i q^{\sss CU}_i -
p^{\sss \mathbf{V}^k \mathbf{D}^k}_i 
q^{\sss \mathbf{D}^k U}_i 
\end{equation}
is surjective. To prove this let
$a \in \oplus_{i \in I}
\Hom_{\mathbb{C}}
(\mathbf{V}^k_i , U_i)$ be orthogonal
to the image of the above map with 
respect to the $\tr$ pairing, that is
\begin{equation}\nonumber
\sum_{\substack{
h \in H \\ \In (h) = i}}
\varepsilon (h)
\tr ( x^{\sss \mathbf{V}^k \mathbf{V}^k}_h
x^{\sss \mathbf{V}^k U}_{\Bar{h}} a_i +
x^{\sss \mathbf{V}^k U}_h
x^{\sss U U}_{\Bar{h}} a_i) - 
\tr p^{\sss \mathbf{V}^k C}_i 
q^{\sss CU}_i a_i -
\tr p^{\sss \mathbf{V}^k \mathbf{D}^k}_i 
q^{\sss \mathbf{D}^k U}_i a_i = 0
\end{equation}
for any
$x^{\sss \mathbf{V}^k U}$,
$p^{\sss \mathbf{V}^k C}$,
$q^{\sss \mathbf{D}^k U}$,
and any $i \in I$.
It follows that 
\begin{equation}\nonumber
\begin{split}
a_{\In (h)} 
x^{\sss \mathbf{V}^k \mathbf{V}^k}_h &=
x^{\sss U U}_h a_{\Out (h)} \; ,
\\ 
q^{\sss CU}_i a_i &= 0 \; ,
\\
a_i p^{\sss \mathbf{V}^k 
\mathbf{D}^k}_i &= 0 \; ,
\end{split}
\end{equation}
for any $h \in H$, $i \in I$. In particular
the kernel of $a$ contains the image of 
$p^{\sss \mathbf{V}^k \mathbf{D}^k}$ and
is $x^{\sss \mathbf{V}^k 
\mathbf{V}^k}$-invariant. Since
$(x^{\sss \mathbf{V}^k \mathbf{V}^k} ,
p^{\sss \mathbf{V}^k \mathbf{D}^k},
q^{\sss \mathbf{D}^k \mathbf{V}^k})$
is stable it follows that $a = 0$ and
hence the map \eqref{DimPiSurjMap} is
surjective.
The proposition is proven.
\end{proof}

\subsection{Dimensions of tensor
product and multiplicity varieties}
\label{Dimensions}

The following proposition
gives the dimensions of various
varieties defined above.

\begin{proposition}Let $d = \dim D$, 
$v = \dim V$. Then
\begin{alphenum}
\item\label{DimPiDDVV}
${^n \Pi}^s_{D, \mathbf{D}, V, \mathbf{V}}$ 
is empty or has pure dimension
\begin{multline}\nonumber
\dim {^n \Pi}^s_{D, \mathbf{D}, V, 
\mathbf{V}} =
\frac{1}{2} <X v, v> + <d , v> -
\frac{1}{2} <v , v> 
+ \\ +
\sum_{s=1}^n \biggl(
\frac{1}{2} <X \mathbf{v}^s, \mathbf{v}^s> + 
<\mathbf{d}^s , \mathbf{v}^s > -
\frac{1}{2} <\mathbf{v}^s , \mathbf{v}^s> +
\frac{1}{2} <\Tilde{\mathbf{v}}^s , 
\Tilde{\mathbf{v}}^s> 
\biggr) \; ,
\end{multline}
where $\mathbf{v}^k = 
\dim \mathbf{V}^k -
\dim \Tilde{\mathbf{V}}^k$,
$\Tilde{\mathbf{v}}^k = 
\dim \Tilde{\mathbf{V}}^k -
\dim \mathbf{V}^{k-1}$;
\item\label{DimPiDDVvv}
${^n \Pi}^s_{D, \mathbf{D}, V, \mathbf{v}, 
\Tilde{\mathbf{v}}}$ 
is empty or has pure dimension
\begin{multline}\nonumber
\dim 
{^n \Pi}^s_{D, \mathbf{D}, V, \mathbf{v}, 
\Tilde{\mathbf{v}}} =
\frac{1}{2} <X v, v> + <d , v>  
+ \\ +
\sum_{s=1}^n \biggl(
\frac{1}{2} <X \mathbf{v}^s, \mathbf{v}^s> + 
<\mathbf{d}^s , \mathbf{v}^s > -
<\mathbf{v}^s , \mathbf{v}^s>
\biggr) \; ;
\end{multline}
\item\label{DimPiDDVv}
${^n \Pi}^s_{D, \mathbf{D}, V, \mathbf{v}}$ 
is empty or has pure dimension
\begin{multline}\nonumber
\dim {^n \Pi}^s_{D, \mathbf{D}, V, \mathbf{v}} =
\frac{1}{2} <X v, v> + <d , v>  
+ \\ +
\sum_{s=1}^n \biggl(
\frac{1}{2} <X \mathbf{v}^s, 
\mathbf{v}^s> + 
<\mathbf{d}^s , \mathbf{v}^s > -
<\mathbf{v}^s , \mathbf{v}^s>
\biggr) \; ;
\end{multline}
\item\label{DimPiDDVvSS}
${^n \Pi}^{s, \ast s}_{D, \mathbf{D}, V, 
\mathbf{v}}$ 
is empty or has pure dimension
\begin{multline}\nonumber
\dim {^n \Pi}^{s, \ast s}_{D, \mathbf{D}, V, 
\mathbf{v}} =
\frac{1}{2} <X v, v> + <d , v>  
+ \\ +
\sum_{s=1}^n \biggl(
\frac{1}{2} <X \mathbf{v}^s, 
\mathbf{v}^s> + 
<\mathbf{d}^s , \mathbf{v}^s > -
<\mathbf{v}^s , \mathbf{v}^s>
\biggr) \; ;
\end{multline}
\item\label{DimTddvv}
the tensor product variety
${^n \mathfrak{T}} (d, \mathbf{d}, v, 
\mathbf{v})$ 
is empty or has pure dimension
\begin{multline}\nonumber
\dim {^n \mathfrak{T}} (d, \mathbf{d}, v, 
\mathbf{v}) =
\frac{1}{2} <X v, v> + <d , v> -
<v , v> 
+ \\ +
\sum_{s=1}^n \biggl(
\frac{1}{2} <X \mathbf{v}^s, 
\mathbf{v}^s> + 
<\mathbf{d}^s , \mathbf{v}^s > -
<\mathbf{v}^s , \mathbf{v}^s> \biggr)
= \\ =
\frac{1}{2} \biggl(
\dim \mathfrak{M}^s (d, v) +
\sum_{s = 1}^n \dim 
\mathfrak{M}^{s , \ast s} 
(\mathbf{d}^s, \mathbf{v}^s) 
\biggr) \; ;
\end{multline}
\item\label{DimSddvv}
the multiplicity variety
${^n \mathfrak{S}} (d, \mathbf{d}, v, 
\mathbf{v})$ 
is empty or has pure dimension
\begin{multline}\nonumber
\dim {^n \mathfrak{S}} (d, \mathbf{d}, v, 
\mathbf{v}) =
\frac{1}{2} <X v, v> + <d , v> -
<v , v> 
+ \\ +
\sum_{s=1}^n \biggl(
\frac{1}{2} <X \mathbf{v}^s, 
\mathbf{v}^s> + 
<\mathbf{d}^s , \mathbf{v}^s > -
<\mathbf{v}^s , \mathbf{v}^s> \biggr)
= \\ =
\frac{1}{2} \biggl(
\dim \mathfrak{M}^{s, \ast s} (d, v) +
\sum_{s = 1}^n \dim 
\mathfrak{M}^{s , \ast s} 
(\mathbf{d}^s, \mathbf{v}^s) 
\biggr) \; ;
\end{multline}
\end{alphenum}
\end{proposition}
\begin{proof}
\ref{DimPiDDVV} follows by induction
in $n$ using Proposition
\ref{InductionOfPi}. The 
base of the induction is provided by
\ref{LambdaDVUDimBig}.

\ref{DimPiDDVvv} follows from 
\ref{DimPiDDVV} using the fact that 
${^n \Pi}^s_{D, \mathbf{D}, V, \mathbf{v}, 
\Tilde{\mathbf{v}}}$ is a fibration with
fibers isomorphic to 
${^n \Pi}^s_{D, \mathbf{D}, V, \mathbf{V}}$
over the variety of graded
$2n$-step partial flags in $V$
with dimensions of the subfactors given
by $\mathbf{v}$ and $\Tilde{\mathbf{v}}$.
The dimension of this flag variety
is equal to
\begin{equation}\nonumber
\frac{1}{2} \biggl(
<v , v> -
\sum_{s = 1}^n \bigl(
<\mathbf{v}^s , \mathbf{v}^s> +
<\Tilde{\mathbf{v}}^s , 
\Tilde{\mathbf{v}}^s>
\bigr) \biggr) \: .
\end{equation}

\ref{DimPiDDVv} follows from the fact
that ${^n \Pi}^s_{D, \mathbf{D}, 
V, \mathbf{v}}$ 
is a union of a finite number of locally 
closed subsets ${^n \Pi}^s_{D, \mathbf{D}, V, 
\mathbf{v}, \Tilde{\mathbf{v}}}$ 
(for different $\Tilde{\mathbf{v}}$) 
having identical dimensions 
(cf. \ref{DimPiDDVvv}).

\ref{DimPiDDVvSS} follows from the
fact that 
${^n \Pi}^{s, \ast s}_{D, \mathbf{D}, V, 
\mathbf{v}}$ is open in  
${^n \Pi}^s_{D, \mathbf{D}, V, \mathbf{v}}$
and from \ref{DimPiDDVv}.

To prove \ref{DimTddvv} (resp. \ref{DimSddvv})
note that 
${^n \mathfrak{T}} (d, \mathbf{d}, v, 
\mathbf{v}) = 
{^n \Pi}^s_{D, \mathbf{D}, V, 
\mathbf{v}} /
G_V$
(resp.
${^n \mathfrak{S}} (d, \mathbf{d}, v, 
\mathbf{v}) = 
{^n \Pi}^{s, \ast s}_{D, \mathbf{D}, 
V, \mathbf{v}} / G_V$)
and the action of the group $G_V$ on 
${^n \Pi}^s_{D, \mathbf{D}, 
V, \mathbf{v}}$ (resp.
${^n \Pi}^{s, \ast s}_{D, \mathbf{D}, 
V, \mathbf{v}}$) is free. Now
\ref{DimTddvv} (resp. \ref{DimSddvv})
follows from 
\ref{DimPiDDVv} (resp.
\ref{DimPiDDVvSS}).
\end{proof}

\subsection{Irreducible components
of tensor product varieties}
\label{IrreducibleComponentsOfTS}

Let ${^n \mathcal{T}} 
(d, \mathbf{d}, v, \mathbf{v})$ 
(resp. ${^n \mathcal{S}} 
(d , \mathbf{d}, v, \mathbf{v})$)
denote the set of irreducible 
components of the tensor product 
variety ${^n \mathfrak{T}} 
(d, \mathbf{d}, v, \mathbf{v})$ 
(resp. the multiplicity
variety ${^n \mathfrak{S}} 
(d , \mathbf{d}, v, \mathbf{v})$).

Since ${^n \Pi}^s_{D,\mathbf{D}, V, 
\mathbf{v}}$ is
the total space of a principal 
$G_V$-bundle over
${^n \mathfrak{T}} 
(d, \mathbf{d}, v, \mathbf{v})$,
the set of irreducible components 
of ${^n \Pi}^s_{D, \mathbf{D}, V, 
\mathbf{v}}$
can be naturally identified with
${^n \mathcal{T}}
(d, \mathbf{d}, v, \mathbf{v})$.
On the other hand,
${^n \Pi}^s_{D, \mathbf{D}, V ,
\mathbf{v}}$
is a union of locally closed subsets
${^n \Pi}^s_{D, \mathbf{D}, V, 
\mathbf{v}, \Tilde{\mathbf{v}}}$ 
(for different
$\Tilde{\mathbf{v}}$) having the
same dimension (cf. Proposition
\ref{Dimensions}). Hence
\begin{equation}\nonumber
{^n \mathcal{T}}
(d, \mathbf{d}, v, \mathbf{v}) =
\bigsqcup_{\Tilde{\mathbf{v}}}
{^n \mathcal{T}}
(d, \mathbf{d}, v, \mathbf{v}, 
\Tilde{\mathbf{v}}) \; ,
\end{equation}
where
$\Tilde{\mathbf{v}}$ ranges over
all $n$-tuples of elements of
$\mathbb{Z} [I]$ such that
$\sum_{s=1}^n \mathbf{v}^s + 
\sum_{s=1}^n \Tilde{\mathbf{v}}^s =
v$, and 
${^n \mathcal{T}}
(d, \mathbf{d}, v, \mathbf{v}, 
\Tilde{\mathbf{v}})$ denotes the set
of irreducible components of
${^n \Pi}^s_{D, \mathbf{D}, V,
\mathbf{v}, \Tilde{\mathbf{v}}}$.
The variety ${^n \Pi}^s_{D, \mathbf{D},
V, \mathbf{v}, \Tilde{\mathbf{v}}}$ 
is a locally trivial fibration
over the (simply connected
$G_V$-homogeneous) 
variety of all graded $2n$-step 
partial flags in 
$V$ with dimensions of
the subfactors given by 
$\mathbf{v}$ and $\Tilde{\mathbf{v}}$
(as in \ref{DefinitionOfTS}).
A fiber of this fibration is
isomorphic to 
${^n \Pi}^s_{D , \mathbf{D}, 
V, \mathbf{V}}$ , where
$\mathbf{V}$ is a point of the base
(a $2n$-step partial flag).
It follows that the set
of irreducible components of
${^n \Pi}^s_{D , \mathbf{D}, 
V, \mathbf{V}}$
can be naturally identified with 
the set
${^n \mathcal{T}}
(d, \mathbf{d}, v, \mathbf{v}, 
\Tilde{\mathbf{v}})$.
To summarize,
${^n \mathcal{T}}
(d, \mathbf{d}, v, \mathbf{v})$
is the set of irreducible
components of 
${^n \Pi}^s_
{D, \mathbf{D}, V, \mathbf{v}}$ and
${^n \mathfrak{T}}
(d, \mathbf{d}, v, \mathbf{v})$,
${^n \mathcal{T}}
(d, \mathbf{d}, v, \mathbf{v}, 
\Tilde{\mathbf{v}})$
is the set of irreducible components
of ${^n \Pi}^s_{D, \mathbf{D}, V,
\mathbf{v}, \Tilde{\mathbf{v}}}$
and
${^n \Pi}^s_{D, \mathbf{D}, 
V, \mathbf{V}}$,
and ${^n \mathcal{T}}
(d, \mathbf{d}, v, \mathbf{v}) =
\bigsqcup_{\Tilde{\mathbf{v}}}
{^n \mathcal{T}}
(d, \mathbf{d}, v, \mathbf{v}, 
\Tilde{\mathbf{v}})$.

Similarly
${^n \mathcal{S}}
(d, \mathbf{d}, v, \mathbf{v})$
is the set of irreducible
components of 
${^n \mathfrak{S}}
(d, \mathbf{d}, v, \mathbf{v})$
and
${^n \Pi}^{s, \ast s}_{D, \mathbf{D}, 
V, \mathbf{v}}$.

Finally, let 
${^n \mathcal{T}}
(d, \mathbf{d}, \mathbf{v}) =
\bigsqcup_{v \in 
\mathbb{Z}_{\geq 0} [I]}
{^n \mathcal{T}}
(d, \mathbf{d}, v, \mathbf{v})$.

\subsection{The first bijection for a
tensor product variety}
\label{DefinitionOfAlpha}

Let $1 < k < n$.
The vector bundle $\rho_2$  
introduced in \ref{InductionOfPi}
induces a bijection of
the sets of irreducible components
\begin{equation}\nonumber
\alpha_{k, n-k} :
{^n \mathcal{T}} 
(d, \mathbf{d}, v, \mathbf{v},
\Tilde{\mathbf{v}}) 
\xrightarrow{\sim}
{^k \mathcal{T}} (d', \mathbf{d}', 
v', \mathbf{v}', \Tilde{\mathbf{v}}') 
\times
{^{n-k}}\mathcal{T} (d'', \mathbf{d}'', 
v'', \mathbf{v}'', 
\Tilde{\mathbf{v}}'') \; ,
\end{equation}
where $d$, $d'$, $d''$, $v$, $v'$, 
$v'' \in \mathbb{Z}_{\geq 0} [I]$,
$d = d' +d''$, $v=v'+v''$,
$\mathbf{d}$,
$\mathbf{v}$, $\Tilde{\mathbf{v}}$
are $n$-tuples of elements
of $\mathbb{Z}_{\geq 0} [I]$,
$\mathbf{d}'$,
$\mathbf{v}'$, $\Tilde{\mathbf{v}}'$
are $k$-tuples of elements
of $\mathbb{Z}_{\geq 0} [I]$,
$\mathbf{d}''$,
$\mathbf{v}''$, $\Tilde{\mathbf{v}}''$
are $(n-k)$-tuples of elements
of $\mathbb{Z}_{\geq 0} [I]$,
$\mathbf{d} = (\mathbf{d}' , 
\mathbf{d}'')$,
$\mathbf{v} = (\mathbf{v}' , 
\mathbf{v}'')$,
$\Tilde{\mathbf{v}} = (\Tilde{\mathbf{v}}', 
\Tilde{\mathbf{v}}'')$,
$v' = \sum_{s=1}^k (\mathbf{v}^{\prime s} +
\Tilde{\mathbf{v}}^{\prime s})$,
$v'' = \sum_{s=1}^{n-k} 
(\mathbf{v}^{\prime \prime s} +
\Tilde{\mathbf{v}}^{\prime \prime s})$.
Taking union over $\Tilde{\mathbf{v}}$ one
obtains a bijection (denoted again by
$\alpha_{k , n-k}$)
\begin{equation}\nonumber
\alpha_{k, n-k} :
{^n \mathcal{T}} 
(d, \mathbf{d}, v, \mathbf{v}) 
\xrightarrow{\sim}
\bigsqcup_{\substack{v' , v''
\in \mathbb{Z}_{\geq 0} [I] \\
v' + v'' = v}}
{^k \mathcal{T}} (d', \mathbf{d}', 
v', \mathbf{v}') 
\times
{^{n-k} \mathcal{T}} (d'', 
\mathbf{d}'', v'', \mathbf{v}'') \; .
\end{equation}

Generalizing
the map $\rho_2$ (cf. \ref{InductionOfPi})
one can consider a map
(given by restriction of $x$, $p$, and
$q$) 
\begin{multline}\nonumber
\rho_3 : \quad
{^n \Pi}^s_{D, \mathbf{D}, V,
\mathbf{V}}
\rightarrow
\\
\rightarrow
{^k \Pi}^s_{\mathbf{D}^k, \mathbf{D}', 
\mathbf{V}^k, \mathbf{V}'}
\times
{^{l-k} \Pi}^s_{\mathbf{D}^l/\mathbf{D}^k, 
\mathbf{D}'', 
\mathbf{V}^l/\mathbf{V}^k, \mathbf{V}''}
\times
{^{n-l} \Pi}^s_{D/\mathbf{D}^l, 
\mathbf{D}''', 
V/\mathbf{V}^l, \mathbf{V}'''} \; ,
\end{multline}
where $1 < k < l < n$, 
$\mathbf{D}'$ (resp. $\mathbf{D}''$,
$\mathbf{D}'''$) is the $k$-step
partial flag in $\mathbf{D}^k$
(resp. $(l-k)$-step partial flag in
$\mathbf{D}^l/\mathbf{D}^k$,
$(n-l)$-step partial flag in
$D/\mathbf{D}^l$) induced by 
the $n$-step flag $\mathbf{D}$ in $D$, 
and, similarly,
$\mathbf{V}'$ (resp. $\mathbf{V}''$,
$\mathbf{V}'''$) is the $2k$-step
partial flag in $\mathbf{V}^k$
(resp. $2(l-k)$-step partial flag in
$\mathbf{V}^l/\mathbf{V}^k$,
$2(n-l)$-step partial flag in
$V/\mathbf{V}^l$) induced by 
the $2n$-step flag $\mathbf{V}$ in $V$.

The map $\rho_3$ can be represented
in two ways
as a composition of two maps $\rho_2$:
\begin{equation}\label{Rho3Diagram}
\xymatrix{
*\txt{${^n \Pi}^s_{D, \mathbf{D}, V,
\mathbf{V}}$}
\ar[ddrr]_-{\rho_3}
\ar[rr]^-{\rho_2}
\ar[dd]^-{\rho_2}
&&
*\txt{${^l \Pi}^s_{\mathbf{D}^l, 
\Bar{\Bar{\mathbf{D}}}, 
\mathbf{V}^l, \Bar{\Bar{\mathbf{V}}}}
\quad \times$ \\ $
{^{n-l} \Pi}^s_{D/\mathbf{D}^l, 
\mathbf{D}''', 
V/\mathbf{V}^l, \mathbf{V}'''}$}
\ar[dd]^-{\rho_2 \times \Id}
\\ \\
*\txt{${^k \Pi}^s_{\mathbf{D}^k, 
\mathbf{D}', 
\mathbf{V}^k, \mathbf{V}'}
\quad \times$ \\ $
{^{n-k} \Pi}^s_{D/\mathbf{D}^k, 
\Bar{\mathbf{D}}, 
V/\mathbf{V}^k, \Bar{\mathbf{V}}}$}
\ar[rr]_-{\Id \times \rho_2}
&&
*\txt{${^k \Pi}^s_{\mathbf{D}^k, 
\mathbf{D}', 
\mathbf{V}^k, \mathbf{V}'}
\quad \times$ \\ $
{^{l-k} \Pi}^s_{\mathbf{D}^l/\mathbf{D}^k, 
\mathbf{D}'', 
\mathbf{V}^l/\mathbf{V}^k, \mathbf{V}''}
\quad \times$ \\ $
{^{n-l} \Pi}^s_{D/\mathbf{D}^l, 
\mathbf{D}''', 
V/\mathbf{V}^l, \mathbf{V}'''}$}
}
\end{equation}
where  
$\Bar{\mathbf{D}}$ (resp. 
$\Bar{\Bar{\mathbf{D}}}$) 
is the partial flag in $D/\mathbf{D}^k$
(resp. partial flag in
$\mathbf{D}^l$) induced by 
the flag $\mathbf{D}$ in $D$, 
and, similarly,
$\Bar{\mathbf{V}}$ (resp. 
$\Bar{\Bar{\mathbf{V}}}$) 
is the partial flag in $V/\mathbf{V}^k$
(resp. partial flag in
$\mathbf{V}^l$) induced by 
the flag $\mathbf{V}$ in $V$.
It follows from the commutativity
of the diagram \eqref{Rho3Diagram} and
Proposition \ref{InductionOfPi}
that $\rho_3$ is
a locally trivial fibration with 
a smooth connected fiber, and therefore,
induces a bijection
\begin{multline}\nonumber
\alpha_{k, l-k, n-l} : \quad
{^n \mathcal{T}} (
(d' +d''+d'''), 
(\mathbf{d}' , 
\mathbf{d}'', \mathbf{d}'''),
v,
(\mathbf{v}' , 
\mathbf{v}'', \mathbf{v}''')) 
\xrightarrow{\sim}
\\
\xrightarrow{\sim}
\bigsqcup_{\substack{v' , v'', v'''
\in \mathbb{Z}_{\geq 0} [I] \\
v' + v'' + v''' = v}}
{^k \mathcal{T}} (d', \mathbf{d}', 
v', \mathbf{v}') 
\times
{^{l-k} \mathcal{T}} (d'', \mathbf{d}'', 
v'', \mathbf{v}'')
\times
{^{n-l} \mathcal{T}} (d''', \mathbf{d}''', 
v''', \mathbf{v}''')  
\end{multline}
Moreover, the commutativity of the
diagram \eqref{Rho3Diagram} implies
that
\begin{equation}\nonumber
\alpha_{k , l-k , n-l}=
(\alpha_{k , l-k} \times \Id) \circ
\alpha_{l , n-l} =
(\Id \times \alpha_{l-k , n-l}) \circ
\alpha_{k , n-k} \; .
\end{equation}

One can consider analogues of 
the fibrations $\rho_2$ and $\rho_3$
taking values in the product of 
any number ($\leq n$) of  
varieties $\Pi^s$. These maps can be
represented as compositions of the
maps $\rho_2$ in several ways
(cf. \eqref{Rho3Diagram}), hence they
are fibrations with smooth connected 
fibers and induce isomorphisms
of the sets of irreducible components.
The most important case is the map
$\rho_n$ that takes value in the product
of $n$ varieties ${^1 \Pi}^s$
(recall that a variety ${^1 \Pi}^s$ is
the same as the variety $\Lambda^s$). 
More precisely, let
$\Hat{\mathbf{D}}^s$
(resp. $\Hat{\mathbf{V}}^s$) be
a $\mathbb{Z} [I]$-graded subspace in 
$\mathbf{D}^s$ complimentary 
to $\mathbf{D}^{s-1}$
(resp. 
a $\mathbb{Z} [I]$-graded subspace in 
$\mathbf{V}^s$ complimentary 
to $\mathbf{V}^{s-1}$)
and let
\begin{equation}\nonumber
\rho_n : \quad
{^n \Pi}^s_{D, \mathbf{D},
V, \mathbf{V}} 
\rightarrow
\Lambda^s_{\Hat{\mathbf{D}}^1 , 
\Hat{\mathbf{V}}^1} 
\times \ldots \times
\Lambda^s_{\Hat{\mathbf{D}}^n , 
\Hat{\mathbf{V}}^n}
\end{equation}
be a regular map given by
\begin{equation}\nonumber
\rho_n (x, p, q) =
((x^{\sss \Hat{\mathbf{V}}^1 
\Hat{\mathbf{V}}^1},
p^{\sss \Hat{\mathbf{V}}^1 
\Hat{\mathbf{D}}^1},
q^{\sss \Hat{\mathbf{D}}^1 
\Hat{\mathbf{V}}^1}),
\ldots ,
(x^{\sss \Hat{\mathbf{V}}^n 
\Hat{\mathbf{V}}^n},
p^{\sss \Hat{\mathbf{V}}^n 
\Hat{\mathbf{D}}^n},
q^{\sss \Hat{\mathbf{D}}^n 
\Hat{\mathbf{V}}^n})) \; .
\end{equation}
The map $\rho_n$ can be represented
(in many ways) as a composition of
maps $\rho_2$. Therefore it is
a fibration with smooth connected
fibers and it induces a bijection
\begin{equation}\nonumber
\alpha_n : \quad
{^n \mathcal{T}} (d, \mathbf{d},
v, \mathbf{v}, \Tilde{\mathbf{v}})
\xrightarrow{\sim}
\mathcal{M} (\mathbf{d}^1, \mathbf{v}^1,
\mathbf{v}^1 + \Tilde{\mathbf{v}}^1)
\times \ldots \times
\mathcal{M} (\mathbf{d}^n, \mathbf{v}^n,
\mathbf{v}^n + \Tilde{\mathbf{v}}^n)
\end{equation}
or, after taking union 
over $\Tilde{\mathbf{v}}$,
\begin{equation}\nonumber
\alpha_n : \quad
{^n \mathcal{T}} (d, \mathbf{d},
v, \mathbf{v})
\xrightarrow{\sim}
\bigsqcup_{\substack{
\mathbf{u}^1, \ldots , \mathbf{u}^n 
\in \mathbb{Z}_{\geq 0} [I] \\
\mathbf{u}^1 + \ldots + 
\mathbf{u}^n = v }}
\mathcal{M} (\mathbf{d}^1, \mathbf{v}^1,
\mathbf{u}^1)
\times \ldots \times
\mathcal{M} (\mathbf{d}^n, \mathbf{v}^n,
\mathbf{u}^n) \; .
\end{equation}
Finally, union over $v$ gives a bijection
\begin{equation}\label{AlphaN}
\alpha_n : \quad
{^n \mathcal{T}} 
(d, \mathbf{d}, \mathbf{v})
\xrightarrow{\sim}
\mathcal{M} (\mathbf{d}^1, \mathbf{v}^1)
\times \ldots \times
\mathcal{M} (\mathbf{d}^n, 
\mathbf{v}^n) \; .
\end{equation}
Recall that $\mathcal{M} (d,v,u)$ 
denotes the set of irreducible 
components of a certain
locally closed subset of a 
quiver variety 
(cf. \ref{DefinitionOfMCal}).
The bijection $\alpha_n$ justifies 
the name ``tensor product variety''
given to ${^n \mathfrak{T}} 
(d, \mathbf{d}, v, \mathbf{v})$,
since it is known that the set 
$\mathcal{M} (d, v)$ can
be equipped with the structure
of a $\mathfrak{g}$-crystal, which
is isomorphic to the crystal of
the canonical basis of an
irreducible finite dimensional 
representation of $\mathfrak{g}$
(see \ref{DefinitionOfMCal} for
references and \ref{CrystalMDV0}
for a proof).

\subsection{The multiplicity variety 
for two multiples}
\label{MultiplicityTwo}

Let $D = D^1 \oplus D^2$, and
$\mathbf{D} = ( 0 \subset D^1 
\subset D ) \quad$ (a flag in $D$).
Similarly, let $V=V^1 \oplus 
U \oplus V^2$, and
$\mathbf{V} = ( 0 \subset V^1 
\subset V^1 \oplus U \subset V )$.
Then one has the following map
\begin{equation}\label{SigmaPiLLL}
\sigma_2 :
{^2 \Pi}^{s, \ast s}_{D, \mathbf{D}, 
V, \mathbf{V}} \rightarrow
\Lambda^{s, \ast s}_{D^1, V^1} 
\times
\Lambda_{U} 
\times
\Lambda^{s, \ast s}_{D^2, V^2} \; , 
\end{equation}
given by 
\begin{equation}\nonumber
\sigma_2 (x,p,q) =(
(x^{\sss V^1 V^1},  
p^{\sss V^1 D^1} ,
q^{\sss D^1 V^1}),
x^{\sss UU},
(x^{\sss V^2 V^2},  
p^{\sss V^2 D^2} ,
q^{\sss D^2 V^2})) \; .
\end{equation}

\begin{proposition}Let
$v^1 = \dim V^1$, $v^2 = \dim V^2$,
$u = \dim U$, 
$d^1 = \dim D^1$, $d^2 = \dim D^2$,
$v = \dim V = v^1 + v^2 + u$,
$d = \dim D = d^1 + d^2$.
Then:
\begin{alphenum}
\item\label{ImageOfSigma}
the image of $\sigma_2$ is 
$\Lambda^{s, \ast s}_{D^1, V^1} 
\times
\Lambda'_U
\times
\Lambda^{s, \ast s}_{D^2, V^2}$
where $\Lambda'_U$ is an 
open subset of $\Lambda_U$,
\item
$\sigma_2$ is a fibration with smooth
connected fibers of dimension
\begin{multline}\nonumber
\frac{1}{2}(<Xv,v> - <Xv^1,v^1> -
<Xu,u> - <Xv^2,v^2>) + \\ +
<d,u>-<d^1,v^1>-<d^2,v^2> + \\ +
\frac{1}{2}(<v,v> - <v^1,v^1> - 
<u,u> - <v^2,v^2>) \; .
\end{multline}
\end{alphenum}
\end{proposition}
\begin{proof}
The proof is analogous to the proofs of
Propositions 
\ref{DefinitionOfMCal} 
and 
\ref{InductionOfPi}. Namely 
given 
\begin{equation}\nonumber
((x^{\sss V^1 V^1},  
p^{\sss V^1 D^1} ,
q^{\sss D^1 V^1}),
x^{\sss UU},
(x^{\sss V^2 V^2},  
p^{\sss V^2 D^2} ,
q^{\sss D^2 V^2}))
\in
\Lambda^{s, \ast s}_{D^1, V^1} 
\times
\Lambda_{U} 
\times
\Lambda^{s, \ast s}_{D^2, V^2}
\end{equation} 
the fiber of
$\sigma_2$ over this point is 
vector bundle $\sigma'$
over an affine space 
consisting of all linear maps
\begin{equation}\nonumber
x^{\sss V^1 U},
x^{\sss U V^2},
p^{\sss U D^2},
q^{\sss D^1 U}
\end{equation}
subject to the 
equations
\begin{equation}\nonumber
\begin{split}
\sum_{\substack{
h \in H \\ \In (h) = i }}
\varepsilon (h)
( x^{\sss V^1 V^1}_h
x^{\sss V^1 U}_{\Bar{h}} +
x^{\sss V^1 U}_h
x^{\sss U U}_{\Bar{h}}) - 
p^{\sss V^1 D^1}_i 
q^{\sss D^1 U}_i &= 0
\\
\sum_{\substack{
h \in H \\ \In (h) = i }}
\varepsilon (h)
( x^{\sss U U}_h
x^{\sss U V^2}_{\Bar{h}} +
x^{\sss U V^2}_h
x^{\sss V^2 V^2}_{\Bar{h}}) - 
p^{\sss U D^2}_i 
q^{\sss D^2 V^2}_i &= 0 \; .
\end{split} 
\end{equation}

Now \ref{ImageOfSigma}
follows from the condition that
$(x,p,q) \in 
\Pi^{s, \ast s}_{D, \mathbf{D}, 
V, \mathbf{V}}$. 
The proof is similar to
\ref{LambdaDVUOpen} 
(cf. 
\cite[Proposition 1.16]{Lusztig2000a}), 
or one can deduce
\ref{ImageOfSigma} from
\ref{LambdaDVUOpen} and its dual
using the fact that
$(x,p,q) \in 
\Pi^{s, \ast s}_{D, \mathbf{D}, 
V, \mathbf{V}}$ if and only if
\begin{equation}\nonumber
(x^{\sss 
(V^1 \oplus U) (V^1 \oplus U)},
p^{\sss (V^1 \oplus U) D^1},
q^{\sss D^1 (V^1 \oplus U)})
\in 
\Lambda^{\ast s}_{D^1, V^1 \oplus U}
\end{equation}
and
\begin{equation}\nonumber
(x^{\sss 
(U \oplus V^2) (U \oplus V^2)},
p^{\sss (U \oplus V^2) D^2},
q^{\sss D^2 (U \oplus V^2)})
\in 
\Lambda^s_{D^2, U \oplus V^2} \; .
\end{equation}

The fiber of $\sigma'$ over a point
\begin{equation}\nonumber
((x^{\sss V^1 V^1},  
p^{\sss V^1 D^1} ,
q^{\sss D^1 V^1}),
x^{\sss UU},
(x^{\sss V^2 V^2},  
p^{\sss V^2 D^2} ,
q^{\sss D^2 V^2})), 
(x^{\sss V^1 U},
x^{\sss U V^2},
p^{\sss U D^2},
q^{\sss D^1 U}))
\end{equation}
is an affine space of all linear maps
\begin{equation}\nonumber
x^{\sss V^1 V^2},
p^{\sss V^1 D^2},
q^{\sss D^1 V^2}
\end{equation}
subject to the 
equations
\begin{equation}\nonumber
\sum_{\substack{
h \in H \\ \In (h) = i }}
\varepsilon (h)
( x^{\sss V^1 V^1}_h
x^{\sss V^1 V^2}_{\Bar{h}} +
x^{\sss V^1 V^2}_h
x^{\sss V^2 V^2}_{\Bar{h}}) - 
p^{\sss V^1 D^1}_i 
q^{\sss D^1 V^2}_i -
p^{\sss V^1 D^2}_i 
q^{\sss D^2 V^2}_i = 0 \; .
\end{equation}

One has to check that the
systems of linear equations used in the proof
are not overdetermined, which follows from
the fact that
$(x^{\sss V^1 V^1},  
p^{\sss V^1 D^1} ,
q^{\sss D^1 V^1}) \in 
\Lambda^{s, \ast s}_{D^1 , V^1}$,
and
$(x^{\sss V^2 V^2},  
p^{\sss V^2 D^2} ,
q^{\sss D^2 V^2}) \in 
\Lambda^{s, \ast s}_{D^2 , V^2}$
(cf. the proof of 
Proposition \ref{InductionOfPi}).
\end{proof}

Note that because 
$\Lambda^{s, \ast s}_{D,V}$ is smooth
and connected 
(cf. \ref{IrrStableAstStable})
the above Proposition implies
that the set of irreducible components
of ${^2 \Pi}^{s, \ast s}_{
D, \mathbf{D}, V, \mathbf{V}}$ is in
a natural bijection with the set
of irreducible components
of $\Lambda'_{U}$.

\subsection{The first bijection for
a multiplicity variety}
\label{InductionOfMultiplicities}

Let $1 < k < n$.
Consider the variety
${^n \Pi}^{s, \ast s}_{D, \mathbf{D},
V, \mathbf{V}}$. Let
$D^1 = \mathbf{D}^k$, $D^2$ be a 
$\mathbb{Z} [I]$-graded
subspace in $D$ complimentary to
$\mathbf{D}^k$. Similarly let
$V^1 = \mathbf{V}^k$,
$U$ be a $\mathbb{Z} [I]$-graded 
subspace in $\Tilde{\mathbf{V}}^{k+1}$ 
complimentary to $\mathbf{V}^k$, and
$V^2$ be a $\mathbb{Z} [I]$-graded 
subspace in $V$ complimentary to
$\Tilde{\mathbf{V}}^{k+1}$. 
Let $\Bar{\mathbf{D}} = 
( 0 \subset D^1 \subset D ) \quad$ 
(a subflag of $\mathbf{D}$),
$\Bar{\mathbf{V}} = ( 0 \subset V^1 
\subset V^1 \oplus U \subset D ) \quad$
(a subflag of $\mathbf{V}$),
$\mathbf{D}'$ (resp. $\mathbf{D}''$)
be the $k$-step flag in $D^1$ 
(resp. the $(n-k)$-step flag in $D^2$) 
obtained by considering the first $k$ 
subspaces of $\mathbf{D}$ (resp. by taking 
intersections of the last $n-k$ subspaces 
of $\mathbf{D}$ with $D^2$), and
$\mathbf{V}'$ (resp. $\mathbf{V}''$)
be the $2k$-step flag in $V^1$ 
(resp. the $2(n-k)$-step flag in $V^2$) 
obtained by considering the first $2k$ 
subspaces of $\mathbf{V}$ (resp. by taking 
intersections of the last $2(n-k)$ subspaces 
of $\mathbf{V}$ with $V^2$).
Then ${^n \Pi}^{s, \ast s}_{D, \mathbf{D},
V, \mathbf{V}} \subset
{^2 \Pi}^{s, \ast s}_{D, \Bar{\mathbf{D}}, 
V, \Bar{\mathbf{V}}}$,
${^k \Pi}^{s, \ast s}_{D^1, \mathbf{D}',
V^1, \mathbf{V}'} \subset
\Lambda^{s, \ast s}_{D^1, V^1}$,
${^{n-k} \Pi}^{s, \ast s}_{D^2, \mathbf{D}'',
V^2, \mathbf{V}''} \subset
\Lambda^{s, \ast s}_{D^2, V^2}$
and the map $\sigma_2$
(cf. \eqref{SigmaPiLLL}) restricts to a
fibration (with the same fibers
as those of $\sigma_2$) 
\begin{equation}\nonumber
\sigma_{k, n-k} : \quad
{^n \Pi}^{s, \ast s}_{D, \mathbf{D}, 
V, \mathbf{V}} \rightarrow
{^k \Pi}^{s, \ast s}_{D^1, \mathbf{D}',
V^1, \mathbf{V}'} 
\times
\Lambda_{U}' 
\times
{^{n-k} \Pi}^{s, \ast s}_{D^2, \mathbf{D}'', 
V^2 , \mathbf{V}''} \; . 
\end{equation} 
This fibration having smooth
connected fibers it 
induces a bijection
\begin{multline}\nonumber
\eta_{k, n-k} : \quad
{^n \mathcal{S}} (d' +d'', 
(\mathbf{d}', \mathbf{d}''), v, 
(\mathbf{v}', \mathbf{v}''), 
(\Tilde{\mathbf{v}}', 
v-v'-v'', \Tilde{\mathbf{v}}'')) 
\xrightarrow{\sim}
\\
\xrightarrow{\sim}
{^k \mathcal{S}} (d', \mathbf{d}', v', 
\mathbf{v}', \Tilde{\mathbf{v}}')
\times
{^{n-k} \mathcal{S}} (d'', \mathbf{d}'', v'', 
\mathbf{v}'', \Tilde{\mathbf{v}}'')
\times \\ \times
{^2 \mathcal{S}} (d'+d'', (d', d''),
v, (v', v''), (0, v-v'-v'')) \; ,
\end{multline}
where the last multiple in the RHS 
represents
the set of irreducible components of
$\Lambda'$, which is in a natural 
bijection (induced by the fibration
$\sigma_2$ -- cf. \ref{MultiplicityTwo})
with the set of irreducible components of 
${^2 \Pi}^{s, \ast s}_{D, \Bar{\mathbf{D}}, 
V, \Bar{\mathbf{V}}}$.

Taking union over 
$\Tilde{\mathbf{v}}'$, and
$\Tilde{\mathbf{v}}''$
one obtains a bijection
\begin{multline}\nonumber
\eta_{k,n-k} : \quad
{^n \mathcal{S}} (d' +d'', 
(\mathbf{d}', \mathbf{d}''), v, 
(\mathbf{v}', \mathbf{v}'')) 
\xrightarrow{\sim}
\\
\xrightarrow{\sim}
\bigsqcup_{v', v'' 
\in \mathbb{Z}_{\geq 0} [I]}
{^k \mathcal{S}} (d', \mathbf{d}', v', 
\mathbf{v}')
\times 
{^{n-k} \mathcal{S}} (d'', 
\mathbf{d}'', v'', \mathbf{v}'')
\times \\ \times 
{^2 \mathcal{S}} (d'+d'', (d', d''),
v, (v', v'')) \; ,
\end{multline}
which is an analogue of the recurrence 
relation between multiplicities
in tensor products of representations
of $\mathfrak{g}$.

As a generalization of the map
$\sigma_{k,n-k}$ one can consider
a map 
\begin{multline}\label{SigmaN}
\sigma_n : \quad 
{^n \Pi}^{s, \ast s}_{
D, \mathbf{D}, V, \mathbf{V}} 
\rightarrow \\ \rightarrow
\Lambda^{s, \ast s}_{
\mathbf{D}^1/\mathbf{D}^0 , 
\mathbf{V}^1/\Tilde{\mathbf{V}}^1}
\times \ldots \times
\Lambda^{s, \ast s}_{
\mathbf{D}^n/\mathbf{D}^{n-1}, 
\mathbf{V}^n/\Tilde{\mathbf{V}}^n}
\times
\Lambda_{
\Tilde{\mathbf{V}}^1/\mathbf{V}^0}
\times \ldots \times
\Lambda_{
\Tilde{\mathbf{V}}^n/\mathbf{V}^{n-1}} 
\end{multline}
defined by restricting $x$, $p$,
and $q$. The map $\sigma_n$ can
be represented (in many ways) as
a composition of the maps of
type $\sigma_{n,n-k}$ for
different $n$, and $k$
(cf. a similar statement for
the map $\rho_n$ in 
\ref{DefinitionOfAlpha}). Therefore
the image of $\sigma_n$ is
equal to
\begin{equation}\nonumber
\Lambda^{s, \ast s}_{
\mathbf{D}^1/\mathbf{D}^0, 
\mathbf{V}^1/\Tilde{\mathbf{V}}^1}
\times \ldots \times
\Lambda^{s, \ast s}_{
\mathbf{D}^n/\mathbf{D}^{n-1}, 
\mathbf{V}^n/\Tilde{\mathbf{V}}^n}
\times X \; ,
\end{equation}
where $X$ is an open subset of
$\Lambda_{
\Tilde{\mathbf{V}}^1/\mathbf{V}^0}
\times \ldots \times
\Lambda_{
\Tilde{\mathbf{V}}^n/\mathbf{V}^{n-1}}$,
and 
${^n \Pi}^{s, \ast s}_{
D, \mathbf{D}, V, \mathbf{V}}$ is
an open subset of the total space of 
a locally-trivial
fibration $\sigma'$ over the
image of $\sigma$ with a smooth
connected fiber, and such that
$\sigma'$ restricted to 
${^n \Pi}^{s, \ast s}_{
D, \mathbf{D}, V, \mathbf{V}}$
is equal to $\sigma$. In particular
the set 
${^n \mathcal{S}} 
(d, \mathbf{d}, v, \mathbf{v})$ 
of irreducible components
of ${^n \Pi}^{s, \ast s}_{
D, \mathbf{D}, V, \mathbf{V}}$
is in a natural bijection with the 
set of irreducible components of $X$.

If $\sum_{k=1}^n \mathbf{v}^k =v$ the
map $\sigma_n$ is a vector bundle
\begin{equation}\nonumber
\sigma_n : \quad
{^n \Pi}^{s, \ast s}_{
D, \mathbf{D}, V, \mathbf{V}} 
\rightarrow
\Lambda^{s, \ast s}_{
\mathbf{D}^1/ \mathbf{D}^0, 
\mathbf{V}^1/\mathbf{V}^0}
\times \ldots \times
\Lambda^{s, \ast s}_{
\mathbf{D}^n/\mathbf{D}^{n-1}, 
\mathbf{V}^n/\mathbf{V}^{n-1}} \; ,
\end{equation}
which implies the following
proposition
\begin{proposition}\label{SNonEmpty}
Assume that $\sum_{k=1}^n \mathbf{v}^k =v$,
and 
$\Lambda^{s, \ast s}_{
\mathbf{D}^k/\mathbf{D}^{k-1}, 
\mathbf{V}^k/\mathbf{V}^{k-1}}$
is non-empty for all $k = 1 \ldots n$.
Then the set
${^n \mathcal{S}} (
d, \mathbf{d}, v, \mathbf{v})$
is a one-element set. 
\end{proposition}

\subsection{The second bijection}
\label{SecondBijection}

Recall that the tensor product variety
${^n \Pi}^s_{D, \mathbf{D}, V, \mathbf{v}}$
is a locally closed subset of 
$\Lambda^s_{D,V}$, and
that $\Lambda^s_{D,V,U}$ denotes
the locally closed subset of
$\Lambda^s_{D,V}$ consisting of
all $(x,p,q) \in \Lambda^s_{D,V}$
such that $\underline{q^{-1} (0)}=U$
(cf. \ref{DefinitionOfMCal}).
Let 
${^n \Pi}^s_{D, \mathbf{D}, V, U, 
\mathbf{v}} = 
{^n \Pi}^s_{D, \mathbf{D}, V, \mathbf{v}}
\cap \Lambda^s_{D,V,U}$.
Let $T$ be a $\mathbb{Z} [I]$-graded
subspace of $V$ complimentary to 
$U$.
The vector bundle $\gamma :
\Lambda_{D,V,U} \rightarrow
\Lambda_U \times 
\Lambda^{\ast s}_{D,T}$ 
(cf. Proposition \ref{DefinitionOfMCal})
restricted to
${^n \Pi}^s_{D, \mathbf{D}, V, U, 
\mathbf{v}}$ has the image equal to 
${^\delta \Lambda}_U \times 
{^n \Pi}^{s, \ast s}_{D, \mathbf{D}, T, 
\mathbf{v}}$, where 
${^\delta \Lambda}_U$ is as
in Proposition \ref{LambdaDVUOpen},
and moreover
\begin{equation}\nonumber
\gamma^{-1} ( {^\delta \Lambda}_U \times 
{^n \Pi}^{s, \ast s}_{D, \mathbf{D}, T, 
\mathbf{v}}) \cap 
\Lambda^s_{D,V,U} = 
{^n \Pi}^s_{D, \mathbf{D}, V, U, 
\mathbf{v}} \; .
\end{equation}
Hence
${^n \Pi}^s_{D, \mathbf{D}, V, U, 
\mathbf{v}}$ is open in 
$\gamma^{-1} ({^\delta \Lambda}_U \times 
{^n \Pi}^{s, \ast s}_{D, \mathbf{D}, T, 
\mathbf{v}})$. Since fibers of $\gamma$
have dimension
$<d, u> - <u,v - u> + <Xu , v-u>$
(cf. Proposition \ref{DefinitionOfMCal}),
it follows that
${^n \Pi}^s_{D, \mathbf{D}, V, U, 
\mathbf{v}}$ has pure dimension
\begin{multline}\nonumber
\dim {^n \Pi}^s_{D, \mathbf{D}, V, U, 
\mathbf{v}} 
= \\ = 
\dim {^\delta \Lambda}_U +
\dim 
{^n \Pi}^{s, \ast s}_{D, \mathbf{D}, T, 
\mathbf{v}} +
<d, u> - <u,v - u> + <Xu , v-u>
= \\ =
\dim {^n \Pi}^s_{D, \mathbf{D}, V, 
\mathbf{v}} - <u, v-u> \; ,
\end{multline}
and the set of irreducible
components of 
${^n \Pi}^s_{D, \mathbf{D}, V, U, 
\mathbf{v}}$ is in a natural bijection
(induced by the restriction of the
vector bundle $\gamma$ to
${^n \Pi}^s_{D, \mathbf{D}, V, U, 
\mathbf{v}}$) with the set
\begin{equation}\nonumber
\mathcal{M} (d, v-u, v) \times
{^n \mathcal{S}} (d, \mathbf{d}, v-u, 
\mathbf{v}) \; ,
\end{equation}
where 
$\mathcal{M} (d, v-u, v)$ represents
the set of irreducible components of 
${^\delta \Lambda}_U$, which is
in a natural bijection (induced by
$\gamma$ -- cf. \ref{DefinitionOfMCal}) 
with the set of irreducible
components of the quiver variety
$\mathfrak{M} (d, v-u, v)$, and
${^n \mathcal{S}} (d, \mathbf{d}, v-u, 
\mathbf{v})$ represents 
the set of irreducible
components of the variety
${^n \Pi}^{s, \ast s}_{D, \mathbf{D}, T,
\mathbf{v}}$ (or the multiplicity
variety 
${^n \mathfrak{S}} (d, \mathbf{d}, v-u, 
\mathbf{v})$ -- cf. 
\ref{IrreducibleComponentsOfTS}). 

Let 
${^n \Pi}^s_{D, \mathbf{D}, V, u, 
\mathbf{v}}$ be the set of 
all $(x,p,q) \in 
{^n \Pi}^s_{D, \mathbf{D}, V, 
\mathbf{v}}$ such that 
$\dim \underline{q^{-1} (0)}=u$.
It is the total space of a fibration
over the graded Grassmannian
$Gr_u^V$ with the fiber over
$U \in Gr_u^V$ equal to
${^n \Pi}^s_{D, \mathbf{D}, V, U,
\mathbf{v}}$. Hence
${^n \Pi}^s_{D, \mathbf{D}, V, u, 
\mathbf{v}}$ has pure dimension
\begin{equation}\nonumber
\dim {^n \Pi}^s_{D, \mathbf{D}, V, u, 
\mathbf{v}} =
\dim {^n \Pi}^s_{D, \mathbf{D}, V, 
\mathbf{v}} \; .
\end{equation}
In particular the dimension 
does not depend on $u$. Therefore
one obtains the following
natural bijection 
of sets of
irreducible components
\begin{equation}\label{BetaN}
\beta_n : \;
{^n \mathcal{T}} (d, \mathbf{d},
v, \mathbf{v}) 
\xrightarrow{\sim}
\bigsqcup_{v_0 \in \mathbb{Z}_{\geq 0}}
{^n \mathcal{S}} (d, \mathbf{d}, v_0, 
\mathbf{v})
\times 
\mathcal{M} (d, v_0, v) \; .
\end{equation}

The bijection $\beta_n$ is an analogue
of the direct sum decomposition for
the tensor product of $n$
representations of $\mathfrak{g}$.
In particular \eqref{BetaN}
explains the
name ''multiplicity variety'' given to
${^n \mathfrak{S}} (d, \mathbf{d}, v_0, 
\mathbf{v})$.

\subsection{The tensor decomposition 
bijection}
\label{TensorDecomposition}

Let $\tau_n$ be a bijection
\begin{multline}\label{TauNWeight}
\tau_n : \quad
\bigsqcup_{\substack{
\mathbf{u}^1 , \ldots , \mathbf{u}^n 
\in \mathbb{Z} [I] \\
\mathbf{u}^1 + \ldots + 
\mathbf{u}^n = v }}
\mathcal{M} (\mathbf{d}^1, \mathbf{v}^1,
\mathbf{u}^1)
\times \ldots \times
\mathcal{M} (\mathbf{d}^n, \mathbf{v}^n,
\mathbf{u}^n)
\xrightarrow{\sim} 
\\
\xrightarrow{\sim}
\bigsqcup_{v_0 \in \mathbb{Z}_{\geq 0}}
{^n \mathcal{S}} (d, \mathbf{d},
v_0, \mathbf{v})
\times
\mathcal{M} (d, v_0, v)
\end{multline}
given by 
\begin{equation}\nonumber
\tau_n = \beta_n 
\circ \alpha_n^{-1} \; .
\end{equation}
Taking union over $v$ one obtains
a bijection
\begin{equation}\label{TauN}
\tau_n :
\mathcal{M} (\mathbf{d}^1, \mathbf{v}^1)
\times \ldots \times
\mathcal{M} (\mathbf{d}^n, \mathbf{v}^n)
\xrightarrow{\sim}
\bigsqcup_{v_0 \in \mathbb{Z}_{\geq 0}}
{^n \mathcal{S}} (d, \mathbf{d},
v_0, \mathbf{v})
\times
\mathcal{M} (d, v_0) \; .
\end{equation}

It follows from the definitions of the 
bijections $\alpha_n$ (cf.
\ref{DefinitionOfAlpha}), 
$\beta_n$ 
(cf. \ref{SecondBijection}),
and
$\eta_{k, n-k}$
(cf. \ref{InductionOfMultiplicities}),
that 
\begin{equation}
\tau_n =
( \eta_{k, n-k}^{-1} \times \tau_2) 
\circ 
(\tau_k \times \tau_{n-k}) \; .
\end{equation}
More precisely, the following
diagram of bijections is commutative:
\begin{equation}\label{TauEta}
\xymatrix{
*\txt{ } \\
& 
*\txt{$
\mathcal{M} (\mathbf{d}^1, \mathbf{v}^1)
\times \ldots \times
\mathcal{M} (\mathbf{d}^k, \mathbf{v}^k)
\times
$ \\ $
\times
\mathcal{M} (\mathbf{d}^{k+1}, 
\mathbf{v}^{k+1})
\times \ldots \times
\mathcal{M} (\mathbf{d}^n, \mathbf{v}^n)
$}
\ar[dd]_-{\tau_k \times \tau_{n-k}}
&
\\ \\
*\txt{$\bigsqcup$ \\ 
$\substack{v_0',v_0''}$}
&
*\txt{$
{^k \mathcal{S}}
(d' , \mathbf{d}', v_0', \mathbf{v}')
\times
\mathcal{M} (d' , v_0')
\times
$ \\ $
{^{n-k} \mathcal{S}}
(d'' , \mathbf{d}'', v_0'', \mathbf{v}'')
\times
\mathcal{M} (d'' , v_0'')$}
\ar[dd]_-{\Id \times P_{23}
\times \Id}
\\ \\
*\txt{$\bigsqcup$ \\
$\substack{v_0', v_0''}$}
&
*\txt{${^k \mathcal{S}}
(d' , \mathbf{d}', v_0', \mathbf{v}')
\times
{^{n-k} \mathcal{S}}
(d'' , \mathbf{d}'', v_0'', \mathbf{v}'')
\times
$ \\ $
\mathcal{M} (d' , v_0')
\times
\mathcal{M} (d'' , v_0'')$}
\ar[dd]_-{\Id \times \Id 
\times \tau_2}
\\ \\
*\txt{$\bigsqcup$ \\
$\substack{v_0, v_0', v_0''}$}
&
*\txt{${^k \mathcal{S}}
(d' , \mathbf{d}', v_0', \mathbf{v}')
\times
{^{n-k} \mathcal{S}}
(d'' , \mathbf{d}'', v_0'', \mathbf{v}'')
\times
$ \\ $
{^2 \mathcal{S}}
(d'+d'' , (d', d''), 
v_0, (v_0', v_0''))
\times
\mathcal{M}
(d , v_0)$}
\ar[dd]_-{\eta^{-1}_{k, n-k}
\times \Id}
\\ \\
*\txt{$\bigsqcup$ \\
$\substack{v_0}$}
&
*\txt{${^n \mathcal{S}}
(d'+d'' , (\mathbf{d}', \mathbf{d}''), 
v_0, (\mathbf{v}', \mathbf{v}''))
\times
\mathcal{M}
(d , v_0)$}
\ar@{<-} '[r] '[ruuuuuuuu] '[uuuuuuuu]_-{
\tau_n}
& \\ *\txt{}
}
\end{equation}
where 
$\mathbf{v}' = (\mathbf{v}^1, \ldots ,
\mathbf{v}^k)$,
$\mathbf{v}'' = (\mathbf{v}^{k+1}, \ldots ,
\mathbf{v}^n)$,
$\mathbf{d}' = (\mathbf{d}^1, \ldots ,
\mathbf{d}^k)$,
and
$\mathbf{d}'' = (\mathbf{d}^{k+1}, \ldots ,
\mathbf{d}^n)$.

In the next section it will be shown 
that the set $\mathcal{M} (d, v_0)$ can
be equipped with a structure of a
$\mathfrak{g}$-crystal, and that the 
bijection $\tau_n$ is a crystal
isomorphism (if one replaces direct
products of sets with tensor
products of crystals and considers the
set ${^n \mathcal{S}}$ as a trivial
crystal).

\section{Levi restriction and the crystal 
structure on quiver varieties}
\label{QuiverCrystals}

\subsection{A subquiver $Q'$}
\label{QPrime}
Let $Q'= (I',H') \subset Q = (I,H)$
be a (full) subquiver of $Q$ (i.e. 
the set $I'$ of vertices of $Q'$ is
a subset of the set $I$ of vertices
of $Q$, and two vertices of $Q'$ are
connected by an oriented edge 
$h \in H'$ of $Q'$ if and only 
if they are connected by an 
oriented edge of $Q$).

Let $H^{\sss QQ'}$ 
(resp. $H^{\sss Q'Q}$)
denote the set of edges $h \in H$ 
of $Q$ such that 
$\Out (h) \in Q'$ and 
$\In (h) \notin Q'$
(resp. $\In (h) \in Q'$ and 
$\Out (h) \notin Q'$).

In this section (until 
\ref{MainTheoremForQ}) 
subscripts $Q$ and $Q'$ are used to 
distinguish the varieties
defined using the two quivers.
For example,
${_{\sss Q} \mathfrak{M}}^s_{D,V}$ is a
quiver variety associated to $Q$
(in particular, $D$ and $V$ are
$\mathbb{Z} [I]$-graded 
$\mathbb{C}$-linear spaces),
whereas 
${_{\sss Q'} \mathfrak{M}}^s_{D,V}$ is a
quiver variety associated to $Q'$
(and, correspondingly, $D$ and $V$ are
$\mathbb{Z} [I']$-graded 
$\mathbb{C}$-linear spaces, or
$\mathbb{Z} [I]$-graded spaces with
zero components in degrees 
$i \in I \backslash I'$),
${_{\sss Q} \mathcal{F}}$ is the 
path algebra of $Q$, whereas 
${_{\sss Q'} \mathcal{F}}$ is 
the path algebra of $Q'$. Note that 
${_{\sss Q'} \mathcal{F}}$
is a subalgebra of 
${_{\sss Q} \mathcal{F}}$.

\subsection{The set 
${_{\sss QQ'} 
\mathcal{M}} (d, v_0, v)$}
\label{DefinitionOfZeta}
Let $D$ and $V$ be 
$\mathbb{Z} [I]$-graded 
$\mathbb{C}$-linear spaces
with graded dimensions $d$ and $v$
respectively, and let 
$D'$, $V'$ be $\mathbb{Z} [I']$-graded 
$\mathbb{C}$-linear spaces defined as
follows
\begin{equation}\nonumber
\begin{split}
D'_i &= D_i \oplus 
(\oplus^{h \in H^{QQ'}}_{\Out (h) = i}
V_{\In (h)}) \; ,
\\
V'_i &= V_i \; ,
\end{split}
\end{equation}
for $i \in I' \subset I$.
One has a regular map
$\zeta_{\sss QQ'} : \;
{_{\sss Q} \Lambda}_{D,V}
\rightarrow
{_{\sss Q'} \Lambda}_{D',V'}$
given by
\begin{equation}\nonumber
\zeta_{\sss QQ'} ((x,p,q)) = 
(x', p' , q' ) \; ,
\end{equation}
where
\begin{equation}\nonumber
\begin{split}
x'_h &= x_h \text{ for }
h \in H' \subset H \; ,
\\
p'_i &= \bigl( p_i , 
\{ x_h \}^{
h \in H^{\sss Q'Q}}_{\In (h) = i}
\bigr) \; ,
\\
q'_i &= \bigl( q_i , 
\{ (\varepsilon (h))^{-1} x_h \}^{
h \in H^{\sss QQ'}}_{\Out (h) = i}
\bigr) \; .
\end{split}
\end{equation}
Let ${_{\sss QQ'} 
\Lambda}^{\ast s}_{D,V}$
be the set of all 
$(x, p, q) 
\in {_{\sss Q} \Lambda}_{D,V}$
such that
$\zeta_{\sss QQ'} ((x, p, q)) \in
{_{\sss Q'} 
\Lambda}^{\ast s}_{D',V'}$.
Let ${_{\sss QQ'} 
\Lambda}^{s, \ast s}_{D,V} =
{_{\sss QQ'} 
\Lambda}^{\ast s}_{D,V} \cap
{_{\sss Q} 
\Lambda}^s_{D,V}$.
Roughly speaking, 
${_{\sss QQ'} 
\Lambda}^{s, \ast s}_{D,V}$
is the open subset of 
$\Lambda_{D,V}$ consisting
of all stable points that are also
$\ast$-stable ``at vertices of $Q'$''.

Let $v_0 \in \mathbb{Z} [I]$ and
let ${_{\sss QQ'} 
\Lambda}^{s, \ast s}_{D,V,v_0}=
{_{\sss QQ'} 
\Lambda}^{s, \ast s}_{D,V}
\cap
{_{\sss Q} \Lambda}^s_{D,V,v_0}$.
Since 
${_{\sss QQ'} 
\Lambda}^{s, \ast s}_{D,V,v_0}$ is
an open $G_V$-invariant subset of
${_{\sss Q} \Lambda}^s_{D,V,v_0}$ one
has the following proposition
(cf. \ref{LambdaDVUDimSmall},
\ref{StableFreeG}).
\begin{proposition}
${_{\sss QQ'} 
\Lambda}^{s, \ast s}_{D,V,v_0}$ is empty
or has pure dimension 
\begin{multline}\nonumber
\dim {_{\sss QQ'} 
\Lambda}^{s, \ast s}_{D,V,v_0} 
= \\ =
\frac{1}{2}<X v, v> + 
\frac{1}{2}<X v_0, v_0> +
<d, v> + <d, v_0> - <v_0, v_0> \; ,  
\end{multline}
and $G_V$-action on
${_{\sss QQ'} 
\Lambda}^{s, \ast s}_{D,V,v_0}$
is free.
\end{proposition}
It follows that
${_{\sss QQ'} 
\mathfrak{M}}^{s, \ast s} (d, v_0 , v) =
{_{\sss QQ'} 
\Lambda}^{s, \ast s}_{D,V,v_0} / G_V$
is naturally a quasi-projective variety
of pure dimension
\begin{multline}\nonumber
\dim 
{_{\sss QQ'} 
\mathfrak{M}}^{s, \ast s} (d ,v_0 ,v) = 
\frac{1}{2}<X v, v> + 
\frac{1}{2}<X v_0, v_0> +
\\ + <d, v> + <d, v_0> - <v_0, v_0> - 
<v, v> \; .
\end{multline}
Let 
${_{\sss QQ'} \mathcal{M}} (d, v_0, v)$
be the set of irreducible components of
${_{\sss QQ'} 
\mathfrak{M}}^{s, \ast s} (d, v_0, v)$.
This set is also in a natural bijection 
(cf. the end of Section
\ref{DefinitionOfMCal}) with
the sets of irreducible components of
${_{\sss QQ'} 
\Lambda}^{s, \ast s}_{D,V,v_0}$, 
and of
${_{\sss QQ'} 
\Lambda}^{s, \ast s}_{D,V,U} = 
{_{\sss QQ'} 
\Lambda}^{s, \ast s}_{D,V}
\cap
{_{\sss Q} \Lambda}^s_{D,V,U}$
(where $\dim U = v-v_0$).

\subsection{Levi restriction}
\label{DefinitionOfTheta}

Given $(x,p,q) \in 
{_{\sss Q} \Lambda}^s_{D,V}$
let ${_{\sss QQ'} \mathcal{K}} (x,q)$ 
denote the maximal graded subspace 
$U \subset V$, satisfying the following
conditions:
\begin{align}\nonumber
U_i &= \{ 0 \}
&&\text{ for any } i \notin I' \; , 
\\\label{RestrictionSmallKernel}
x_h U &\subset U 
&&\text{ for any } h \in H' \; ,
\\\nonumber
x_h U &= 0 
&&\text{ for any } h \in H^{\sss QQ'} \; ,
\\\nonumber
q_i U &= 0 
&&\text{ for any } i \in I' \; .
\end{align}
Let $w \in \mathbb{Z} [I'] \subset
\mathbb{Z} [I]$, and
$W$ be a graded subspace of $V$
with $\dim W \in \mathbb{Z} [I'] \subset
\mathbb{Z} [I]$. Then 
${_{\sss QQ'} 
\Lambda}^s_{D,V,v_0,w}$ 
(resp. 
${_{\sss QQ'} 
\Lambda}^s_{D,V,v_0,W}$)
denotes the set of all
$(x,p,q) \in
{_{\sss Q} 
\Lambda}^s_{D,V,v_0}$
such that
$\dim {_{\sss QQ'} \mathcal{K}} 
(x,q) = w$
(resp.
${_{\sss QQ'} \mathcal{K}} 
(x,q) = W$).
Note that 
${_{\sss QQ'} 
\Lambda}^s_{D,V,v_0,w}$ 
is a fibration over a graded
Grassmannian $Gr_w^V$ with 
fibers isomorphic to
${_{\sss QQ'} 
\Lambda}^s_{D,V,v_0,W}$,
where $\dim W = w$.

Let $W$ be as above and $T$ be
a graded subspace of $V$ complimentary
to $W$. Then one has a regular map 
(cf. \ref{DefinitionOfMCal})
\begin{equation}\nonumber
\nu_{\sss QQ'} : \;
{_{\sss QQ'} 
\Lambda}^s_{D,V,v_0,W}
\rightarrow
{_{\sss Q'} \Lambda}_W \times
{_{\sss QQ'} 
\Lambda}^{s, \ast s}_{D,T,v_0}
\end{equation}
given by 
\begin{equation}\nonumber
\nu_{\sss QQ'} 
((x,p,q)) = (x^{\sss WW},
(x^{\sss TT}, p^{\sss TD} ,
q^{\sss DT})) \; .
\end{equation}

The fiber of $\nu_{\sss QQ'}$ 
over a point $(x, (y,p,q)) \in 
{_{\sss Q'} \Lambda}_W \times
{_{\sss QQ'} 
\Lambda}^{s, \ast s}_{D,T,v_0}$
is the same
as the fiber of the
map $\gamma$ 
(cf. \ref{DefinitionOfMCal})
for the quiver $Q'$ over
the point 
$(x, \zeta_{\sss QQ'} ((y, p, q))) \in 
{_{\sss Q'} \Lambda}_W \times
{_{\sss Q'} 
\Lambda}^{s, \ast s}_{D',T'}$,
where $T'_i = T_i$ for 
$i \in I' \subset I$, and
$D'$ and $\zeta_{\sss QQ'}$ 
are defined in 
\ref{DefinitionOfZeta}. 
Therefore an analogue
of the Proposition 
\ref{DefinitionOfMCal}
holds and one obtains
the following bijection
between sets of
irreducible components:
\begin{multline}\nonumber
\theta_{\sss QQ'} : \quad
{_{\sss Q} \mathcal{M}}
(d, v_0, v) 
\xrightarrow{\sim}
\\
\xrightarrow{\sim}
\bigsqcup_{u \in 
\mathbb{Z}_{\geq 0} [I']}
{_{\sss QQ'} \mathcal{M}}
(d, v_0, v-u)
\times
{_{\sss Q'} \mathcal{M}}
(\delta_{\sss QQ'} (d,v) , 
\rho_{\sss QQ'} (v-u), 
\rho_{\sss QQ'} (v))
\end{multline}
where
\begin{align}\nonumber 
(\rho_{\sss QQ'} (v))_i &= v_i 
&&\text{ for } \quad
v \in \mathbb{Z} [I] \; , \;
i \in I' \subset I \; , 
\\\nonumber
(\delta_{\sss QQ'} (d, v))_i &=
d_i + \sum_{\substack{h \in H^{QQ'} \\
\Out (h) = i}}
v_{\In (h)}
&&\text{ for } \quad 
d, v \in \mathbb{Z} [I] \; , \;
i \in I' \subset I \; .
\end{align} 
Union over $v$ gives a bijection
\begin{equation}\nonumber
\theta_{\sss QQ'} : \;
{_{\sss Q} \mathcal{M}}
(d, v_0) 
\xrightarrow{\sim}
\bigsqcup_{x \in 
\mathbb{Z}_{\geq 0} [I]}
{_{\sss QQ'} \mathcal{M}}
(d, v_0, x)
\times
{_{\sss Q'} \mathcal{M}}
(\delta_{\sss QQ'} (d,x) , 
\rho_{\sss QQ'} (x))
\end{equation}

The bijection $\theta_{\sss QQ'}$
is an analogue of the restriction
of a representation of $\mathfrak{g}$
to a Levi subalgebra of a parabolic
subalgebra of $\mathfrak{g}$. This
construction is a straightforward
generalization of the reduction
to $\mathfrak{sl}_2$ subalgebras
introduced by Lusztig 
\cite[12]{Lusztig1991a}
and used by Nakajima
\cite{Nakajima1998},
Kashiwara and Saito
\cite{KashiwaraSaito}.

\subsection{Levi restriction and the second
bijection for tensor product varieties}
(cf. \ref{SecondBijection})
\label{LeviTensorSecond}

Let $D$, $V$ be $\mathbb{Z} [I]$-graded
$\mathbb{C}$-linear spaces,
$\mathbf{D}$ be an $n$-step partial
flag in $D$, and 
$\mathbf{v}$ be an $n$-tuple
of elements of $\mathbb{Z} [I]$.
Recall 
(cf. \ref{DefinitionOfTS}) that
${^n_{\sss Q} \Pi}^s_{
D, \mathbf{D}, V, \mathbf{v}}$
is a locally closed subset of
${_{\sss Q} \Lambda}^s_{D, V}$.
Let 
\begin{equation}\nonumber
{_{\sss QQ'}^n \Pi}^{s, \ast s}_{
D, \mathbf{D}, V, \mathbf{v}} =
{_{\sss Q}^n \Pi}^s_{
D, \mathbf{D}, V, \mathbf{v}}
\cap
{_{\sss QQ'} \Lambda}^{
s, \ast s}_{D, V} \; .
\end{equation}
Then ${_{\sss QQ'}^n \Pi}^{s, \ast s}_{
D, \mathbf{D}, V, \mathbf{v}}$ is an 
open subset of
${_{\sss Q}^n \Pi}^s_{
D, \mathbf{D}, V, \mathbf{v}}$ and
therefore it is empty or 
has pure dimension
equal to that of
${_{\sss Q}^n \Pi}^s_{
D, \mathbf{D}, V, \mathbf{v}}$.
Let
${_{\sss QQ'}^n \mathcal{T}}
(d, \mathbf{d}, v, \mathbf{v})$
be the set of irreducible components of 
${_{\sss QQ'}^n \Pi}^{s, \ast s}_{
D, \mathbf{D}, V, \mathbf{v}}$.

The restriction of the map
$\gamma: \; 
{_{\sss Q} \Lambda}_{D,V,U}
\rightarrow 
{_{\sss Q}\Lambda}_U
\times 
{_{\sss Q} \Lambda}_{D,T}^{\ast s}$
(cf. \ref{DefinitionOfMCal})
to ${_{\sss QQ'} 
\Lambda}^{s, \ast s}_{D,V,U}=
{_{\sss QQ'} \Lambda}^{s, 
\ast s}_{D,V} \cap 
{_{\sss Q} \Lambda}_{D,V,U}$
has the image equal to
${_{\sss Q}^{\delta} \Lambda}'_U
\times
{_{\sss Q} \Lambda}^{s, \ast s}_{D,T}$,
where 
$\delta$ is as in 
\ref{LambdaDVUOpen}, and
${_{\sss Q}^{\delta} \Lambda}'_U$
denotes the open subset of
${_{\sss Q}^{\delta} \Lambda}_U$
consisting of all $x \in
{_{\sss Q}^{\delta} \Lambda}_U$
such that 
$\cap_{\substack{h \in H \\ 
\Out (h) =i }} \ker x_h = \{ 0 \}$
for any $i \in I'$. 
It follows that the restriction
of the map $\gamma$ to 
${_{\sss QQ'}^n \Pi}^{s, \ast s}_{
D, \mathbf{D}, V, \mathbf{v}}
\cap {_{\sss Q} \Lambda}_{D,V,U}$
has the image equal to
${_{\sss Q}^{\delta} \Lambda}'_U
\times
{_{\sss Q}^n \Pi}^{s, \ast s}_{
D, \mathbf{D}, T, \mathbf{v}}$, and
one can repeat the argument in
\ref{SecondBijection} 
to get a bijection
\begin{equation}\nonumber
{_{\sss QQ'} \beta}_n : \;
{_{\sss QQ'}^n \mathcal{T}} (d, \mathbf{d},
v, \mathbf{v}) 
\xrightarrow{\sim}
\bigsqcup_{v_0 \in \mathbb{Z}_{\geq 0} [I]}
{_{\sss Q}^n \mathcal{S}} 
(d, \mathbf{d}, v_0, \mathbf{v})
\times 
{_{\sss QQ'} \mathcal{M}} (d, v_0, v) \; .
\end{equation}

\subsection{Levi restriction and the first
bijection for tensor product varieties}
(cf. \ref{DefinitionOfAlpha})
\label{LeviTensorFirst}

Let $\mathbf{v}$, $\Tilde{\mathbf{v}}$ 
be $n$-tuples of elements of 
$\mathbb{Z} [I]$,
$\Hat{\mathbf{v}}$ be an 
$n$-tuple of elements of 
$\mathbb{Z} [I'] \subset 
\mathbb{Z} [I]$, such that
$\sum_{k=1}^n \mathbf{v}^k +
\sum_{k=1}^n \Tilde{\mathbf{v}}^k +
\sum_{k=1}^n \Hat{\mathbf{v}}^k =
\dim V$,
and let
${_{\sss QQ'} ^n \Pi}^{s, \ast s}_{D, 
\mathbf{D}, V, \mathbf{v}, 
\Tilde{\mathbf{v}}, \Hat{\mathbf{v}}}$ be
a locally closed subset of 
${_{\sss QQ'} ^n \Pi}^{s, \ast s}_{D, 
\mathbf{D}, V, \mathbf{v}} \cap
{_{\sss Q} ^n \Pi}^s_{D, 
\mathbf{D}, V, \mathbf{v}, 
\Tilde{\mathbf{v}}+\Hat{\mathbf{v}}}$, 
consisting of all $(x, p, q) \in 
{_{\sss QQ'} ^n \Pi}^{s, \ast s}_{D, 
\mathbf{D}, V, \mathbf{v}} \cap
{_{\sss Q} ^n \Pi}^s_{D, 
\mathbf{D}, V, \mathbf{v}, 
\Tilde{\mathbf{v}}+\Hat{\mathbf{v}}}$ 
such that
the dimension of the maximal graded
subspace $\Hat{\mathbf{V}}^k$ of 
$\overline{p (\mathbf{D}^k)}$
satisfying the following conditions
(cf. \ref{RestrictionSmallKernel})
\begin{align}\nonumber
\Hat{\mathbf{V}}^k_i &\subset
\bigl( \overline{p (\mathbf{D}^{k-1})} 
\bigr)_i
&&\text{ for any } i \notin I' \; , 
\\\label{RestrictionThree}
x_h \Hat{\mathbf{V}}^k &\subset 
\Hat{\mathbf{V}}^k 
&&\text{ for any } h \in H' \; ,
\\\nonumber
x_h \Hat{\mathbf{V}}^k &\subset 
\overline{p (\mathbf{D}^{k-1})} 
&&\text{ for any } h \in H^{\sss QQ'} \; ,
\\\nonumber
q_i ( \Hat{\mathbf{V}}^k_i ) &\subset 
\mathbf{D}^{k-1}_i 
&&\text{ for any } i \in I' \; ,
\end{align}
is equal
$\sum_{l=1}^{k-1} \mathbf{v}^l +
\sum_{l=1}^{k-1} \Tilde{\mathbf{v}}^l +
\sum_{l=1}^k \Hat{\mathbf{v}}^l$
(cf. the definition of
${_{\sss Q} ^n \Pi}^s_{D, 
\mathbf{D}, V, \mathbf{v}, 
\Tilde{\mathbf{v}}}$
in \ref{DefinitionOfTS}).
Note that the maximal subspace
$\Hat{\mathbf{V}}^k$
satisfying conditions
\ref{RestrictionThree}
contains 
$\overline{p (\mathbf{D}^{k-1})}$,
which has dimension
$\sum_{l=1}^{k-1} \mathbf{v}^l +
\sum_{l=1}^{k-1} \Tilde{\mathbf{v}}^l +
\sum_{l=1}^{k-1} \Hat{\mathbf{v}}^l$,
and is contained in
$\underline{q^{-1} (\mathbf{D}^{k-1})}$,
which has dimension
$\sum_{l=1}^{k-1} \mathbf{v}^l +
\sum_{l=1}^{k} \Tilde{\mathbf{v}}^l +
\sum_{l=1}^{k} \Hat{\mathbf{v}}^l$.

Let $\mathbf{V} =
( 0 = \mathbf{V}^0 =
\Hat{\mathbf{V}}^1
\subset
\Tilde{\mathbf{V}}^1
\subset
\mathbf{V}^1
\subset
\Hat{\mathbf{V}}^2
\subset
\Tilde{\mathbf{V}}^2
\subset
\ldots
\subset
\Hat{\mathbf{V}}^n
\subset
\Tilde{\mathbf{V}}^n
\subset
\mathbf{V}^n = V )$ be
a $3n$-step $\mathbb{Z} [I]$-graded
partial flag in $V$, 
such that
$\dim \Hat{\mathbf{V}}^k -
\dim \mathbf{V}^{k-1} =
\Hat{\mathbf{v}}^k$,
$\dim \Tilde{\mathbf{V}}^k -
\dim \Hat{\mathbf{V}}^k =
\Tilde{\mathbf{v}}^k$,
$\dim \mathbf{V}^k -
\dim \Tilde{\mathbf{V}}^k =
\mathbf{v}^k$,
and let
${_{\sss QQ'}^n \Pi}^{s, \ast s}_{
D, \mathbf{D}, V, \mathbf{V}}$
be a locally closed subset of 
${_{\sss QQ'}^n \Pi}^{s, \ast s}_{
D, \mathbf{D}, V, \mathbf{v}}$ 
consisting of all
$(x,p,q) \in
{_{\sss QQ'}^n \Pi}^{s, \ast s}_{
D, \mathbf{D}, V, \mathbf{v}}$ 
such that
$\mathbf{V}^k = 
\overline{p (\mathbf{D}^k)}$,
$\Tilde{\mathbf{V}}^k =
\underline{q^{-1} (\mathbf{D}^{k-1})}
\cap \mathbf{V}^k$, and
$\Hat{\mathbf{V}}^k$ is the maximal
graded subspace of 
$\mathbf{V}^k$ satisfying conditions
\ref{RestrictionThree}
(cf. the definition of
${_{\sss Q} ^n \Pi}^s_{D, 
\mathbf{D}, V, \mathbf{V}}$
in \ref{DefinitionOfTS}).

The variety 
${_{\sss QQ'}^n \Pi}^{s, \ast s}_{
D, \mathbf{D}, V, \mathbf{v},
\Tilde{\mathbf{v}}, \Hat{\mathbf{v}}}$
is a fibration over the variety
of $3n$-step partial flags in $V$ with
dimensions of the subfactors
given by $\mathbf{v}$,
$\Tilde{\mathbf{v}}$, and
$\Hat{\mathbf{v}}$, and
a fiber of this fibration 
is isomorphic to
${_{\sss QQ'}^n \Pi}^{s, \ast s}_{
D, \mathbf{D}, V, \mathbf{V}}$.

One has a regular map (given by
restrictions of $x$, $p$, and $q$)
\begin{multline}\nonumber
{_{\sss QQ'} \sigma}_n : \quad 
{_{\sss QQ'}^n \Pi}^{s, \ast s}_{
D, \mathbf{D}, V, \mathbf{V}} 
\rightarrow \\ \rightarrow
{_{\sss QQ'} \Lambda}^{s, \ast s}_{
\mathbf{D}^1 / \mathbf{D}^0, 
\mathbf{V}^1/\Hat{\mathbf{V}}^1,
\Tilde{\mathbf{V}}^1/\Hat{\mathbf{V}}^1}
\times \ldots \times
{_{\sss QQ'} \Lambda}^{s, \ast s}_{
\mathbf{D}^n / \mathbf{D}^{n-1}, 
\mathbf{V}^n/\Hat{\mathbf{V}}^n,
\Tilde{\mathbf{V}}^n/\Hat{\mathbf{V}}^n}
\times \\ \times
{_{\sss Q'} \Lambda}_{
\Hat{\mathbf{V}}^1/\mathbf{V}^0}
\times \ldots \times
{_{\sss Q'} \Lambda}_{
\Hat{\mathbf{V}}^n/\mathbf{V}^{n-1}} \; .
\end{multline}

Let $D'$, $V'$ be 
$\mathbb{Z} [I']$-graded 
$\mathbb{C}$-linear 
spaces defined in \ref{DefinitionOfZeta},
$\mathbf{D}'=
(0 = \mathbf{D}^{'0} \subset 
\mathbf{D}^{'1}
\subset \ldots \subset 
\mathbf{D}^{'n} =
D')$ (resp. $\mathbf{V}' =
(0 = \mathbf{V}^{'0} =
\Tilde{\mathbf{V}}^{'1} \subset
\mathbf{V}^{'1} \subset 
\Tilde{\mathbf{V}}^{'2}
\subset \ldots \subset
\Tilde{\mathbf{V}}^{'n} \subset
\mathbf{V}^{'n} = V'$) be
an $n$-step partial flag in
$D'$ (resp. a $2n$-step partial flag in
$V'$) induced by the flags
$(0 = \mathbf{D}^0 \subset \mathbf{D}^1
\subset \ldots \subset \mathbf{D}^n =
D)$ and $(0 = \mathbf{V}^0 \subset
\mathbf{V}^1 \subset \ldots \subset
\mathbf{V}^n = V)$ (resp. by the flag
$(0 = \mathbf{V}^0 =
\Hat{\mathbf{V}}^1 \subset
\mathbf{V}^1 \subset \Hat{\mathbf{V}}^2
\subset \ldots \subset
\Hat{\mathbf{V}}^n \subset
\mathbf{V}^n = V$)).

Recall (cf. 
\ref{InductionOfMultiplicities})
that 
\begin{multline}\nonumber
{_{\sss Q'} \sigma}_n : \quad 
{_{\sss Q'}^n \Pi}^{s, \ast s}_{
D', \mathbf{D}', V', \mathbf{V}'} 
\rightarrow \\ \rightarrow
{_{\sss Q'} \Lambda}^{s, \ast s}_{
\mathbf{D}^{'1} / \mathbf{D}^{'0}, 
\mathbf{V}^{'1}/\Tilde{\mathbf{V}}^{'1}}
\times \ldots \times
{_{\sss Q'} \Lambda}^{s, \ast s}_{
\mathbf{D}^{'n} / \mathbf{D}^{'n-1}, 
\mathbf{V}^{'n}/\Tilde{\mathbf{V}}^{'n}}
\times \\ \times
{_{\sss Q'} \Lambda}_{
\Tilde{\mathbf{V}}^{'1}/\mathbf{V}^{'0}}
\times \ldots \times
{_{\sss Q'} \Lambda}_{
\Tilde{\mathbf{V}}^{'n}/\mathbf{V}^{'n-1}}
\end{multline}
is the map obtained by restricting
$x$, $p$, and $q$. Let $\zeta'_{\sss QQ'}$
denote the restriction of the map 
$\zeta_{\sss QQ'}$ (cf. \ref{DefinitionOfZeta})
to ${_{\sss QQ'}^n \Pi}^{s, \ast s}_{
D, \mathbf{D}, V, \mathbf{V}} \subset
{_{\sss Q} \Lambda}_{D,V}$.

The following proposition can be proven
by an inductive (in $n$) argument similar 
to the ones used in subsections 
\ref{InductionOfPi}--
\ref{InductionOfMultiplicities}.
Since the proof is completely
analogous to the proofs in these 
subsections it is omitted.

\begin{proposition}
\begin{alphenum}
\item
The image of the map
${_{\sss QQ'} \sigma}_n$ is equal to
\begin{equation}\nonumber
{_{\sss QQ'} \Lambda}^{s, \ast s}_{
\mathbf{D}^1 / \mathbf{D}^0, 
\mathbf{V}^1/\Hat{\mathbf{V}}^1,
\Tilde{\mathbf{V}}^1/\Hat{\mathbf{V}}^1}
\times \ldots \times
{_{\sss QQ'} \Lambda}^{s, \ast s}_{
\mathbf{D}^n / \mathbf{D}^{n-1}, 
\mathbf{V}^n/\Hat{\mathbf{V}}^n,
\Tilde{\mathbf{V}}^n/\Hat{\mathbf{V}}^n}
\times X \; ,
\end{equation}
where $X$ is the open subset of
${_{\sss Q'} \Lambda}_{
\Hat{\mathbf{V}}^1/\mathbf{V}^0}
\times \ldots \times
{_{\sss Q'} \Lambda}_{
\Hat{\mathbf{V}}^n/\mathbf{V}^{n-1}}$
such that
\begin{equation}\nonumber
{_{\sss Q'} \Lambda}^{s, \ast s}_{
\mathbf{D}^{'1} / \mathbf{D}^{'0}, 
\mathbf{V}^{'1}/\Tilde{\mathbf{V}}^{'1}}
\times \ldots \times
{_{\sss Q'} \Lambda}^{s, \ast s}_{
\mathbf{D}^{'n} / \mathbf{D}^{'n-1}, 
\mathbf{V}^{'n}/\Tilde{\mathbf{V}}^{'n}}
\times X
\end{equation}
is the image of the map
${_{\sss Q'} \sigma}_n \circ 
\zeta'_{\sss QQ'}$. 
\item
The set
${_{\sss QQ'}^n \Pi}^{s, \ast s}_{
D, \mathbf{D}, V, \mathbf{V}}$ is an
open dense subset of the
total space of a locally trivial
fibration over the image of
${_{\sss QQ'} \sigma}_n$ with a
smooth connected fiber and such
that the restriction of the projection
map onto  
${_{\sss QQ'}^n \Pi}^{s, \ast s}_{
D, \mathbf{D}, V, \mathbf{V}}$
is equal to
${_{\sss QQ'} \sigma}_n$.
\item
The dimension of 
${_{\sss QQ'}^n \Pi}^{s, \ast s}_{
D, \mathbf{D}, V, \mathbf{v},
\Tilde{\mathbf{v}}, 
\Hat{\mathbf{v}}}$ does
not depend on
$\Tilde{\mathbf{v}}$ and
$\Hat{\mathbf{v}}$ (it only
depends on $d$, $\mathbf{d}$,
$v$, and $\mathbf{v}$.
\end{alphenum}
\end{proposition}
The above proposition 
together with 
\ref{InductionOfMultiplicities}
implies that
the map ${_{\sss QQ'} \sigma}_n$ 
induces the following bijection
of sets of irreducible components 
\begin{multline}
{_{\sss QQ'} \alpha}_n : \quad
{_{\sss QQ'}^n \mathcal{T}} 
(d, \mathbf{d}, v, \mathbf{v})
\xrightarrow{\sim} 
\\
\xrightarrow{\sim}
\bigsqcup_{\substack{
\mathbf{u}^1 , \ldots , \mathbf{u}^n 
\in \mathbb{Z}_{\geq 0} [I] \\
v - \mathbf{u}^1 - \ldots -
\mathbf{u}^n \in
\mathbb{Z}_{\geq 0} [I']
}}
{_{\sss Q'}^n \mathcal{S}}
(\delta_{\sss QQ'} (d,v),
\delta_{\sss QQ'} (\mathbf{d},\mathbf{u}),
\rho_{\sss QQ'} (v),
\rho_{\sss QQ'} (\mathbf{u}))
\times \\ \times 
{_{\sss QQ'} \mathcal{M}}
(\mathbf{d}^1 , \mathbf{v}^1, 
\mathbf{u}^1) 
\times \ldots \times
{_{\sss QQ'} \mathcal{M}}
(\mathbf{d}^n , \mathbf{v}^n, 
\mathbf{u}^n) \; ,
\end{multline}
where $\delta_{\sss QQ'} (d,v)$, and
$\rho_{\sss QQ'} (v)$ are as in
\ref{DefinitionOfTheta}, and
$\bigl( \delta_{\sss QQ'} 
(\mathbf{d},\mathbf{v})\bigr)^k =
\delta_{\sss QQ'} 
(\mathbf{d}^k, \mathbf{v}^k)$,
$\bigl( \rho_{\sss QQ'} 
(\mathbf{v})\bigr)^k =
\rho_{\sss QQ'} 
(\mathbf{v}^k)$.

\subsection{Levi restriction and the 
tensor product decomposition}

Let ${_{\sss QQ'}\tau}_n$ be the bijection
\begin{multline}\nonumber
\bigsqcup_{\substack{
\mathbf{u}^1 , \ldots , \mathbf{u}^n 
\in \mathbb{Z}_{\geq 0} [I] \\
v - \mathbf{u}^1 - \ldots -
\mathbf{u}^n \in
\mathbb{Z}_{\geq 0} [I']
}}
{_{\sss Q'}^n \mathcal{S}}
(\delta_{\sss QQ'} (d,v),
\delta_{\sss QQ'} (\mathbf{d},\mathbf{u}),
\rho_{\sss QQ'} (v),
\rho_{\sss QQ'} (\mathbf{u}))
\times \\ \times
{_{\sss QQ'} \mathcal{M}}
(\mathbf{d}^1 , \mathbf{v}^1, 
\mathbf{u}^1) 
\times \ldots \times
{_{\sss QQ'} \mathcal{M}}
(\mathbf{d}^n , \mathbf{v}^n, 
\mathbf{u}^n)
\xrightarrow{\sim} \\ \xrightarrow{\sim}
\bigsqcup_{v_0 \in \mathbb{Z}_{\geq 0} [I]}
{_{\sss Q}^n \mathcal{S}} 
(d, \mathbf{d}, v_0, \mathbf{v})
\times 
{_{\sss QQ'} \mathcal{M}} (d, v_0, v) \; .
\end{multline}
given by 
\begin{equation}\nonumber
{_{\sss QQ'}\tau}_n = 
{_{\sss QQ'} \beta}_n 
\circ 
{_{\sss QQ'} \alpha}_n^{-1} \; .
\end{equation}

It follows from the definitions
of the bijections 
${_{\sss QQ'} \alpha}_n$
(cf. \ref{LeviTensorFirst}),
${_{\sss QQ'} \beta}_n$
(cf. \ref{LeviTensorSecond}),
and
$\theta_{\sss QQ'}$
(cf. \ref{DefinitionOfTheta}), 
that the following diagram of
bijections is commutative:

\begin{equation}\label{TauTheta}
\xymatrix{
*\txt{$\bigsqcup$ 
\\ 
$\substack{
\mathbf{u}^1 , \ldots , \mathbf{u}^n 
\in \mathbb{Z}_{\geq 0} [I] \\
\mathbf{u}^1 + \ldots + \mathbf{u}^n = v
}$}
&
*\txt{$
{_{\sss Q} \mathcal{M}} 
(\mathbf{d}^1, \mathbf{v}^1, \mathbf{u}^1)
\times \ldots \times
{_{\sss Q} \mathcal{M}}
(\mathbf{d}^n, \mathbf{v}^n, \mathbf{u}^n)
$}
\ar[dd]^-{\theta_{QQ'} \times \ldots
\times \theta_{QQ'}}
& 
\\ \\
*\txt{$\bigsqcup$ 
\\ 
$\substack{
\mathbf{u}^1 , \ldots , \mathbf{u}^n 
\in \mathbb{Z}_{\geq 0} [I]  \\
\mathbf{u}^1 + \ldots + \mathbf{u}^n = v \\
\mathbf{w}^1 , \ldots , \mathbf{w}^n  
\in \mathbb{Z}_{\geq 0} [I'] }$}
&
*\txt{$
{_{\sss QQ'} \mathcal{M}} 
(\mathbf{d}^1, \mathbf{v}^1, 
\mathbf{u}^1 - \mathbf{w}^1)
\times$ \\ $
{_{\sss Q'} \mathcal{M}} 
(\delta_{\sss QQ'} (\mathbf{d}^1, \mathbf{u}^1), 
\rho_{\sss QQ'} (\mathbf{u}^1 - \mathbf{w}^1), 
\rho_{\sss QQ'} (\mathbf{u}^1)) 
\times $ \\ $\vdots$ \\ $\times
{_{\sss QQ'} \mathcal{M}} 
(\mathbf{d}^n, \mathbf{v}^n, 
\mathbf{u}^n - \mathbf{w}^n)
\times $ \\ $
{_{\sss Q'} \mathcal{M}} 
(\delta_{\sss QQ'} (\mathbf{d}^n, \mathbf{u}^n), 
\rho_{\sss QQ'} (\mathbf{u}^n - \mathbf{w}^n), 
\rho_{\sss QQ'} (\mathbf{u}^n))
$}
\ar[dd]^-{(\Id \times \ldots \times \Id 
\times {_{\sss Q'} \tau}_n) \circ P}
&
\\ \\ 
*\txt{$\bigsqcup$
\\
$\substack{
\mathbf{x}^1 , \ldots , \mathbf{x}^n 
\in \mathbb{Z}_{\geq 0} [I] \\
v - \mathbf{x}^1 - \ldots - \mathbf{x}^n 
\in \mathbb{Z}_{\geq 0} [I'] \\
w \in \mathbb{Z}_{\geq 0} [I'] }$}
&
*\txt{$
{_{\sss QQ'} \mathcal{M}} 
(\mathbf{d}^1, \mathbf{v}^1, \mathbf{x}^1)
\times \ldots \times
{_{\sss QQ'} \mathcal{M}} 
(\mathbf{d}^n, \mathbf{v}^n, \mathbf{x}^n)
\times $ \\ $
{_{\sss Q'}^n \mathcal{S}} 
(\delta_{\sss QQ'} (d,v), 
\delta_{\sss QQ'} (\mathbf{d}, \mathbf{x}), 
\rho_{\sss QQ'} (v - w), 
\rho_{\sss QQ'} (\mathbf{x}))
\times$ \\ $
{_{\sss Q'} \mathcal{M}} 
(\delta_{\sss QQ'} (d,v), 
\rho_{\sss QQ'} (v-w), 
\rho_{\sss QQ'} (v))
$}
\ar[dd]^-{{_{\sss QQ'} \tau}_n \times \Id} 
&
\\ \\ 
*\txt{$\bigsqcup$
\\
$\substack{
v_0 \in \mathbb{Z}_{\geq 0} [I] \\
w \in \mathbb{Z}_{\geq 0} [I'] }$}
&
*\txt{$
{_{\sss Q}^n \mathcal{S}} 
(d, \mathbf{d}, v_0, v)
\times
{_{\sss QQ'} \mathcal{M}} 
(d, v_0, v-w) 
\times$ \\ $
{_{\sss Q'} \mathcal{M}} 
(\delta_{\sss QQ'} (d,v), 
\rho_{\sss QQ'} (v-w), 
\rho_{\sss QQ'} (v))
$}
\ar[dd]^-{\Id \times 
\theta^{-1}_{\sss QQ'}}
&
\\ \\ 
*\txt{$\bigsqcup$ 
\\
$\substack{
v_0 \in \mathbb{Z}_{\geq 0} [I]}$}
&
*\txt{$
{_{\sss Q}^n \mathcal{S}} 
(d, \mathbf{d}, v_0, v)
\times
{_{\sss Q} \mathcal{M}} 
(d, v_0, v) 
$}
\ar@{<-} '[r] '[ruuuuuuuu] '[uuuuuuuu]_-{
{_{\sss Q} \tau}_n}
&
}
\end{equation}

where $P$ is the permutation that moves all 
even terms in the direct product to the right.
Roughly speaking the commutativity of the above
diagram means that ``the Levi restriction 
commutes with the 
tensor product decomposition''.

\subsection{Digression:  
$\mathfrak{sl}_2$ case}
\label{SL2Case}

The Levi restriction bijections allow one
to reduce generic ADE case to the case 
of a quiver $R$ with one vertex and no edges
(which corresponds to $\mathfrak{sl}_2$).
This subsection contains a complete
description of quiver varieties 
and tensor product
varieties for the quiver $R$. 

Since $R$ has only one vertex and no edges
the vertex and edge indices are omitted
in the notation. Thus given two
(non graded) $\mathbb{C}$-linear spaces
$D$ and $V$ with dimensions
$d$ and $v$ respectively,
${_{\sss R} \Lambda}_{D,V}$ is the variety
of pairs $(p,q)$, where
$p \in \Hom_{\mathbb{C}} (D,V)$,
$q \in \Hom_{\mathbb{C}} (V,D)$,
and $pq =0$.
The open set
${_{\sss R} \Lambda}^{s}_{D,V}$ (resp.
${_{\sss R} \Lambda}^{\ast s}_{D,V}$)
consists of all 
$(p,q) \in {_{\sss R} \Lambda}_{D,V}$
such that $p$ is surjective
(resp. $q$ is injective).

The map $(p,q) \rightarrow qp \in
\End_{\mathbb{C}} (D)$ provides
an isomorphism between the quiver
variety 
${_{\sss R} \mathfrak{M}}^{s, *s}_{D,V}=
{_{\sss R} \Lambda}^{s, \ast s}_{D,V}/{G_V}$
(cf. \ref{Stability}) and the
$GL_D$-orbit in 
$\End_{\mathbb{C}} (D)$ consisting 
of all $t \in \End_{\mathbb{C}} (D)$
such that $t^2 =0$ and 
$\rank t = v$. Note that
this orbit is empty (if $2v > d$) 
or is a smooth connected quasi-projective
variety of 
dimension $2v(d-v)$ (if $2v \leq d$), 
which confers the 
corresponding statements in 
\ref{Stability}.

Similarly, the map 
$(p,q) \rightarrow (qp, \ker p)$ 
provides an isomorphism between the quiver
variety 
${_{\sss R} \mathfrak{M}}^{s} (d, v_0, v)=
{_{\sss R} \Lambda}^{s}_{D,V, v_0}/{G_V}$
(cf. \ref{MdvvSection})
and the variety of pairs
$(t, B)$, where 
$t \in \End_{\mathbb{C}} (D)$,
$t^2 =0$, $\rank t = v_0$,
and $B$ is a subspace of $D$ such that
$\im t \subset B \subset \ker t$, 
$\dim B = d -v$. It follows that 
${_{\sss R} \mathfrak{M}}^{s} (d, v_0, v)$
is smooth connected 
if $v_0 \leq v \leq d-v_0$ and 
empty otherwise, and hence
its set of irreducible components
${_{\sss R} \mathcal{M}} (d, v_0, v)$
is a one-element or the empty set
respectively.
Therefore
\begin{equation}\nonumber
{_{\sss R} \mathcal{M}} (d, v_0)=
\bigsqcup_{v \in \mathbb{Z}_{\geq 0}}
{_{\sss R} \mathcal{M}} (d, v_0, v)=
\bigsqcup_{v = v_0}^{d-v_0}
{_{\sss R} \mathcal{M}} (d, v_0, v)=
\bigsqcup_{v = v_0}^{d-v_0}
\{ {_{\sss R} \mathfrak{M}}^s 
(d, v_0, v) \} \; .
\end{equation}
Endow this set with the following structure 
of an $sl_2$-crystal:
\begin{equation}
\label{GL2CrystalMaps}
\begin{split}
{_{\sss R} wt} ({_{\sss R} \mathfrak{M}}^s 
(d, v_0, v) ) &= (d - 2v) \; ,
\\
{_{\sss R} \varepsilon}
({_{\sss R} \mathfrak{M}}^s 
(d, v_0, v)) &= v - v_0 \; ,
\\
{_{\sss R} \varphi ({_{\sss R}} 
\mathfrak{M}}^s
(d, v_0, v)) &= d - v - v_0  \; ,
\\ 
{_{\sss R} \Tilde{e}} (
{_{\sss R} \mathfrak{M}}^s 
(d, v_0, v)) &= 
\begin{cases} 
{_{\sss R} \mathfrak{M}}^s
(d, v_0, v-1)
&\text{ if $v > v_0$} \; , \\
0 &\text{ if $v \leq v_0$} \; ,
\end{cases}
\\
{_{\sss R} \Tilde{f}} (
{_{\sss R} \mathfrak{M}}^s 
(d, v_0, v)) &= 
\begin{cases} 
{_{\sss R} \mathfrak{M}}^s
(d, v_0, v+1)
&\text{ if $v < d - v_0$} \; , \\
0 &\text{ if $v \geq d - v_0$} \; .
\end{cases}
\end{split}
\end{equation}
Here the weight lattice of
$\mathfrak{sl}_2$ is identified
with $\mathbb{Z}$,
and the indexes of $\varepsilon$, $\varphi$,
$\Tilde{e}$, and $\Tilde{f}$ are omitted
because $\mathfrak{sl}_2$ has only one root.
The set ${_{\sss R} \mathcal{M}} (d, v_0)$
equipped with the above structure is a highest
weight normal $sl_2$-crystal with the highest
weight $(d - 2v_0)$. In other words, it is 
isomorphic (as a crystal) to 
$\mathcal{L} ((d-2v_0))$ (the crystal
of the canonical basis of the highest
weight irreducible representation
with highest weight $(d-2v_0)$). 

Let $D^1$ be a subspace of $D$,
$\mathbf{D}$ be the 2-step partial
flag $(0 \subset D^1 \subset D)$,
$\mathbf{d}$ be the pair 
$(\mathbf{d}^1, \mathbf{d}^2)$,
where
$\mathbf{d}^1 = \dim D^1$, 
$\mathbf{d}^2 = d - \dim D^1$,
and 
$\mathbf{v}=(\mathbf{v}^1, \mathbf{v}^2)$ 
be an element of
$\mathbb{Z} \oplus \mathbb{Z}$.
The map $(p,q) \rightarrow qp \in
\End_{\mathbb{C}} (D)$ provides
an isomorphism between the multiplicity
variety 
${_{\sss R}^2\mathfrak{S}}_{
D,\mathbf{D},V, \mathbf{v}}$
(cf. \ref{DefinitionOfTS}) and the
variety of all 
$t \in \End_{\mathbb{C}} (D)$ 
such that 
\begin{equation}\label{TForGL2}
\begin{split}
&t^2   = 0 \; ,
\\
&t D^1  \subset D^1 \; ,
\\
&\rank t   = v \; ,
\\ 
&\rank t|_{D^1} = \mathbf{v}^1 \; ,
\\ 
&\rank t|_{(D/D^1)} = \mathbf{v}^2 \; .
\end{split}
\end{equation}
Similarly, the map 
$(p,q) \rightarrow (qp, \ker p)$ 
provides an isomorphism between the 
tensor product variety
${_{\sss R}^2 \mathfrak{T}}_{
D, \mathbf{D}, V, \mathbf{v}}$
(cf. \ref{DefinitionOfTS})
and the variety of pairs
$(t, B)$, where 
$t \in \End_{\mathbb{C}} (D)$
satisfies conditions
\eqref{TForGL2}
and $B$ is a subspace of $D$ such that
$\im t \subset B \subset \ker t$, 
$\dim B = d - v$. 

A straightforward linear algebra
considerations 
(cf. \cite{Malkin2000a})
show that the variety of all 
$t \in \End_{\mathbb{C}} (D)$ satisfying
conditions \eqref{TForGL2} is
empty if either 
$v < \mathbf{v}^1 + \mathbf{v}^2$, or
$v > \mathbf{d}^2 - \mathbf{v}^2 + 
\mathbf{v}^1$, or 
$v > \mathbf{d}^1 - \mathbf{v}^1 + 
\mathbf{v}^2$, and
is a smooth connected 
quasi-projective variety otherwise.
Hence the RHS of the tensor
decomposition bijection
${_{\sss R} \tau_2}$
(cf. \ref{TensorDecomposition})
becomes
\begin{equation}\nonumber
\bigsqcup_{v_0 \in \mathbb{Z}_{\geq 0}}
{_{\sss R}^2 \mathcal{S}} (
d, \mathbf{d}, v_0, \mathbf{v})
\times 
{_{\sss R} \mathcal{M}} (d, v_0)
=
\bigsqcup_{v_0 = \mathbf{v}^1 + 
\mathbf{v}^2}^{
\min ( \mathbf{d}^2 - \mathbf{v}^2 + 
\mathbf{v}^1 , \mathbf{d}^1 - 
\mathbf{v}^1 + \mathbf{v}^2 )}
{_{\sss R} \mathcal{M}} (d, v_0) \; .
\end{equation}
The bijection 
\begin{equation}\nonumber
{_{\sss R} \tau_2} : \;
{_{\sss R} \mathcal{M}} 
(\mathbf{d}^1, \mathbf{v}^1)
\times
{_{\sss R} \mathcal{M}} 
(\mathbf{d}^2, \mathbf{v}^2)
\xrightarrow{\sim}
\bigsqcup_{v_0 = 
\mathbf{v}^1 + \mathbf{v}^2}^{
\min ( \mathbf{d}^2 - \mathbf{v}^2 + 
\mathbf{v}^1 , \mathbf{d}^1 - 
\mathbf{v}^1 + \mathbf{v}^2 )}
{_{\sss R} \mathcal{M}} (d, v_0) 
\end{equation}
can be written as
\begin{equation}\label{ActionOfTauGL2}
{_{\sss R} \tau_2} \bigl(
{_{\sss R} \mathfrak{M}} 
(\mathbf{d}^1, \mathbf{v}^1, \mathbf{u}^1),
{_{\sss R} \mathfrak{M}} 
(\mathbf{d}^2, \mathbf{v}^2, \mathbf{u}^2)
\bigr) =
{_{\sss R} \mathfrak{M}} 
(d, v_0, \mathbf{u}^1 + \mathbf{u}^2) \; ,
\end{equation}
where $v_0$ is described as follows.
Fix a subspace $B$ in $D$ such that
$\dim B = d - \mathbf{u}^1 - \mathbf{u}^2$ and
$\dim D^1 \cap B = \mathbf{d}^1 - \mathbf{u}^1$. 
Then $v_0$ is the (unique) integer such
that the set of operators 
$t \in \End_{\mathbb{C}} (D)$ 
with $\rank t = v_0$ form an open
subset in the set of all $t$
satisfying the following conditions:
\begin{equation}\label{TForTauGL2}
\begin{split}
&t^2   = 0 \; ,
\\
&t D^1  \subset D^1 \; ,
\\ 
&\rank t|_{D^1} = \mathbf{v}^1 \; ,
\\ 
&\rank t|_{(D/D^1)} = \mathbf{v}^2 \; ,
\\
&\im t \subset B \subset \ker t \; . 
\end{split}
\end{equation}
Again elementary linear algebra
(cf. \cite{Malkin2000a})
shows that
\begin{equation}\label{VZero}
v_0 = \min (
\mathbf{u}^2 + \mathbf{v}^1 , 
\mathbf{d}^1 - \mathbf{u}^1 + 
\mathbf{v}^2 ) \; .
\end{equation}
The equality \eqref{VZero} 
together with definition
of the tensor product of
crystals (cf. \ref{CrystalTensorProduct})
imply the following theorem 
(cf. \cite{Malkin2000a}).
\begin{theorem}
The bijection ${_{\sss R} \tau_2}$
is a crystal isomorphism
\begin{equation}\nonumber
{_{\sss R} \tau_2} : \quad
{_{\sss R} \mathcal{M}} 
(\mathbf{d}^1, \mathbf{v}^1)
\otimes
{_{\sss R} \mathcal{M}} 
(\mathbf{d}^2, \mathbf{v}^2)
\xrightarrow{\sim}
\bigoplus_{v_0 \in \mathbb{Z}_{\geq 0}}
{_{\sss R}^2 \mathcal{S}} (
d, \mathbf{d}, v_0, \mathbf{v})
\otimes 
{_{\sss R} \mathcal{M}} (d, v_0) \; ,
\end{equation}
where $\mathbf{d} = 
(\mathbf{d}^1 , \mathbf{d}^2)$,
$\mathbf{v} =
(\mathbf{v}^1 , \mathbf{v}^2 )$,
$d = \mathbf{d}^1 + \mathbf{d}^2$,
and the set 
${_{\sss R}^2 \mathcal{S}} (
d, \mathbf{d}, v_0, \mathbf{v})$
is considered as a trivial 
$\mathfrak{sl}_2$-crystal.
\end{theorem}
\begin{remark}
In \cite{Malkin2000a} definitions
of tensor product and multiplicity
varieties for $\mathfrak{sl}_2$ are
slightly different. Namely the varieties
in \cite{Malkin2000a} are fibrations over
the Grassmannian of all subspaces $D^1$ in
$D$ with given dimension, and the fibers of
these fibrations are the corresponding 
(tensor product or multiplicity)
varieties as they defined in this paper.
The sets of irreducible components and
the tensor decomposition bijection are
the same here and in \cite{Malkin2000a}. 
\end{remark}

\subsection{A crystal structure
on ${_{\sss Q} \mathcal{M}} (d, v_0)$}
\label{CrystalMDV0}

In this subsection it is shown that
the set 
${_{\sss Q} \mathcal{M}} (d, v_0)$
(cf. \ref{DefinitionOfMCal}) can be
endowed with a structure of 
$\mathfrak{g}$-crystal. 
This structure was introduced by
Nakajima \cite{Nakajima1998}
following an idea of Lusztig
\cite[Section 12]{Lusztig1991a}.

The first step is to define the 
weight function. Identify 
the weight lattice of 
$\mathfrak{g}$ with $\mathbb{Z}[I]$
(i.e. 
$i \in I \subset \mathbb{Z} [I]$
is the $i$-th fundamental weight).
Then the weight function
\begin{equation}\nonumber
{_{\sss Q} wt} : \; 
{_{\sss Q} \mathcal{M}} (d, v_0)
\rightarrow
\mathbb{Z} [I]
\end{equation}
is given by
\begin{equation}\nonumber
{_{\sss Q} wt} (Z) = 
d - 2v + {_{\sss Q} X} v 
\end{equation}
for $Z \in 
{_{\sss Q} \mathcal{M}} (d, v_0, v)
\subset
{_{\sss Q} \mathcal{M}} (d, v_0)$.
Here ${_{\sss Q} X}$ 
is the matrix defined in
\eqref{DefinitionOfX}.
The crucial property of the
function ${_{\sss Q} wt}$
is its good behavior with respect
to the Levi restriction.
Namely
\begin{equation}\nonumber
\rho_{\sss QQ'} 
(d - 2v + {_{\sss Q} X} v) =
\delta_{\sss QQ'} (d, v) -
2 \rho_{\sss QQ'} (v) +
{_{\sss Q'} X} 
\rho_{\sss QQ'} (v) \; ,
\end{equation}
where 
$\delta_{\sss QQ'}$ and 
$\rho_{\sss QQ'}$ are as in
\ref{DefinitionOfTheta}.
This allows one to use Levi
restriction to $\mathfrak{sl}_2$
subalgebras to define the crystal
structure on 
${_{\sss Q} \mathcal{M}} (d, v_0)$.
More explicitly, let
$R_i$ be the subquiver of $Q$
consisting of the vertex $i$ and
no edges. Recall 
(cf. \ref{DefinitionOfTheta}) 
that there is
a bijection 
\begin{equation}\label{ThetaRi}
\theta_{\sss QR_i} : \;
{_{\sss Q} \mathcal{M}}
(d, v_0) 
\xrightarrow{\sim}
\bigsqcup_{
v \in \mathbb{Z}_{\geq 0} [I]}
{_{\sss QR_i} \mathcal{M}}
(d, v_0, v)
\times
{_{\sss R_i} \mathcal{M}}
(d_i + ({_{\sss Q} X} v)_i, v_i) \; .
\end{equation}
Consider the RHS of 
\eqref{ThetaRi} as 
an $\mathfrak{sl}_2$-crystal with the 
crystal structure coming from the 
second multiple (on which it is defined
as in \ref{SL2Case}), and let
\begin{equation}\label{EpsilonPhiEF}
\begin{split}
{_{\sss Q} \varepsilon}_i &=
{_{\sss R_i} \varepsilon} \circ
\theta_{\sss QR_i} \; ,
\\
{_{\sss Q} \varphi}_i &=
{_{\sss R_i} \varphi} \circ
\theta_{\sss QR_i} \; ,
\\
{_{\sss Q} \Tilde{e}}_i &=
\theta_{\sss QR_i}^{-1} \circ
{_{\sss R_i} \Tilde{e}} \circ
\theta_{\sss QR_i} \; ,
\\
{_{\sss Q} \Tilde{f}}_i &=
\theta_{\sss QR_i}^{-1} \circ
{_{\sss R_i} \Tilde{f}} \circ
\theta_{\sss QR_i} \; .
\end{split}
\end{equation}
These formulas together with the
weight function ${_{\sss Q} wt}$
provide a structure of 
$\mathfrak{g}$-crystal on the set
${_{\sss Q} \mathcal{M}} (d, v_0)$.
By abuse of notation this crystal
is also denoted by 
${_{\sss Q} \mathcal{M}} (d, v_0)$.

\begin{proposition}
\begin{alphenum}
\item\label{MEmptyCrystal}
If $d - 2v_0 + {_{\sss Q} X} v_0 \notin
\mathbb{Z}_{\geq 0} [I]$ then
${_{\sss Q} \mathcal{M}} (d, v_0)$
is an empty set.
\item\label{MHighestWeightCrystal}
If $d - 2v_0 + {_{\sss Q} X} v_0 \in
\mathbb{Z}_{\geq 0} [I]$ then
${_{\sss Q} \mathcal{M}} (d, v_0)$ is
a highest weight normal
$\mathfrak{g}$-crystal with the highest
weight $d - 2v_0 + {_{\sss Q} X} v_0$.
\item\label{MMPrime}
If $d - 2v_0 + {_{\sss Q} X} v_0 =
d' - 2v'_0 + {_{\sss Q} X} v'_0$ then
the crystals
${_{\sss Q} \mathcal{M}} (d, v_0)$
and
${_{\sss Q} \mathcal{M}} (d', v'_0)$
are isomorphic.
\end{alphenum}
\end{proposition}
\begin{proof}
\ref{MEmptyCrystal} follows
from \ref{LambdaSAstSEmpty}.

\ref{MHighestWeightCrystal} follows 
directly from definitions and 
the fact that 
${_{\sss Q} \mathfrak{M}}^s 
(d, v_0, v_0)$ 
(a connected smooth variety)
is the only element of 
${_{\sss Q} \mathcal{M}} (d, v_0)$
which is killed by 
${_{\sss Q} \Tilde{e}}_i$
for all $i \in I$. 
To prove this let $Z \in 
{_{\sss Q} \mathcal{M}} (d, v_0)$
be an irreducible component
of ${_{\sss Q} \Lambda}^s_{D, V, U}$,
and $(x,p,q)$ be a
generic point of $Z$.
Let $\gamma$ be as in 
\ref{DefinitionOfMCal} and
$(x^{\sss UU},
(x^{\sss TT}, p^{\sss TD} ,
q^{\sss DT})) = \gamma ((x, p, q))$.
If $U \neq \{ 0 \}$ then there exists
$i \in I$ such that 
$\bigcap_{\substack{h \in H \\ \Out (h) =i}} 
\ker x^{UU}_h \neq \{ 0 \}$
(cf. \cite[12]{Lusztig1991a}).
It follows that 
$(\bigcap_{\substack{h \in H \\ \Out (h) =i}} 
\ker x_h) \bigcap \ker q_i \neq \{ 0 \}$, and
hence 
${_{\sss Q} \Tilde{e}}_i Z \neq 0$.
Therefore if ${_{\sss Q} \Tilde{e}}_i Z = 0$
for all $i \in I$ then $U = \{ 0 \}$, which
means
$Z = {_{\sss Q} \mathfrak{M}}^s (d, v_0, v_0)$.

To prove \ref{MMPrime} note that as sets
both ${_{\sss Q} \mathcal{M}} 
(d, v_0, v_0 + u)$
and ${_{\sss Q} \mathcal{M}} 
(d', v'_0, v'_0 +u)$ 
are in natural bijections 
(induced by the vector bundles $\gamma$ -- 
cf. \ref{DefinitionOfMCal})
with the set of irreducible components 
of the variety
${^\delta \Lambda}_U$, where $\dim U = u$, and
$\delta = d - 2 v_0 + {_{\sss Q} X} v_0 =
d' - 2 v'_0 + {_{\sss Q} X} v'_0$. 
In this way one obtains a bijection between
${_{\sss Q} \mathcal{M}} (d, v_0, v_0 + u)$
and 
${_{\sss Q} \mathcal{M}} (d', v'_0, v'_0 +u)$,
and it follows from definitions that 
this bijection is a crystal isomorphism.
\end{proof}

\subsection{The main theorem}
\label{MainTheoremForQ}
From now on the quiver $Q$ is fixed
and thus omitted in the notation.
The following is the main result of
this paper.

\begin{theorem}
The bijection $\tau_n$
(cf. \ref{TensorDecomposition})
is a crystal isomorphism
\begin{equation}\nonumber
\tau_n : \quad
\mathcal{M}
(\mathbf{d}^1, \mathbf{v}^1)
\otimes \ldots \otimes
\mathcal{M}
(\mathbf{d}^n, \mathbf{v}^n)
\xrightarrow{\sim}
\bigoplus_{v_0 \in 
\mathbb{Z}_{\geq 0} [I]}
{^n \mathcal{S}} (
d, \mathbf{d}, v_0, \mathbf{v})
\otimes 
\mathcal{M} (d, v_0) \; ,
\end{equation}
where the set 
${^n \mathcal{S}} (
d, \mathbf{d}, v_0, \mathbf{v})$
is considered as a trivial 
$\mathfrak{g}$-crystal.
\end{theorem}
Note that both sides might be 
empty.
\begin{proof}
Because of commutativity of 
diagram \eqref{TauEta} it is
enough to consider the case
$n=2$, and because the definition
of the crystal structure uses
Levi restriction which commutes
with tensor product 
(cf. \eqref{TauTheta}) 
it is enough to consider
the case 
$\mathfrak{g}=\mathfrak{sl}_2$.
Now the theorem follows from
Theorem \ref{SL2Case}.
\end{proof}

\subsection{Corollary}
\label{Corollary}

Recall that $L (\mu )$ denotes the
highest weight irreducible 
representation of
$\mathfrak{g}$ with highest weight
$\mu \in \mathbb{Z}_{\geq 0} [I]$,
and $\mathcal{L} (\mu)$ denotes
the crystal of the canonical basis
of $L (\mu)$.

The following proposition is a corollary
of Theorem \ref{MainTheoremForQ}.

\begin{proposition}
\indent\par
\begin{alphenum}
\item\label{SaitoCorollary}
The $\mathfrak{g}$-crystal
$\mathcal{M} (d, v_0)$
is isomorphic to
$\mathcal{L} (\mu)$, where
$\mu = d - 2v_0 + X v_0$.
\item\label{ADEHallCorollary}
The cardinal of the set
${^n \mathcal{S}} (
d, \mathbf{d}, v_0, \mathbf{v})$
of irreducible components of
the multiplicity variety
${^n \mathfrak{S}} (
d, \mathbf{d}, v_0, \mathbf{v})$
is equal to the multiplicity of
$L (\mu)$ in
$L (\mu^1) \otimes \ldots
\otimes L (\mu^n)$, where
$\mu = d - 2v_0 + X v_0$,
$\mu^k = \mathbf{d}^k - 
2\mathbf{v}^k + 
X \mathbf{v}^k$.
In other words,
\begin{equation}\nonumber
| {^n \mathcal{S}} (
d, \mathbf{d}, v_0, \mathbf{v}) | =
\dim_{\mathbb{C}} \Hom_{\mathfrak{g}}
(L (\mu),
L (\mu^1) \otimes \ldots
\otimes L (\mu^n)) \; .
\end{equation}
\item
\label{CorollaryTensorVariety}
The bijection $\alpha_n$ 
(cf. \ref{DefinitionOfAlpha}) 
identifies the set
${^n \mathcal{T}} (
d, \mathbf{d}, v, \mathbf{v})$
of irreducible components of
the tensor product variety
${^n \mathfrak{T}} (
d, \mathbf{d}, v, \mathbf{v})$
with the weight subset of
weight $d - 2v + X v$
in the crystal
$\mathcal{L} (\mu^1) \otimes \ldots
\otimes \mathcal{L} (\mu^n)$,
where $\mu^1, \ldots , \mu^n$ are
as in \ref{ADEHallCorollary}.
\end{alphenum}
\end{proposition} 
\begin{proof}
The proposition follows from
Theorem \ref{MainTheoremForQ},
Theorem \ref{JosephTheorem},
Proposition \ref{CrystalMDV0},
and Proposition \ref{SNonEmpty}.
\end{proof}
\begin{remark}
Note that \ref{ADEHallCorollary} 
is a generalization of a theorem
due to Hall \cite{Hall1959}
(cf. 
\cite[Chapter II]{Macdonald1995}).
Statement \ref{SaitoCorollary}
was also proven
by Saito using
results of \cite{KashiwaraSaito}.
\end{remark}

\subsection{The extended Lie algebra 
$\mathfrak{g}'$}
\label{GPrimeSection}

Proposition \ref{MMPrime} shows
that in the set of 
$\mathfrak{g}$-crystals
$\{ \mathcal{M} (d, v_0) \}_{
d, v_0 \in \mathbb{Z}_{\geq 0} [I]}$
there are isomorphic elements, 
which suggests that in the context of quiver
varieties it is more natural to
consider a central 
extension of $\mathfrak{g}$ than
$\mathfrak{g}$ itself.

Let $\mathfrak{g}' = \mathfrak{g}
\oplus \mathfrak{t}$, where 
$\mathfrak{t}$ is a Cartan subalgebra
of $\mathfrak{g}$. Then $\mathfrak{g}'$
is a reductive Lie algebra. Identify
the weight lattice of $\mathfrak{g}'$
with $\mathbb{Z} [I] \oplus 
\mathbb{Z} [I]$ in such a way that
the projection $\mathbb{Z}[I] \oplus
\mathbb{Z} [I] \rightarrow
\mathbb{Z} [I]$ given by 
$(v,u) \rightarrow v - u$ is the
projection from the weight lattice 
of $\mathfrak{g}'$ onto the weight
lattice of $\mathfrak{g}$.
A weight $(v,u) \in
\mathbb{Z} [I] \oplus \mathbb{Z} [I]$
is called \emph{integrable} if
$u \in \mathbb{Z}_{\geq 0} [I]$
and $v-u \in \mathbb{Z}_{\geq 0} [I]$.
Let $\mathcal{Q}^+ \subset 
\mathbb{Z}_{\geq 0} [I] \oplus 
\mathbb{Z}_{\geq 0} [I]$ 
be the set of all integrable weights. 
A finite dimensional representation of
$\mathfrak{g}'$ is called 
\emph{integrable} if the 
highest weights of all its
irreducible components are 
integrable.
Note that the category of 
finite dimensional integrable
representations is closed with respect
to tensor products.

One can endow the set 
$\mathcal{M} (d,v_0)$ with a structure
of $\mathfrak{g}'$-crystal as follows.
The maps $\varepsilon_i$, $\varphi_i$,
$\Tilde{e}_i$, and $\Tilde{f}_i$
are given by \eqref{EpsilonPhiEF}
(note that roots of $\mathfrak{g}'$
are roots of $\mathfrak{g}$).
The weight function is given by
(cf. \ref{CrystalMDV0}) 
\begin{equation}\label{WeightPrime}
wt' (Z) = (d - v + X v , v) \in
\mathbb{Z} [I] \oplus \mathbb{Z} [I] 
\end{equation}
for $Z \in 
\mathcal{M} (d, v_0, v)
\subset
\mathcal{M} (d, v_0)$.
Note that this weight function
commutes with the Levi restriction
(cf. an analogous statement for the
weight function $wt$ in
\ref{CrystalMDV0}).

The following theorem shows that 
non-empty elements of
the set 
\begin{equation}\nonumber
\{ \mathcal{M} (d, v_0) \}_{
(d, v_0) \in \mathbb{Z}_{\geq 0} [I]
\oplus \mathbb{Z}_{\geq 0} [I]}
\end{equation}
form a closed (with respect to
tensor product) family of
$\mathfrak{g}'$-crystals,
which is isomorphic to the 
family of crystals of canonical 
bases of irreducible integrable
representations of $\mathfrak{g}'$.

\begin{theorem}
\begin{alphenum}
\item\label{MDVEmptyPrime}
If $(d - v_0 + Xv_0, v_0) 
\notin \mathcal{Q}^+$
then the set $\mathcal{M} (d, v_0)$
is empty.
\item\label{MDVCrystalPrime}
If $(d - v_0 + X v_0 , v_0) 
\in \mathcal{Q}^+$
then the set $\mathcal{M} (d, v_0)$
together with the weight
function $wt'$ and 
the maps $\varepsilon_i$, $\varphi_i$,
$\Tilde{e}_i$, and $\Tilde{f}_i$
given by \eqref{EpsilonPhiEF} is
a highest weight normal
$\mathfrak{g}'$-crystal with the
highest weight $(d - v_0 + X v_0, v_0)$.
\item\label{PositiveWeightsPrime}
The image of the weight function
$wt' : \; \mathcal{M} (d, v_0)
\rightarrow 
\mathbb{Z} [I] \oplus \mathbb{Z} [I]$
is contained in 
$\mathbb{Z}_{\geq 0} [I] \times
\mathbb{Z}_{\geq 0} [I]$.
\item\label{TauCrystalPrime}
The bijection $\tau_n$
(cf. \ref{TensorDecomposition})
is an isomorphism
of $\mathfrak{g}'$-crystals
\begin{equation}\nonumber
\tau_n : \;
\mathcal{M}
(\mathbf{d}^1, \mathbf{v}^1)
\otimes \ldots \otimes
\mathcal{M}
(\mathbf{d}^n, \mathbf{v}^n)
\xrightarrow{\sim}
\bigoplus_{v_0 \in 
\mathbb{Z}_{\geq 0} [I]}
{^n \mathcal{S}} (
d, \mathbf{d}, v_0, \mathbf{v})
\otimes 
\mathcal{M} (d, v_0)
\end{equation}
where the set 
${^n \mathcal{S}} (
d, \mathbf{d}, v_0, \mathbf{v})$
is considered as a trivial 
$\mathfrak{g}'$-crystal.
\item\label{MFamilyPrime}
The family 
$\{ \mathcal{M} (d, v_0) \}$
of $\mathfrak{g}'$-crystals
is isomorphic to the family 
$\{ \mathcal{L} (d , v_0) \}$
of crystals of canonical bases
of irreducible integrable
representations of $\mathfrak{g}'$.
\end{alphenum}
\end{theorem}
In \ref{MFamilyPrime} 
$\mathcal{L} (d , v_0)$ denotes
the crystal of the canonical
basis of the irreducible
integrable representation
$L (d , v_0)$ with the highest
weight $(d - v_0 + X v_0, v_0)$.
\begin{proof}
Proofs of 
\ref{MDVEmptyPrime},
\ref{MDVCrystalPrime},
\ref{TauCrystalPrime},
and \ref{MFamilyPrime} are 
analogous to the proofs
of Propositions 
\ref{MEmptyCrystal},
\ref{MHighestWeightCrystal},
Theorem \ref{MainTheoremForQ}, 
and Corollary \ref{SaitoCorollary},
respectively,
or one can deduce the former statements
from the latter.

To prove \ref{PositiveWeightsPrime}
note that an element
of the set 
$\mathcal{M} (d, v_0, v)$ is
an irreducible component
of the variety
$\mathfrak{M}^s_{D, V, v_0}$
(cf. \ref{MdvvSection}),
where $\dim D = d$, $\dim V = v$.
A point of this variety is a 
stable triple $(x,p,q)$. Since
it is stable
$\im p_i + \sum_ {\substack{
h \in H \\ \In (h) = i}}
\im x_h = V_i$. Therefore
$v_i \leq (d + Xv)_i$, which implies 
\ref{PositiveWeightsPrime}.
\end{proof}

\begin{remark}
It follows from 
\ref{PositiveWeightsPrime}
that all weights of a
finite dimensional
integrable representation of
$\mathfrak{g}'$ belong to
$\mathbb{Z}_{\geq 0} [I] \times
\mathbb{Z}_{\geq 0} [I]$.
\end{remark}

\appendix

\section{Another description of
multiplicity varieties and the 
tensor product diagram}

\subsection{Lusztig's description
of the variety 
$\mathfrak{M}^{s, \ast s}_{D,V}$}
\label{AnotherQuiver}

Lusztig (cf. \cite{Lusztig1998, Lusztig2000a,
Lusztig2000b}) identified varieties 
$\mathfrak{M}^{s, \ast s}_{D,V}$
(for a fixed $D$ and varying  $V$)
with locally closed subsets in a
certain variety $Z_D$. This subsection
contains an overview of the 
Lusztig's construction (see the
above references for more details).

Recall (cf. \ref{SL2Case}) that in the
$\mathfrak{sl}_2$ case the map 
$(p,q) \rightarrow qp$ allows one to
identify 
$\mathfrak{M}^{s, \ast s}_{D, V}$
with a nilpotent orbit in 
$\End_{\mathbb{C}} D$ (namely with
the orbit consisting of all 
operators $t$ such that $t^2 =0$,
and  $\rank t = \dim V$). Lusztig's
construction is the direct 
generalization of this device. 
However for a generic quiver one has
to take care not only of the maps
$p$ and $q$, but also of $x$ 
(i.e one has to include the path 
algebra $\mathcal{F}$ into the play).

Let $\mathbb{C}^I = 
\oplus_{\sss i \in I}
\mathbb{C}$ considered as a
semi-simple $\mathbb{C}$-algebra.
Note that a $\mathbb{Z} [I]$-graded
$\mathbb{C}$-linear space is the same 
as a $\mathbb{C}^I$-module, and
that the path algebra $\mathcal{F}$
is a $\mathbb{C}^I$--$\mathbb{C}^I$
algebra. Let 
$\Tilde{\mathcal{F}}$
be a $\mathbb{C}$-algebra defined as 
follows. As a $\mathbb{C}$-linear
space $\Tilde{\mathcal{F}} =
\mathcal{F} \oplus 
\oplus_{\sss i \in I}
\mathbb{C} u_i$, and the 
multiplication $\circ$ in 
$\Tilde{\mathcal{F}}$ is given by
\begin{equation}\nonumber
\begin{split}
f \circ f' &= \sum_{i \in I}
f \cdot \theta_i \cdot f' \; , 
\\
f \circ u_i &= f \cdot [i] \; ,
\\
u_i \circ f &= [i] \cdot f \; ,
\\
u_i \circ u_j &= 
\delta_{ij} \; u_i \; , 
\end{split}
\end{equation}
where $f, f' \in \mathcal{F}$,
$i, j \in I$, $\theta_i$ is as in 
\ref{DefinitionOfThetaI},
$[i] \in \mathcal{F}$ 
denotes the path of length $0$
starting and ending at a vertex 
$i \in I$, 
and $\cdot$ denotes the multiplication 
in the path algebra $\mathcal{F}$.
Multiplication by $u_i$ on the
left and on the right endows 
$\Tilde{\mathcal{F}}$ with a
structure of $\mathbb{C}^I$--$\mathbb{C}^I$
algebra. This algebra was introduced
by Lusztig \cite[2.1]{Lusztig2000a}
(see also \cite[2.4]{Lusztig1998} and
\cite[2.1]{Lusztig2000b}).
It is an associative algebra 
with unit (given by
$\sum_{i \in I} u_i$).
Since the Dynkin graph is
assumed to be of finite type
$\Tilde{\mathcal{F}}$ is
finitely generated (cf.
\cite[Lemma 2.2]{Lusztig2000a}).

Let $D$ be a $\mathbb{C}^I$-module,
and $Z_D$ be the set of all 
$\mathbb{C}^I$--$\mathbb{C}^I$ algebra 
homomorphisms $\pi : \; \Tilde{\mathcal{F}}
\rightarrow \End_{\mathbb{C}^I} (D)$
(in other words, the set of all
representations of $\Tilde{\mathcal{F}}$
in $D$).
The set $Z_D$ is naturally an affine
variety. Let $\vartheta'$ be a
regular map
\begin{equation}\nonumber
\vartheta' : \quad
\Lambda^{s,\ast s}_{D,V} 
\rightarrow
Z_D
\end{equation}
given by
\begin{equation}\nonumber
\vartheta' ((x,p,q)) = \pi \; ,
\end{equation}
where
\begin{equation}\nonumber
\pi ([h_1 \ldots h_n]) = 
q _{\sss \Out (h_1)} x_{h_1}
\ldots x_{h_n} p_{\sss \In (h_n)}
\end{equation}
for a path $[h_1 \ldots h_n]$ in
$\mathcal{F}$.
The map $\vartheta'$ being constant 
on the orbits of $G_V$, it induces
a regular map $\vartheta : \;
\mathfrak{M}^{s, \ast s}_{D,V} =
\Lambda^{s,\ast s}_{D,V} /G_V 
\rightarrow Z_D$.
Let $Z_{D, V} = \vartheta (
\mathfrak{M}^{s, \ast s}_{D,V})
\subset Z_D$. One has $Z_{D, V} =
Z_{D, V'}$ if $\dim V = \dim V'$ and
thus the notation $Z_{D, v}$ (where
$v = \dim V$) is sometimes used
instead of $Z_{D,V}$. 

The following 
theorem is due to Lusztig.
\begin{theorem}\indent\par
\begin{alphenum}
\item
\cite[Theorem 5.5]{Lusztig1998}
\cite[Lemma 4.12c]{Lusztig2000b}
The map 
$\vartheta : \;
\mathfrak{M}^{s, \ast s}_{D,V}
\rightarrow Z_{D, V}$ is a homeomorphism
both in Zariski and in the smooth 
topologies.
\item
\cite[Lemma 4.12c]{Lusztig2000b}
The set $Z_{D,V}$ is locally closed in 
$Z_D$.
\item
\cite[Lemma 4.12d]{Lusztig2000b}
$Z_D = \bigsqcup_{v \in 
\mathbb{Z}_{\geq 0}[I]} Z_{D,v}$. 
\end{alphenum}
\end{theorem}

The above theorem shows that varieties 
$\mathfrak{M}^{s, \ast s}_{D,V}$ 
for various $V$ 
can be ``glued together'' which is
crucial for a geometric construction
of the tensor product. Nakajima also
considered a union of all
$\mathfrak{M}^{s, \ast s}_{D,V}$ 
for a fixed $D$ 
(cf. \cite{Nakajima1998},
\cite[2.5]{Nakajima2001}). 
He calls the resultant variety 
$\mathfrak{M}_0 (\infty, d)$, 
where $d = \dim D$.
It follows from 
\cite[Theorem 5.5]{Lusztig1998} that
$\mathfrak{M}_0 (\infty, d) = Z_D$.

Given $\pi \in Z_D$ one can find
$v \in \mathbb{Z}_{\geq 0} [I]$ such that
$\pi \in Z_{D, v}$ as follows
(cf. \cite[Section 2]{Lusztig1998}).
Let $\Dot{D} = \mathcal{F} 
\otimes_{\sss\mathbb{C}^I} D$.
Then $\Dot{D}$ is naturally a left
$\mathcal{F}$-module.
Let $\varpi_{\pi} \in
\Hom_{\mathbb{C}^I} (\Dot{D}, D)$ be 
given by 
$\varpi_{\pi} (f \otimes d) = 
\pi (f) d$, and let $\mathcal{K}_{\pi}$  
be the largest $\mathcal{F}$-submodule
of $\Dot{D}$ contained in the kernel
of $\varpi_{\pi}$. Then
$v = \dim (\Dot{D}/\mathcal{K}_{\pi})$
is finite and $\pi \in Z_{D, v}$.

\subsection{Multiplicity varieties}
Let $D$, $V$ be $\mathbb{C}^I$-modules,
and $\mathbf{D} = 
(\{ 0 \} = \mathbf{D}^0 \subset
\mathbf{D}^1 \subset \ldots \subset
\mathbf{D}^n = D)$ be an $n$-step
$\mathbb{C}^I$-filtration of $D$.
Let $\mathbf{v} 
\in (\mathbb{Z}_{\geq 0} [I])^n$, 
and $v = \dim V$.
Recall (cf. \ref{DefinitionOfTS})
that the multiplicity variety 
$\mathfrak{S}_{D, \mathbf{D}, V, 
\mathbf{v}}$ is a subset of
$\mathfrak{M}^{s, \ast s}_{D,V}$.
Consider the restriction of the map
$\vartheta : \; 
\mathfrak{M}^{s, \ast s}_{D, V}
\rightarrow Z_D$ 
(cf. Appendix \ref{AnotherQuiver})
to $\mathfrak{S}_{D, \mathbf{D}, V, 
\mathbf{v}}$.
It follows from the definition
of the multiplicity variety 
that this restriction 
provides a homeomorphism 
between a multiplicity variety
$\mathfrak{S}_{D, \mathbf{D}, V, 
\mathbf{v}}$ and a subvariety
$Z_{D, \mathbf{D}, v, \mathbf{v}}$
of $Z_D$ consisting of all
$\pi \in Z_D$ such that
\begin{itemize}
\item
for any $k =1 , \ldots , n$
the $\mathbb{C}^I$-submodule
$\mathbf{D}^k$ is an 
$\Tilde{\mathcal{F}}$-submodule
of $D$ with respect to the 
representation $\pi :
\Tilde{\mathcal{F}} \rightarrow
\End_{\mathbb{C}^I} D$,
\item
$\pi \in Z_{D,v}$,
\item 
$\pi |_{\sss \mathbf{D}^k /
\mathbf{D}^{k-1}}
\in Z_{\mathbf{D}^k /\mathbf{D}^{k-1},
\mathbf{v}^k} \subset 
Z_{\mathbf{D}^k /\mathbf{D}^{k-1}}$ 
for any $k =1 , \ldots , n$.
\end{itemize}
This description of the multiplicity
varieties is reminiscent of the definition
of the Hall-Ringel algebra (cf. 
\cite{Ringel1988}) associated to 
the algebra $\Tilde{\mathcal{F}}$.
One should be careful however because
the set $Z_{D,v} \subset Z_D$ is
in general a union of orbits of
$\Aut_{\mathbb{C}^I} (D)$, rather
than a single orbit.
Similar situation occurs in the 
geometric (quiver) construction 
of the positive part of a (quantum) 
universal enveloping algebra when the 
underlying quiver is not of
finite type 
(cf. \cite{Lusztig1991a, Schofield}).

\subsection{The tensor product
diagram}
Let $D$, $D^1$, $\ldots$ , $D^n$ be
$\mathbb{C}^I$-modules. 
By analogy with Lusztig's construction
of the canonical basis one can 
consider the following diagram
(cf. \cite[6.1a]{Lusztig1990}).
\begin{equation}
\label{TensorProductDiagram}
\xymatrix{
&
Z' 
\ar[r]^-{p_2}
\ar[ld]_-{p_1} 
&
Z'' 
\ar[rrd]^-{p_3}
\\
Z_{D^1}
\times \ldots \times Z_{D^n}
&&&&
Z_D 
}
\end{equation}
Here the notation is as follows:

$Z''$ is the variety of all pairs
$(\pi , \mathbf{W})$ consisting of
$\pi \in Z_D$ and a 
$\mathbb{C}^I$-filtration 
$\mathbf{W} = (0 = 
\mathbf{W}^0 \subset \mathbf{W}^1 
\subset \ldots \subset 
\mathbf{W}^n = D)$ 
such that $\pi (f) \mathbf{W}^k 
\subset \mathbf{W}^k$
for any $f \in \mathcal{F}$,
$k = 1, \ldots , n$, 
and $\dim (\mathbf{W}^k / 
\mathbf{W}^{k-1})= \dim D^k$,

$Z'$ is the variety of all triples
$(\pi , \mathbf{W}, \mathbf{R})$
where $(\pi , \mathbf{W}) \in
Z''$, and $\mathbf{R}$ is an
$n$-tuple $(\mathbf{R}^1 , \ldots ,
\mathbf{R}^n)$ consisting 
of $\mathbb{C}^I$-isomorphisms
$\mathbf{R}^k : \; 
D^k \xrightarrow{\sim}
\mathbf{W}^k /\mathbf{W}^{k-1}$
for $k = 1 , \ldots , n$,

$p_2 ((\pi, \mathbf{W}, \mathbf{R}))=
(\pi , \mathbf{W})$,

$p_3 ((\pi , \mathbf{W})) = \pi$, 

$p_1 ((\pi , \mathbf{W}, \mathbf{R})) =
(\pi^1 , \ldots , \pi^n)$,
where
$\pi^k = (\mathbf{R}^k)^{-1} 
(\pi |_{\sss \mathbf{W}^k / 
\mathbf{W}^{k-1}})
\mathbf{R}^k$ for any 
$k = 1 , \ldots ,n$.

The diagram 
\eqref{TensorProductDiagram} is 
closely related to tensor product
and multiplicity varieties. In
particular a subset 
\begin{equation}\label{HallSubset}
p_2 \circ p_1^{-1}
(Z_{D^1, \mathbf{v}^1} 
\times \ldots \times
Z_{D^n, \mathbf{v}^n})
\bigcap
p_3^{-1} (Z_{D,V}) 
\subset Z''
\end{equation}
is the total space of a fibration 
over the space of 
$\mathbb{C}^I$-filtrations 
of $D$ with dimensions of the 
subfactors given by
$\dim D^1 , \ldots , \dim D^n$,
and the fiber of this fibration
over a point $\mathbf{D}$
is isomorphic to the multiplicity
variety ${^n \mathfrak{S}}_{
D, \mathbf{D}, V, \mathbf{v}}$.
It follows that the subset 
\eqref{HallSubset} has pure
dimension and the number of
its irreducible components 
is equal to 
\begin{equation}\nonumber
\dim_{\mathbb{C}} \Hom_{\mathfrak{g}'}
\bigl( L ( d,v ),
L ( \mathbf{d}^1 , \mathbf{v}^1 )
\otimes \ldots \otimes 
L ( \mathbf{d}^n , 
\mathbf{v}^n ) \bigr) \; ,
\end{equation}
where $d = \dim D$, 
$\mathbf{d}^k = \dim D^k$,
$v = \dim V$, and
$L ( d,v )$ denotes the highest weight
irreducible representation of 
$\mathfrak{g}'$ with the highest
weight $(d - v + Xv , v)$
(cf. \ref{GPrimeSection}).
Tensor product varieties
also can be described in the context
of the diagram 
\eqref{TensorProductDiagram}. 
They are related to a resolution 
of singularities of $Z''$. 

Let $Z = \bigsqcup Z_D$, where
$D$ ranges over (representatives
of the) isomorphism
classes of $\mathbb{C}^I$-modules.
The diagram 
\eqref{TensorProductDiagram}
can be (conjecturally) 
used to equip the category
of sheaves on $Z$ perverse
with respect to the stratification
$Z_D = \bigsqcup Z_{D,V}$ with a
structure of a Tannakian category.
The results of this paper concerning
crystal tensor product provide 
a step toward finding a relation
between this category and
the category of integrable
finite dimensional representations
of the Lie algebra $\mathfrak{g}'$.

\end{document}